\documentclass[11pt,article,reqno,oneside]{amsart}
\usepackage[utf8]{inputenc}
\usepackage{geometry}[margins=0.5in]
\newtheorem{remark}{Remark}
\usepackage{fullpage}
\usepackage{cite}
\usepackage{amsmath,amsthm,amssymb}
\usepackage{graphicx}
\usepackage{mathabx}
\usepackage{color}
\usepackage{ifpdf}
\usepackage{amsaddr}
\usepackage{lipsum}
\usepackage{amsfonts}
\usepackage{graphicx}
\usepackage{epstopdf}
\usepackage{algorithmic}
\usepackage{tikz}
\usepackage{caption}
\usepackage{pgfplots}
\usepackage{float}
\usepackage{xr}
\usepackage{hyperref} 
\usepackage{cleveref}
\usepackage{xr-hyper} 
\usepackage{cleveref}
\ifpdf
\DeclareGraphicsExtensions{.pdf,.pdf,.png,.jpg}
\else
  \DeclareGraphicsExtensions{.pdf}
\fi
\usepackage{amsopn}

\crefname{figure}{Fig.}{Figs.}
\Crefname{figure}{Fig.}{Figs.}

\usepackage{enumitem}
\setlist[enumerate]{leftmargin=.5in}
\setlist[itemize]{leftmargin=.5in}
\usepackage{subcaption}
\title{Analysis of dry friction dynamics in a vibro-impact energy harvester}
\thanks{This is a preprint. This work was funded by the NSF CMMI/DMS 2009270 and EPSRC EP/V034391/1.}
\author{Christina Athanasouli$^{1,2,*}$, Daniil Yurchenko$^3$, Rachel Kuske$^2$}
\thanks{$*$Corresponding Author}
\thanks{$^1$Department of Mathematics \& Statistics, Williams College, Williamstown,  Massachusetts, USA}
\thanks{$^2$School of Mathematics, Georgia Institute of Technology, Atlanta, Georgia,  USA}
\thanks{$^3$Institute of Sound and Vibration Research, University of Southampton, Southampton, UK}
\email{CA: ca11@williams.edu, DY: d.yurchenko@soton.ac.uk, RK: rachel@math.gatech.edu}

\numberwithin{equation}{section}
\definecolor{darkgreen}{rgb}{0, 0.67, 0.18}
\usepackage{tasks}
\begin{document}
\begin{abstract}
Vibro-impact (VI) systems provide a promising nonlinear mechanism for energy harvesting (EH) in many engineering applications. Here, we consider a VI-EH system that consists of an inclined cylindrical capsule that is externally forced and a bullet that is allowed to move inside the capsule, and analyze its dynamics under the presence of dry friction. Dry friction introduces a switching manifold corresponding to zero relative velocity where the bullet sticks to the capsule, appearing as sliding in the model. We identify analytical conditions for the occurrence of non-stick and sliding motions, and construct a series of nonlinear maps that capture model solutions and their dynamics on the switching and impacting manifolds. An interplay of smooth (period-doubling) and non-smooth (grazing) bifurcations characterizes the transition from periodic solutions with alternating impacts to solutions with an additional impact on one end of the capsule per period. This transition is preceded by a sequence of grazing-sliding, switching-sliding and crossing-sliding bifurcations on the switching manifold that may reverse period doubling bifurcations for larger values of the dry friction coefficient. In general, a larger dry friction coefficient also results in larger sliding intervals, lower impact velocities yielding lower average energy outputs, and a shift in the location of some bifurcations. Surprisingly, we identify parameter regimes in which higher dry friction maintains higher energy output levels, as it shifts the location of grazing bifurcations.
\end{abstract}
\maketitle
\noindent\textbf{Keywords.} {periodic solutions, dry friction, vibro-impact systems, grazing bifurcations, sliding bifurcations, energy harvesting}\\
\textbf{AMS codes.} {37G15, 74H60, 74H45,  74M20, 70K50, 34A36}


\section{Introduction}

Dry friction is the force between two contacting solid bodies and is prevalent in mechanical systems leading to  energy dissipation. There are different friction laws used for modeling behavior of dynamical systems \cite{Awr}. The most common model of dry friction is based on Coulomb's law and assumes that friction is proportional to the normal force and acts in the opposite to the velocity direction. The constant of proportionality, $\mu_k$, (related to kinetic friction) depends on material properties of the bodies in contact and the roughness of the contact surfaces. Systems with dry friction demonstrate a range of interesting phenomena, including self-sustained oscillations, Painlevé paradox, stick-slip and chaotic behavior, jamming and various bifurcations  \cite{Or12,Painleve,Hogan,Stelter1991}.  

In contrast to phenomena related to friction, the type of non-smoothness in vibro-impact (VI) systems is instead related to impacts between colliding bodies. In this work we study the dynamics of a VI system as proposed in \cite{YURCHENKO2017456} for energy harvesting (EH). That is, we consider an externally driven inclined mass (capsule), with a smaller freely moving mass (bullet) in its interior. The bullet impacts the ends of the capsule, which are covered with membranes made of an electrically active
polymer - dielectric elastomer (DE). Two electrodes are located on either side of the DE providing a mechanism for a capacitor with variable capacitance. At an impact the membrane stretches or deforms, and therefore, the capacitance changes. This change leads to a charge differential which in turn can be harvested. We refer to this model as the VI-EH system.

To better capture realistic experimental settings, we consider the influence of dry friction between the bullet and the capsule on the dynamics of the system, as well as on the harvested energy. The VI-EH system exhibits non-smooth phenomena due to the impacts occurring at the ends of the capsule. In particular, Newton's law of restitution is applied at each impact, such that the relative velocity of the system changes direction and is scaled by a factor, $0<r<1$, called the restitution coefficient. Previous studies of this system in the absence of dry friction \cite{SERDUKOVA2021115811,Serdukova2019StabilityAB,Serdukova2022FundamentalCO} have given insight into the types of periodic motions that can be produced in different parameter regimes, as well as the types of bifurcations that result in transitions between the various motions. 
 Specifically, Serdukova et al. obtain analytical expressions for $T$-periodic solutions of the VI-EH system, where $T$ corresponds to one period, and derive nonlinear maps that capture the the impacts associated with these solutions \cite{Serdukova2019StabilityAB}. Linear stability analysis about fixed points of these maps reveals an interplay between smooth and non-smooth bifurcations, such as grazing at the impact surfaces, that can yield periodic solutions with incrementally higher numbers of impacts. Note that the original version of the VI-EH model incorporates a capsule and a ball. Here, we have modified the system to instead include a bullet. Recent experiments suggest this configuration design to avoid rotation at impact \cite{PrivComm_HSDY}.

Dry friction introduces an additional type of non-smoothness in the VI-EH model, as the direction of the force changes depending on the sign of the relative velocity. Our VI-EH system can be rewritten as a Filippov system between impacts, with a switching boundary corresponding to zero relative velocity. Periodic orbits that reach the switching boundary may exhibit sliding as the bullet and the capsule are moving together. 

Several studies have investigated the effect of dry friction on various VI systems \cite{CONE1995659,LI2021104001,LIU201530,ZhangFu2017} with a combination of different types of non-smooth behavior, i.e. both switching and impacting. Cone et al. \cite{CONE1995659} performed a numerical investigation of the dynamics and bifurcations of a horizontal impact pair oscillator with dry friction through a clamping force. For two specific values of the friction coefficient, they identified sequences of smooth and grazing bifurcations at the switching and impact surfaces as the amplitude or frequency of the excitation force were varied. In \cite{LUO20081030} a theory for investigating the dynamics of discontinuous systems based on the switchability of the flow was developed. Zhang and Fu \cite{ZhangFu2017} utilized this theory to find sufficient and necessary conditions for the occurrence of various types of motion in a horizontal vibro-impact pair with dry friction, illustrated with appropriate mappings. Similar techniques were employed in \cite{LI201981} for a two-degree-of-freedom system to obtain corresponding conditions and mappings that capture model solutions with different types of motion. Liu et al. \cite{LIU201530} considered a VI-capsule system that moves on a supporting surface and numerically investigated the effect of four models of friction between the capsule and the surface. These models reflect  various degrees of lubricated vs. dry contact. They found that periodic motion under thicker lubrication exists over a larger range of excitation amplitudes compared to the cases of dry and thinly lubricated contacts. Moreover, in the former case which is associated with higher friction, they observed larger velocities of the capsule, suggesting that higher friction may promote system stability. Li et al. \cite{LI2021104001} considered the influence of gravity and dry friction in an inclined impact pair on targeted energy transfer (TET). They developed and validated a numerical scheme based on the Moreau-Jean time-stepping method to simulate the model. They defined various power output measures to characterize the TET output under various combinations of excitation amplitude and inclination, mass ratio, and friction. According to their numerical results, they argued that the inclination of the system and the friction should be small to maintain high performance. To our knowledge, previous studies do not provide a systematic analysis of  bifurcation sequences and study of the stability of different phenomena in VI systems for energy harvesting, together with the influence of friction on these dynamics.

 Here, we develop a framework for analyzing the dynamics of our modified VI-EH model in the presence of dry friction. In particular, we derive nonlinear maps that combine both switching and impacting types of non-smooth behavior, and apply them to study the bifurcation structure of periodic orbits. We find that changes in the system behavior are described by sequences of smooth and non-smooth bifurcations, such as grazing-sliding, crossing-sliding and switching-sliding bifurcations due to the dynamics close to the switching boundary, and grazing bifurcations on the impact surfaces. The derivation of these maps is more complex in the presence of dry friction as solutions may be influenced by sliding dynamics on the switching boundary, and thus, the associated impact times and velocities of the periodic orbits are dependent on intermediate events at which sliding may occur. The maps must be integrated with conditions for the occurrence of sliding and non-stick motions similar to \cite{ZhangFu2017} to properly capture model solutions and transitions between different types of solutions. We also embed these conditions to eliminate unphysical solutions that contradict the switching or impacting dynamics. Including friction in the VI-EH system increases the number of non-smooth events to track, nevertheless we demonstrate how to fold intermediate times related to switching and sliding efficiently into the linear stability framework of \cite{Serdukova2019StabilityAB}. Then the linear stability analysis of the periodic solutions is comparable to studies without friction, and both stable and unstable solutions can be identified along with the types of bifurcations that lead to their onset or quenching.
 Using our semi-analytical approach, we note the following key application-oriented findings:
 \begin{itemize}
     \item  
 For most of the parameter regime under investigation the addition of dry friction reduces the velocities at impacts, which decrease energy output. In fact, we observe that higher values of the friction coefficient yield lower energy output.
 \item Higher friction can regularize chaotic solutions or limit period doubling observed under  under low friction, thus avoiding these more complex behaviors that are harder to control.
 \item We identify intervals in our parameter regime in which, surprisingly, higher dry friction is beneficial, as it shifts parameter values for the grazing bifurcation, sustaining higher energy outputs that can drop through grazing with and without friction.
 \end{itemize}
While the derivation of the nonlinear maps in this setting requires some additional analysis to handle multiple types of non-smooth bifurcations, we note that it has several advantages over purely numerical approaches: 
\begin{itemize}
    \item Verification that numerical simulations handle events at the switching boundaries accurately (see discussion in \cref{sec:ComparisonAndStabAnalysis}). 
     \item Accelerated alternative computations in sensitive parameter regimes where continuation can be slow (see discussion in \cref{sec:ComparisonAndStabAnalysis}).
    \item Knowledge of unstable branches can be useful as they offer bounds for transitions to different states (see \cref{sec:ComparisonAndStabAnalysis} for examples) and can (intermittently) emerge in stochastic settings, which we pursue in future studies.
    \item Tracking of complete sequences of non-smooth and smooth bifurcations, such as grazing-sliding, switching-sliding and crossing-sliding bifurcations within period-doubling regimes, which have not been explicitly identified in other work. These sequences often result in the quenching of period-doubled solutions via sliding, which has also been observed numerically in other studies of mechanical models of dry friction discussed above.
\end{itemize}
 The paper is organized as follows: in \cref{sec:VIEHModel} we introduce the VI-EH system, possible types of motion, and decompose the phase space into significant regions for handling the Filippov and impacting dynamics; in \cref{sec:FrameworkSection} we present our framework for characterizing periodic orbits of the system and introduce the collection of nonlinear maps. In \cref{sec:ExamplesOfPeriodicMotions} we derive analytical expressions that capture various $T$-periodic solutions, where $T$ is the period of the forcing. In \cref{sec:StabilityAnalysis} we perform a linear stability analysis of these periodic solutions, and in \cref{sec:ComparisonAndStabAnalysis} we compare analytical and numerical solutions for a range of parameter values. Moreover, in \cref{sec:Energy} we provide numerical and analytical results for the influence of dry friction on the output voltage. Finally, \cref{sec:Conclusions} contains a discussion of our results.

\section{The VI-EH model}\label{sec:VIEHModel}
\subsection{Model description}

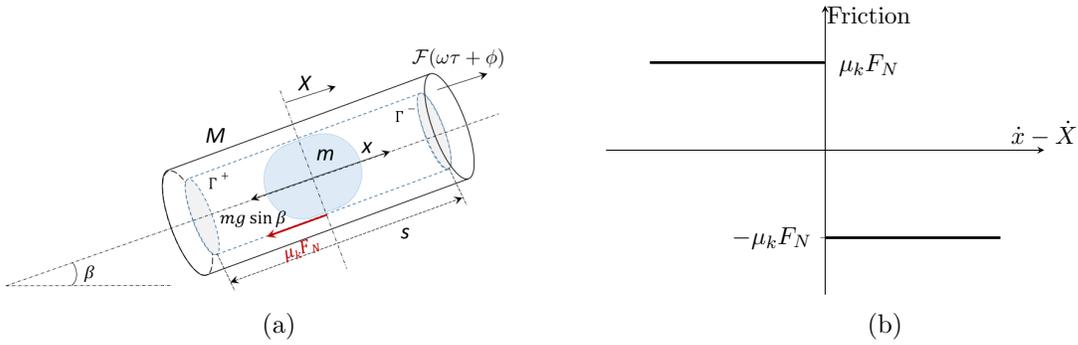
\begin{figure}[htbp!]
\centering
\begin{subfigure}[b]{0.45\textwidth}
     \includegraphics[scale=0.38]{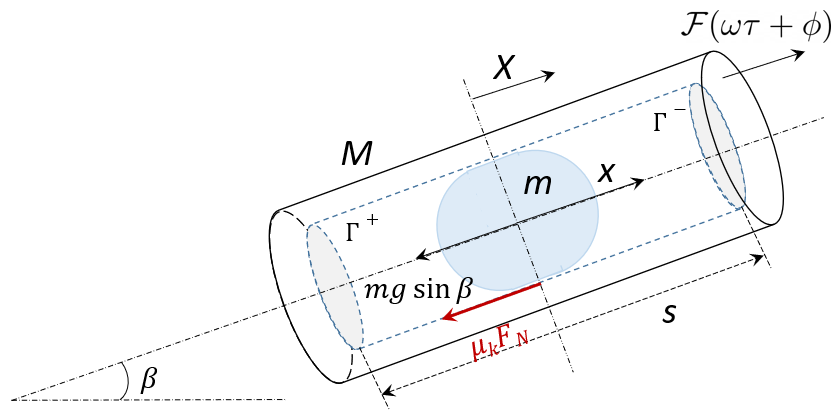}
     \caption{}
\end{subfigure}\hspace{.5cm}
\begin{subfigure}[b]{0.45\textwidth}
   \begin{tikzpicture}[scale=0.85]
  \begin{axis}[
      axis lines=middle,
      xmin=-5, xmax=5,
      ymin=-3.3, ymax=3.3,
      xtick={0},
      xticklabels={},
      ytick={-2},
      yticklabels={$-\mu_k F_N$},
      axis equal image,
      domain=-5:5,
      samples=200,
      smooth,
      clip=false,
      ]  
\addplot[mark=none,very thick,black] coordinates {(-4,2) (0,2)};
\addplot[mark=none,very thick,black] coordinates {(0,-2) (4,-2)};
 \node[above,black] at (1000,320) {$\dot{x}-\dot{X}$};
    \node[above,black] at (600,480)  {$\mu_k F_N$};
     \node[above,black] at (600,600)  {Friction};
  \end{axis}
\end{tikzpicture} 
    \caption{}
 \label{fig:FrictionGraph}   
\end{subfigure}
\caption{(a) Schematic of the VI-EH model with $F_N=mg\cos{\beta}\text{sgn}(\dot{x}-\dot{X})$ (in red) indicating the normal force. (b) Graph of the discontinuous force due to dry friction that is proportional to the normal force $F_N$.}
 \label{fig:Schematic}
\end{figure}

We consider the model of the VI-EH energy harvester based on an impact pair, consisting of a bullet located in the interior of a capsule. The bullet has a cylindrical shape with rounded ends. The motion of the capsule is driven by a harmonic excitation force $\mathcal{F}(\omega\tau +\phi)$ with period $2\pi/\omega$. For concreteness, we consider harmonic forcing  $\mathcal{F}(\omega\tau +\phi)=A\cos{(\omega\tau +\phi)}$ (\cref{fig:Schematic}).
    The motion of the bullet is driven by gravity and dry friction between impacts with the friction coefficient $\mu_k$, and its velocity changes according to the impact condition \eqref{eq:ImpactCond2} when the bullet collides with either end of the capsule. The absolute displacement of the center of the capsule and the bullet are $X(\tau)$ and $x(\tau)$, respectively. A schematic of the friction force is given in \cref{fig:FrictionGraph}. 

Let $m$ and $M$ be the mass of the bullet and the capsule, respectively. Moreover, $\beta$ and $s$ correspond to the angle of inclination and length of the capsule. By Newton's second law we have the equation of motion for $X(\tau)$ and $x(\tau)$ between impacts:
\begin{align}\label{eq:CapsuleDifEq}
\ddot{X} &= \frac{\mathcal{F}(\omega\tau +\phi)}{M}\\\label{eq:BallDifEq}
\ddot{x} &= -g\sin{\beta}-\mu_k g\cos{\beta}\text{sgn}{(\dot{x}-\dot{X})}\,,
\end{align}
where $g=9.8$ kg/s$^2$ is the acceleration due to gravity. Here, we assume that the mass of the bullet is negligible, namely $m\ll M$, and therefore, the force due to friction does not appear in \eqref{eq:CapsuleDifEq}. Specifically, $m=0.0035$ kg, and $M=0.1245$ kg, so the mass of the ball is neglected relative to the mass of the capsule. We also define the mass ratio $\eta = m/M$, and note that $\eta \approx 0.0281$. The system can undergo impacts when the bullet hits the bottom or top membrane at the ends of the capsule. The absolute velocities of the bullet and the capsule at the $k^{\textrm{th}}$ impact are denoted as $\dot{X}_{k}^{\pm}$ and $\dot{x}_{k}^{\pm}$, where the superscripts $+$ and $-$ correspond to velocities before and after the impact, respectively. We define the restitution coefficient to be:
\begin{equation}\label{eq:RestitutionCoefficient}
r = -\dfrac{\dot{X}_k^{+}-\dot{x}_k^{+}}{\dot{X}_k^{-}-\dot{x}_k^{-}}\,.
\end{equation}
Then, at the $k^{\textrm{th}}$ impact we have: $x_k - X_k = \pm \dfrac{s}{2}$, where $s$ is the length of the capsule. By conservation of momentum:
\begin{equation}\label{eq:ConservationOfMomentum}
      m\dot{x}_k^{-} + M\dot{X}_k^{-} = m\dot{x}_k^{+} + M\dot{X}_k^{+}\,.
\end{equation}
We use \eqref{eq:RestitutionCoefficient} to eliminate $\dot{X}_k^{+}$ in  \eqref{eq:ConservationOfMomentum} and get:
\begin{equation}\label{eq:ImpactCond1}
\dot{x}_k^{+} = \dfrac{\eta-r}{1+\eta}\dot{x}_k^{-} +  \dfrac{1+r}{1+\eta}\dot{X}_k^{-}\,.
\end{equation}  
 Since the mass ratio $\eta\ll 1$ the impact condition in \eqref{eq:ImpactCond1} is simplified to:
\begin{equation}\label{eq:ImpactCond2}
\dot{x}_k^{+} = -r\dot{x}_k^{-} + (1+r)\dot{X}_k^{-}\,.
\end{equation}  
\subsection{Non-dimensionalized equations and relative frame}
We reduce
the number of parameters to some key dimensionless quantities by non-dimensionalizing the original system. Furthermore, we introduce a
relative displacement $Z(t)$ in terms of the difference of the non-dimensional displacements of the capsule
$X^{*}(t)$ and of the bullet $x^{*}(t)$. We define the following non-dimensional quantities to obtain the equation of motion \eqref{eq:ZDifEq}:
\begin{align}\label{eq:NonDimQuantities}
    \tau = \dfrac{\pi}{\omega}t \quad X(\tau) = \dfrac{A\pi^2}{M\omega^2}X^{*}(t) \quad x(\tau) = \dfrac{A\pi^2}{M\omega^2}x^{*}(t) \quad A = \|\mathcal{F}\|\,\\
\label{eq:ZDifEq}
    \ddot{Z} = \frac{\mathcal{F}(\pi t + \phi)}{A} + \frac{Mg\sin{\beta}}{A}- \frac{M\mu_{k}g\cos{\beta}\,\text{sgn}(\dot{Z})}{A}
\end{align}
Here, we consider harmonic forcing for $\mathcal{F}(\pi t+\phi)$ in the non-dimensional setting
\begin{align}\label{eq:NonDimForcing}
    f(t) &= \frac{\mathcal{F}(\pi t+\phi)}{A} =  \cos{(\pi t+\phi)}\,, \text{ with }\\
\label{eq:NonDimAntiderivatives}
    F_1(t) &= \int f(t) \,dt =\dfrac{\sin(\pi t+\phi)}{\pi}\,, \qquad F_2(t) =\int F_1(t) \,dt = -\dfrac{\cos(\pi t+\phi)}{\pi^2}\,.
\end{align}
Moreover, let 
\begin{align}
    \bar{g}_1 = \dfrac{Mg\sin{\beta}}{A}\,, \qquad \bar{g}_2 = \dfrac{M\mu_{k}g\cos{\beta}}{A}\,, \qquad 
    L^{\pm}&= -(\bar{g}_1\mp\bar{g}_{2})\,.
\end{align} Then, the equation of motion between impacts can be rewritten as:
\begin{equation}\label{eq:PWSsystem}
    \Ddot{Z} = \begin{cases}
f(t) - L^{+}\,, & \dot{Z}>0\\
f(t) - L^{-}\,, & \dot{Z}<0\,.
\end{cases}
\end{equation}
The impact conditions in the relative frame are:
\begin{align}\label{eq:ImpactCondsRelFrame}
    \dot{Z}_k^{+} = -r \dot{Z}_k^{-} \text{ at } Z_k = \pm \dfrac{d}{2}, \text{ where }  d=\dfrac{M\omega^2}{A\pi^2}s \, .
\end{align}
Note that we consider $s=0.5$ m, $\omega=5\pi$ Hz, $M=0.1245$ kg, throughout the paper, unless otherwise stated. The non-dimensional parameter $d$ \eqref{eq:ImpactCondsRelFrame} depends on several physical parameters influencing the dynamics and design of the system. We note that \eqref{eq:PWSsystem} is a piecewise-smooth dynamical system and can be rewritten in the form of first order differential equations with a jump discontinuity due to the distinct dynamics across the switching boundary $\{\dot{Z}=0\}$ induced by the term corresponding to dry friction ($-M\mu_{k}g\cos{\beta}\,\text{sgn}(\dot{Z})/A$ in \cref{eq:ZDifEq}). This is a Filippov system which may exhibit sliding solutions, namely solutions that evolve along the switching boundary for a nontrivial amount of time. Physically, this sliding that may occur away from impacts due to dry friction corresponds to the capsule and the bullet moving together and ``sticking". During these time intervals the sum of the excitation force $f(t)$ and gravity cannot overcome the static friction between the bullet and the capsule, and thus, the two bodies are moving together. 
After some time, the excitation force becomes sufficiently large so that the bullet and the capsule again move relative to each other.

As introduced in previous related work (e.g. \cite{SERDUKOVA2021115811, Serdukova2019StabilityAB}), to distinguish between orbits involving different numbers of impacts on either end of the capsule, we use the notation $n$:$m$/$pT$ to categorize periodic orbits for the VI-EH device with $T$-periodic external forcing. Here, $n$ ($m$) corresponds to the number of impacts on the bottom (top) of the capsule per time interval $T$, for a $pT$-periodic motion where $p\in\mathbb{Z}$. For $p = 1$, we simplify the notation to $n$:$m$. 

\subsection{Phase space in relative frame and types of motion}\label{sec:PhasePlane}

We define the following sets that compose the phase space $(Z,\dot{Z})$ and use them to describe various types of motion that occur in the VI-EH system:

\begin{align}
\begin{split}
\Sigma &= \{(Z,\dot{Z}) | \dot{Z}=0\, , Z\in(-d/2,d/2) \}\\
    \Sigma^{+(-)} &= \{(Z,\dot{Z}) | \dot{Z}>0\, (<0)\, , Z\in(-d/2,d/2) \}\\
       \Gamma^{\pm}  &= \{(Z,\dot{Z}) |  Z=\pm d/2 \}\\
 \end{split}
\end{align}

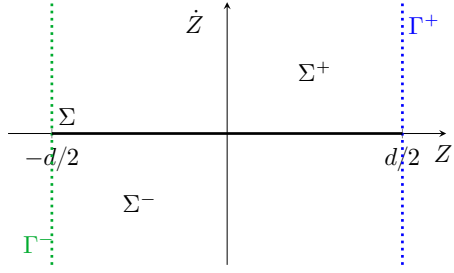
\begin{figure}[htbp]
\centering
   \begin{tikzpicture}[scale=0.85]
  \begin{axis}[
      xlabel={\hspace{0.7cm}$Z$},
      ylabel={$\dot{Z}$},
      axis lines=middle,
      xmin=-5, xmax=5,
      ymin=-3, ymax=3,
      xtick={-4,0,4},
      xticklabels={$-d/2$,0,$d/2$},
      ytick={-4,4},
      axis equal image,
       ylabel style={above}, 
    xlabel style = { right},
      domain=-5:5,
      samples=200,
      smooth,
      clip=false,
      ]  \addplot[mark=none,dotted, color=darkgreen, very thick] coordinates {(-4,-3) (-4,3)};
\addplot[mark=none,dotted,blue, very thick] coordinates {(4,-3) (4,3)};
\addplot[mark=none,very thick,black] coordinates {(-4,0) (4,0)};
 \node[above,black] at (700,400) {$\Sigma^{+}$};
    \node[above,black] at (300,100)  {$\Sigma^{-}$};
    \node[above,black] at (135,300) {$\Sigma$};
    \node[blue] at (950,550) {$\Gamma^{+}$};
    \node[color=darkgreen] at (70,50) {$\Gamma^{-}$};
  \end{axis}
\end{tikzpicture} 
 \label{fig:PhasePlaneSchematic}
    \caption{Schematic of the phase plane in relative frame. The blue and green dotted lines indicate the impact surfaces, $\Gamma^{+}$ (bottom membrane) and $\Gamma^{-}$ (top membrane), respectively. The thick black line corresponds to the switching manifold, $\Sigma$, with $\Sigma^{+}$ ($\Sigma^{-}$) indicating regions of the phase plane with $\dot{Z}>0$ $(<0)$.}
\end{figure}

The switching boundary is denoted by $\Sigma$, and $\Sigma^{+}$ and $\Sigma^{-}$ are regions of the phase space in which the relative velocity, $\dot{Z}$, is positive and negative, respectively (\cref{fig:PhasePlaneSchematic}). The impact surfaces are denoted by $\Gamma^{+}$ and $\Gamma^{-}$. 

Below we describe three main types of motion that may occur in our system and the associated equations of motion, noting four 1:1/$pT$ periodic solutions within these types.

\begin{figure}[hbpt!]
\centering
    \includegraphics[scale=0.19]{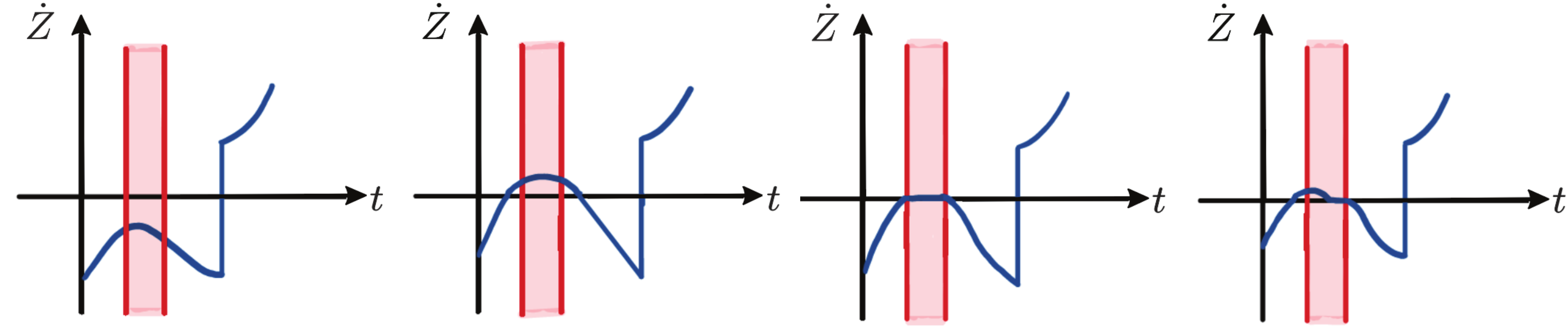}
 \caption{Schematic representations of the motion from $\Gamma^{+}$ to $\Gamma^{-}$ in (from left to right) 1:1, 1:1$^{\Sigma}_{c}$, 1:1$^{\Sigma}_{s}$, and 1:1$^{\Sigma}_{cs}$ periodic motions in the $\dot{Z}$ vs $t$ plane. The red rectangle indicates the time interval over which sliding may occur determined by \eqref{eq:SlidingOnset2a}-\eqref{eq:SlidingOnset2c} and \eqref{eq:NonPassToSemiPassCondToOmega2_1}-\eqref{eq:NonPassToSemiPassCondToOmega2_3}.}  
 \label{fig:Schematic_DifferentMotions}
\end{figure}

\begin{figure}[hbpt!]
    \centering
      \begin{subfigure}[b]{0.3\textwidth}
 \includegraphics[scale=0.38]{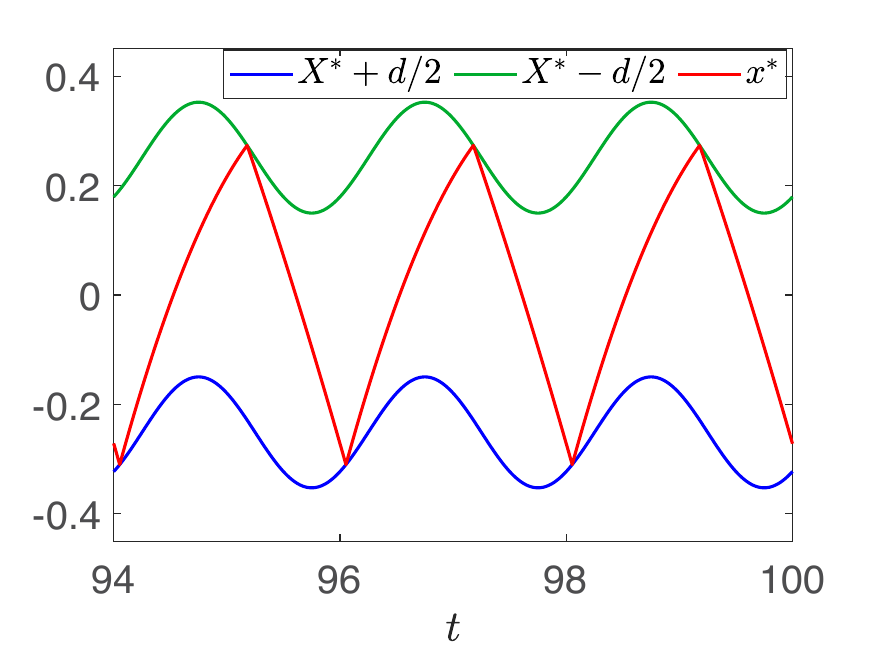}
       \caption{}        \label{fig:NonStickSimple_Timetrace_AbsFrame}
    \end{subfigure}
    \hfill
      \begin{subfigure}[b]{0.3\textwidth}
\includegraphics[scale=0.38]{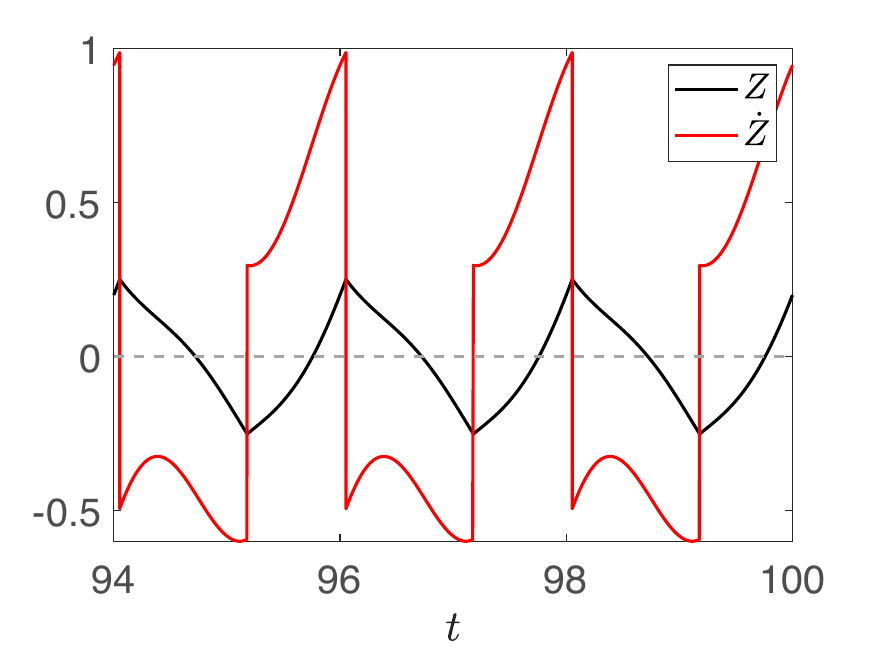}
          \caption{}
\label{fig:NonStickSimple_Timetrace}
    \end{subfigure}
    \hfill
    \begin{subfigure}[b]{0.3\textwidth}
        \includegraphics[scale=0.38]{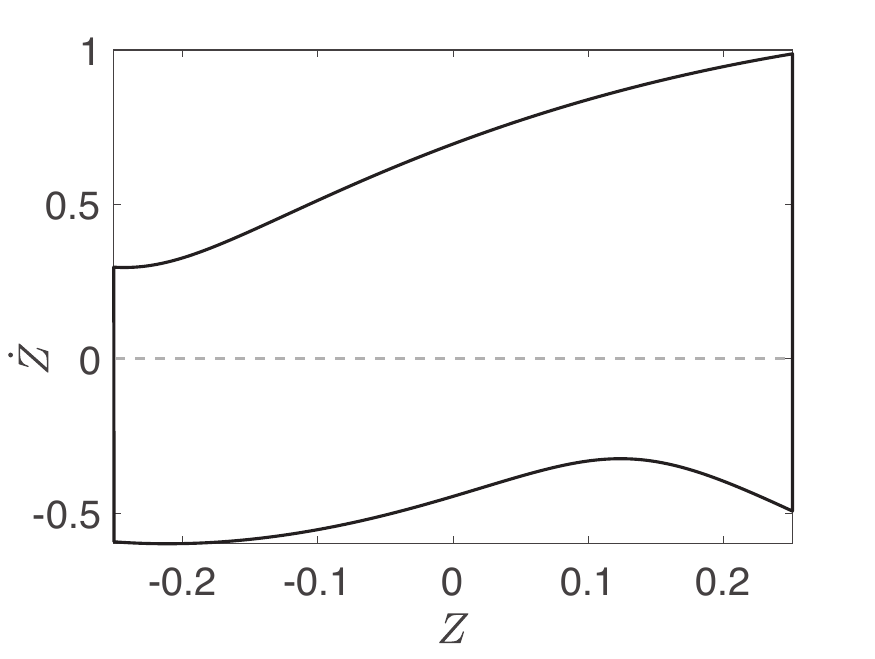}
       \caption{}
        \label{fig:NonStickSimple_PhasePlane}
    \end{subfigure}
    \hfill
    \begin{subfigure}[b]{0.3\textwidth}
        \includegraphics[scale=0.38]{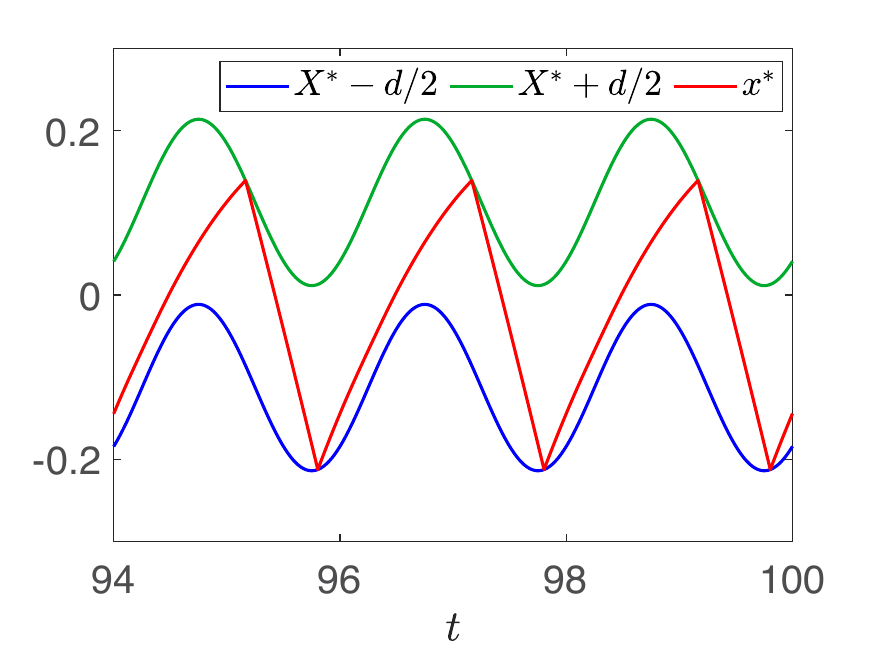}
       \caption{}
        \label{fig:NonStickLoop_Timetrace_AbsFrame}
    \end{subfigure}
    \hfill
       \begin{subfigure}[b]{0.3\textwidth}
\includegraphics[scale=0.38]{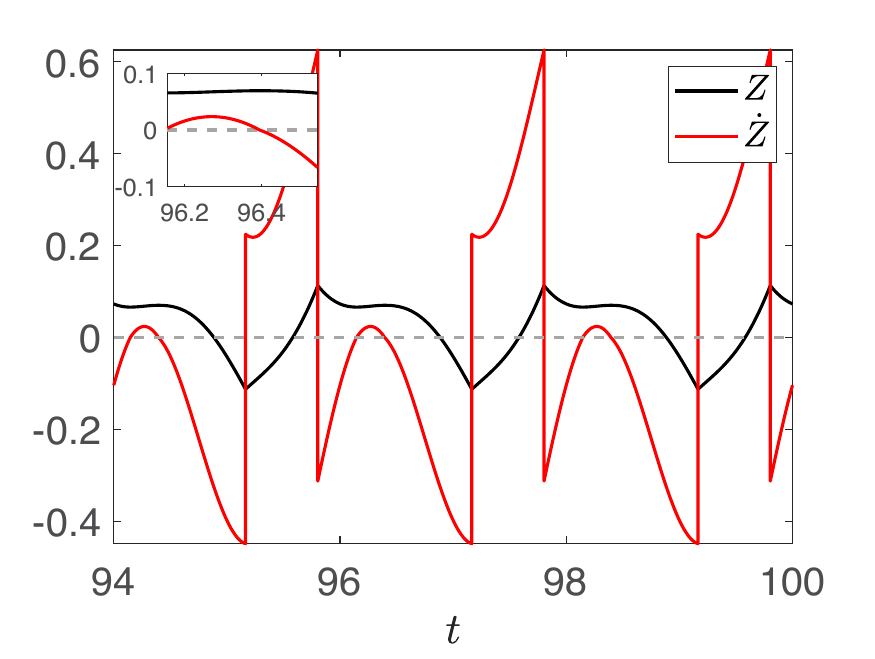}
          \caption{}
           \label{fig:NonStickLoop_Timetrace}
    \end{subfigure}
    \hfill
    \begin{subfigure}[b]{0.3\textwidth}
        \includegraphics[scale=0.38]{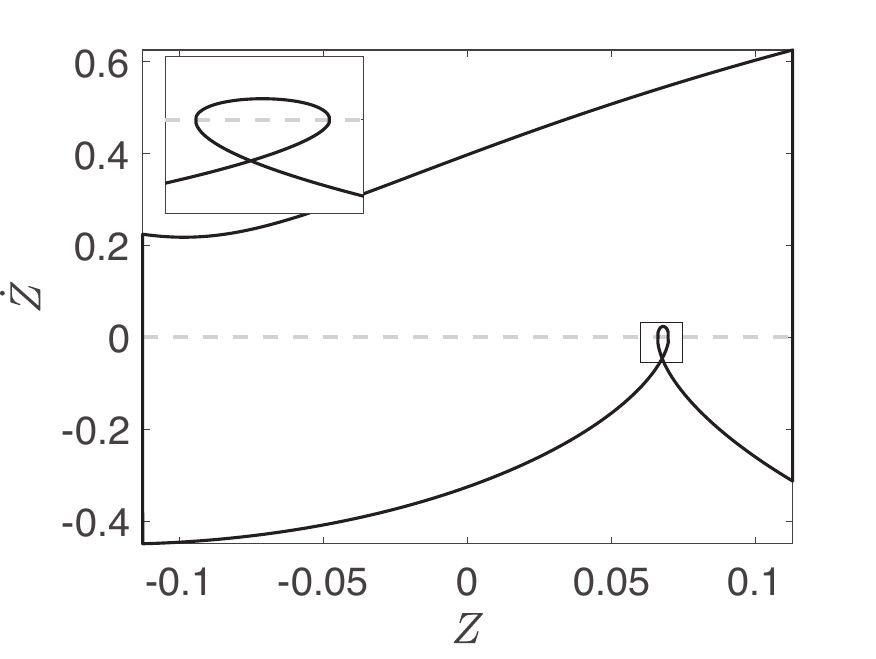}
       \caption{}
        \label{fig:NonStickLoop_PhasePlane}
    \end{subfigure}
          \begin{subfigure}[b]{0.3\textwidth}
\includegraphics[scale=0.38]{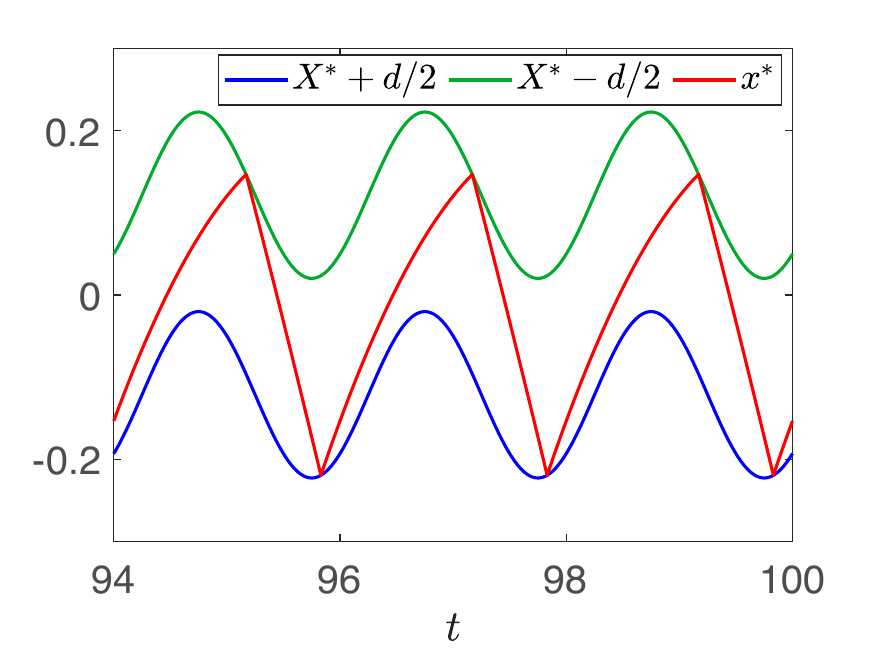}
          \caption{}
           \label{fig:SlidingStick_Timetrace_AbsFrame}
    \end{subfigure}
    \hfill
        \begin{subfigure}[b]{0.3\textwidth}
\includegraphics[scale=0.38]{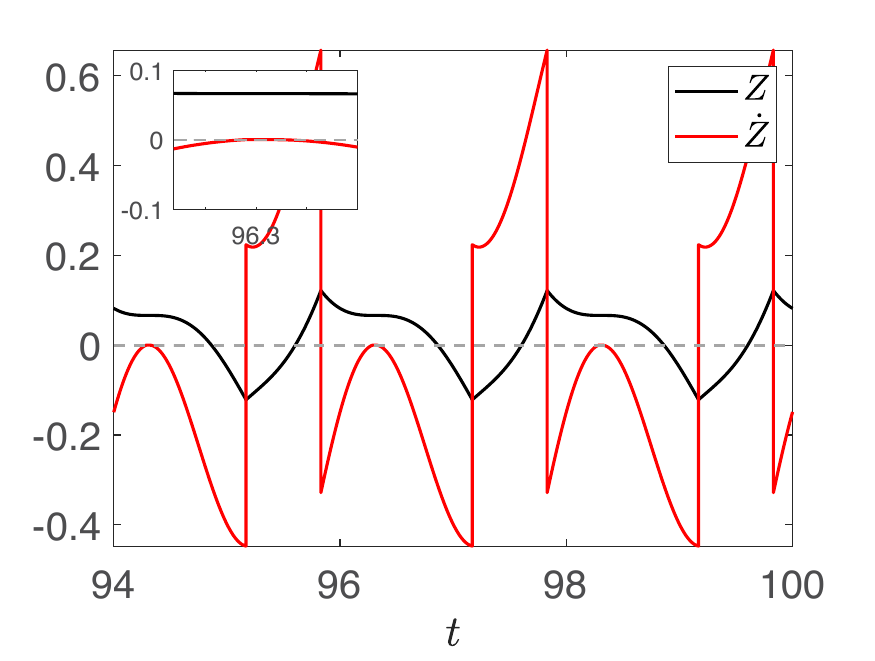}
          \caption{}
           \label{fig:SlidingStick_Timetrace}
    \end{subfigure}
    \hfill
    \begin{subfigure}[b]{0.3\textwidth}
        \includegraphics[scale=0.38]{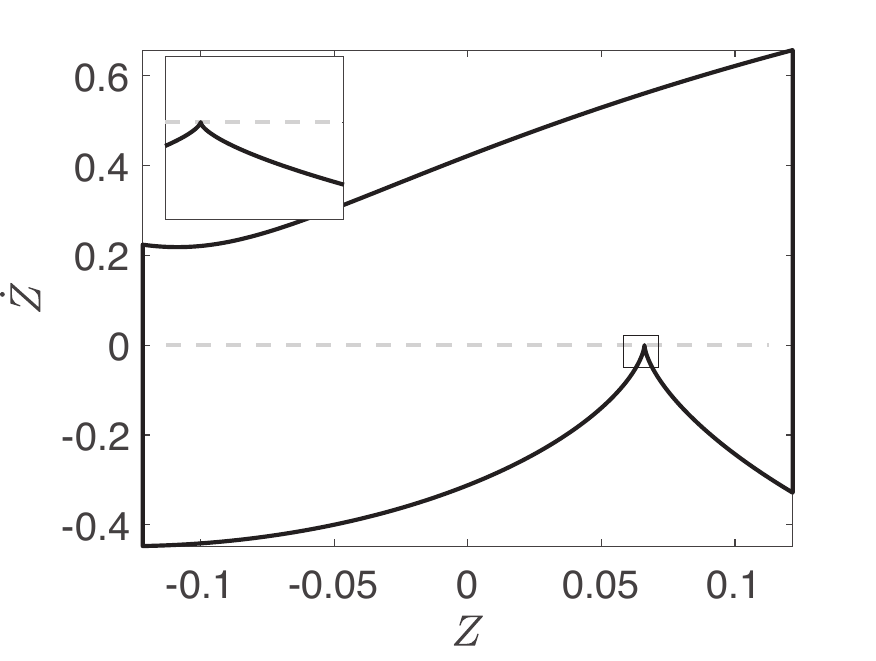}
       \caption{}
        \label{fig:SlidingStick_PhasePlane}
    \end{subfigure}
    \caption{Periodic solutions demonstrating the different types of motion in the VI-EH model with dry friction. The parameter values are $r=0.5$, $\beta=\pi/4$, and $\mu_k=0.5$. Left column: time series in absolute frame; Middle column: times series in relative frame. Right column: $Z-\dot{Z}$ phase plane. Panels (a),(b),(c): $A=3.1$, $d=0.502$ show a non-stick motion in which the relative velocity, $\dot{Z}$, does not change sign. Panels (d),(e),(f): $A=6.9$, $d=0.2255$ show a non-stick motion in which the relative velocity, $\dot{Z}$, changes sign. Panels (g),(h),(i): $A=6.4$, $d=0.2432$ show a sliding motion in which the relative velocity, $\dot{Z}=0$ for a nontrivial period of time. The gray dashed line in the middle- and right-column panels corresponds to $\dot{Z}=0$.}
    \label{fig:DifferentMotions}
\end{figure}

\begin{remark}
    For simplicity, and based on the values of $\beta$ used ($\pi/4$, $\pi/6$ and $\pi/12$), namely $\beta$ bounded away from zero, a trajectory reaches $\Sigma$ during the upward motion of the bullet (from $\Gamma^{+}$ to $\Gamma^{-}$), while the trajectory does not reach $\Sigma$ during the downward motion (from $\Gamma^{-}$ to $\Gamma^{+}$). We describe the upward motion from $\Gamma^{+}$ to $\Gamma^{-}$ which distinguishes between the four 1:1-type periodic solutions.
\end{remark}  
 Schematic representations and actual time series of \eqref{eq:PWSsystem}-\eqref{eq:ImpactCondsRelFrame} of each solution type are shown in \cref{fig:Schematic_DifferentMotions} and \cref{fig:DifferentMotions}, respectively.
\begin{enumerate}
 \item Non-stick motion in $ \Sigma^{\lambda}$, $\lambda=\pm$:  
 There are multiple types of non-stick motion, as trajectories may evolve away from or reach  $\Sigma$ between impacts. In the latter, the sum of the excitation and gravity forces is sufficient to overcome friction and cross $\Sigma$.
In either case, during non-stick motion, we have:
\begin{equation}\label{eq:NonStickMotionEqs}
    \ddot{Z} = f(t) -L^{\lambda}\,, \text{ in } \Sigma^{\lambda}\,\text{ for }\lambda=\pm\,.
\end{equation}
We compute two specific periodic solutions of this type:
\begin{itemize}
    \item 1:1 solution: This is a 1:1 periodic solution in which the discontinuity in behavior occurs only due to the impacts at $\Gamma^{\pm}$, with dynamics following \eqref{eq:NonStickMotionEqs}
    (see top left panel in \cref{fig:Schematic_DifferentMotions} and \cref{fig:NonStickSimple_Timetrace_AbsFrame}, \cref{fig:NonStickSimple_Timetrace}, \cref{fig:NonStickSimple_PhasePlane}).
    \item 1:1$^{\Sigma}_{c}$ solution: From $\Gamma^{+}$ to $\Gamma^{-}$ this solution crosses $\Sigma$ twice, as it moves from $\Sigma^{-}$ to $\Sigma^{+}$ and back to $\Sigma^{-}$ in a non-stick motion following  \eqref{eq:NonStickMotionEqs}, then continuing to $\Gamma^-$ (see top right panel in \cref{fig:Schematic_DifferentMotions} and \cref{fig:NonStickLoop_Timetrace_AbsFrame}, \cref{fig:NonStickLoop_Timetrace}, \cref{fig:NonStickLoop_PhasePlane}). 
\end{itemize}
    \item Sliding motion on $\Sigma$: When 
the sum of the excitation and gravity forces cannot overcome friction, the bullet moves together with the capsule. Then the solution slides along a fixed value of the relative displacement $Z$, with 
    \begin{equation}\label{eq:SlidingStickMotionEqs}
    \ddot{Z} = 0 \text{ and } \dot{Z} = 0 \quad (\text{motion on }\Sigma)\,.
\end{equation}
    During this motion $Z$ remains constant, as bullet and capsule move together ($\dot{x}-\dot{X}=0$). We compute two periodic solutions that include segments of sliding motion:
    \begin{itemize}
        \item {1:1$^{\Sigma}_{s}$ solution: This solution intersects $\Sigma$, and initiates sliding \eqref{eq:SlidingStickMotionEqs}. Following the sliding motion, the solution trajectory continues through $\Sigma^{-}$, and intersects the impact surface $\Gamma^{-}$ (see bottom left panel in \cref{fig:Schematic_DifferentMotions} and \cref{fig:SlidingStick_Timetrace_AbsFrame}, \cref{fig:SlidingStick_Timetrace}, \cref{fig:SlidingStick_PhasePlane}).}
     \item 1:1$^{\Sigma}_{cs}$ periodic solution: This solution combines crossing and sliding, as the trajectory reaches $\Sigma$, crosses to $\Sigma^{+}$, returns to $\Sigma$, slides according to \eqref{eq:SlidingStickMotionEqs}, then evolves
     through $\Sigma^{-}$ to reach the impact surface $\Gamma^{-}$ (see bottom right panel in \cref{fig:Schematic_DifferentMotions}).
    \end{itemize}
    Note that in the above notation we use the superscript $\Sigma$ when periodic solutions reach $\Sigma$, and the subscripts $c$ and $s$, when the solutions cross and/or slide on $\Sigma$, respectively.
    \item Side-stick motion in $ \Sigma\cap\Gamma^{\lambda}$, $\lambda=\pm$: When the bullet reaches the bottom or top membranes $\Gamma^{\pm}$ of the capsule at $\Sigma$ with $\dot{Z}=0$, the vector fields away from $\Gamma^{\pm}$ may point to $\Sigma \, \cap \, \Gamma^{\pm}$. Then the bullet and the capsule must move together. During side-stick motion, we have:
    \begin{equation}\label{eq:SideStickMotionEqs}
    \ddot{Z} = 0 \text{ and } \dot{Z} = 0 \text{ in } \Sigma\,\cap\,\Gamma^{\lambda}\,\text{ for }\lambda=\pm\,.
\end{equation}
\end{enumerate}
We leave the detailed analysis of side-stick motion for future research, as it does not appear in this study. 

Throughout this paper the term ``sliding'' corresponds to $\dot{Z}=\dot{X}-\dot{x}=0$, and indicates sliding on the state-dependent switching boundary $\Sigma$, during which the bullet and the capsule ``stick" and thus  move together. The non-stick and sliding motion are also referred to as slipping and stick-slip motions, respectively, in the literature.
\section{Framework for the analysis of periodic solutions}
\label{sec:FrameworkSection}

Numerical simulations provide insight into the number of impacts characterizing each stable periodic solution. Below we show bifurcation diagrams of the impact velocities with respect to the dimensionless length parameter $d$ \cref{eq:ImpactCondsRelFrame}. The blue (green) dots correspond to impact velocities, $\dot{Z}_{k}$, at $\Gamma^{+}$ ($\Gamma^{-}$), namely the bottom (top) membrane, that were obtained numerically. The numerical results shown throughout the paper have been obtained via a continuation type method for decreasing $d$ (increasing $A$ or decreasing $s$), unless stated otherwise for specific figures.

\begin{figure}[hbtp!]
    \centering
    \begin{subfigure}[b]{0.45\textwidth}
\includegraphics[scale=0.42]{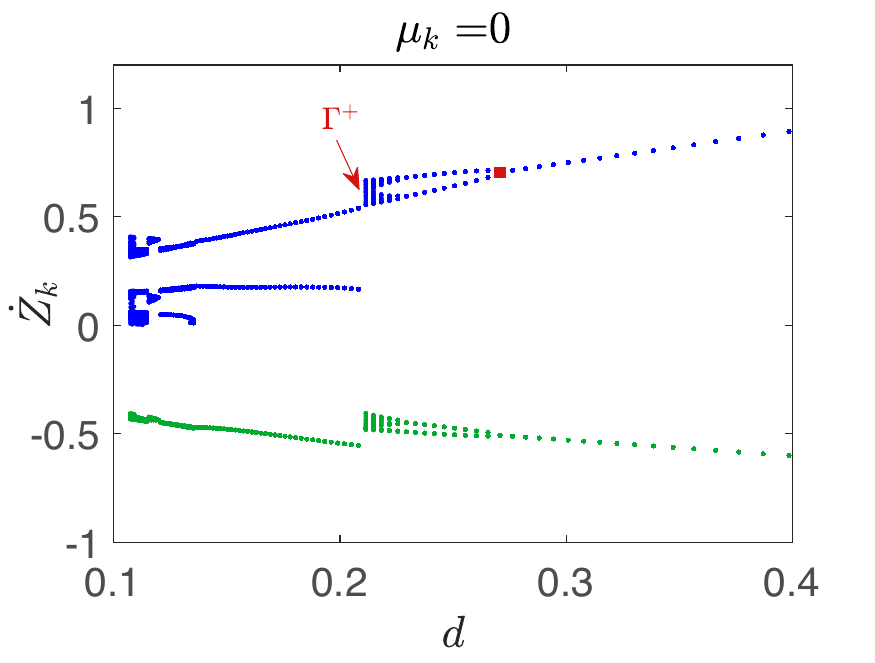}
          \caption{}
           \label{fig:ZdotBifDiagrams_FrictionComp_A}
    \end{subfigure}
    \hfill
    \begin{subfigure}[b]{0.45\textwidth}
        \includegraphics[scale=0.42]{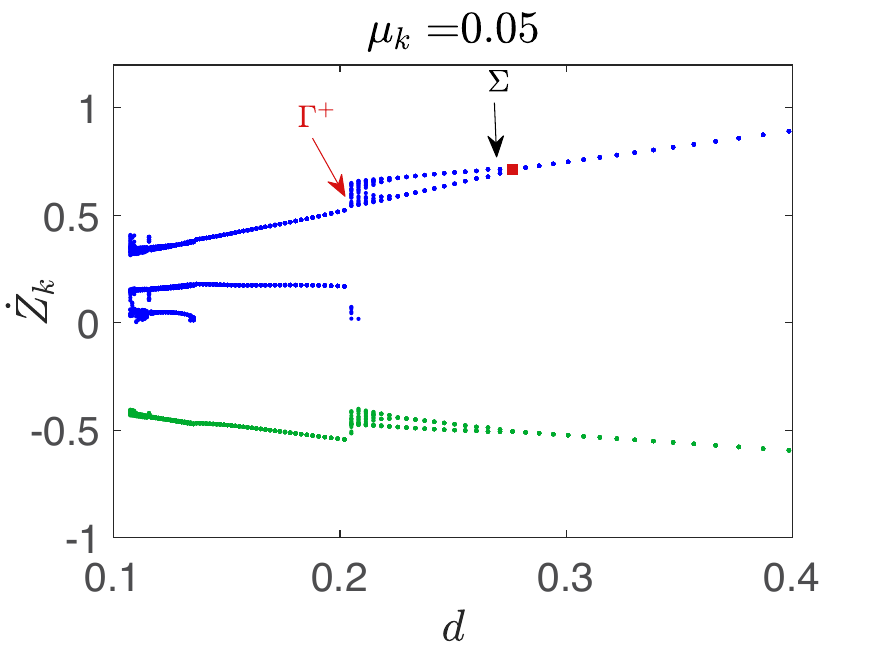}
       \caption{}
        \label{fig:ZdotBifDiagrams_FrictionComp_B}
    \end{subfigure}
    \begin{subfigure}[b]{0.45\textwidth}
\includegraphics[scale=0.42]{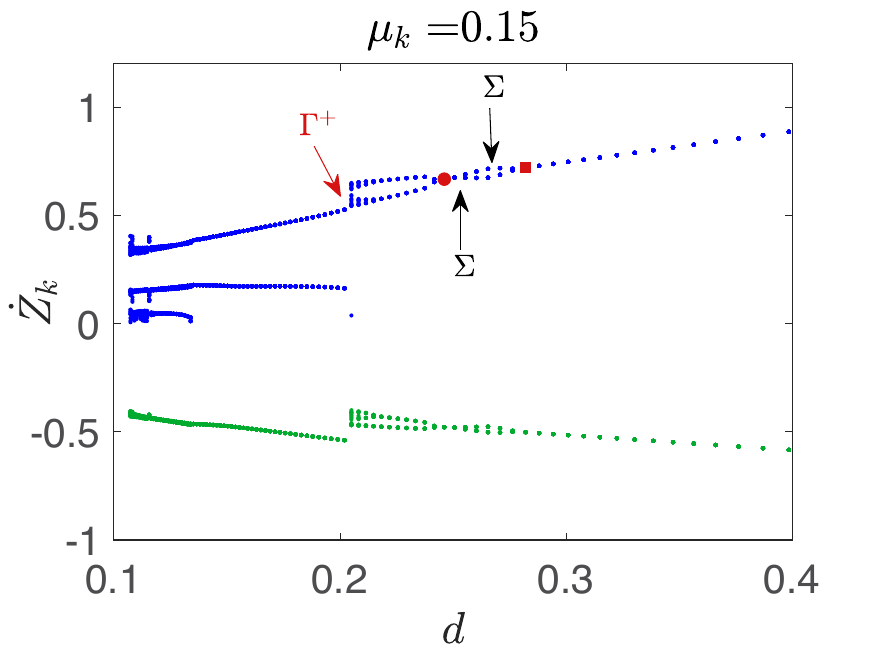}
         \caption{}
          \label{fig:ZdotBifDiagrams_FrictionComp_C}
    \end{subfigure}
    \hfill
     \begin{subfigure}[b]{0.45\textwidth}
        \includegraphics[scale=0.42]{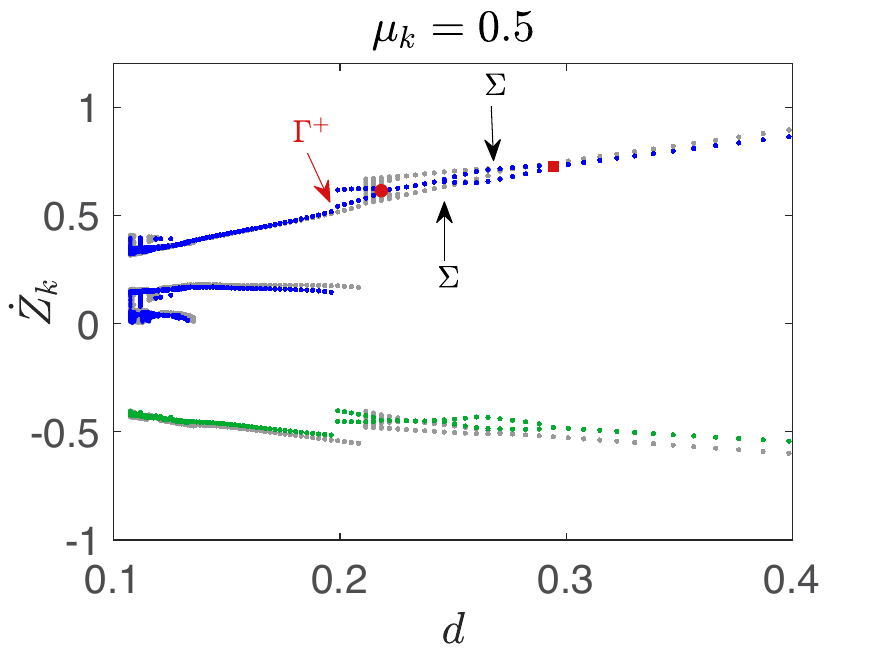}
       \caption{}
         \label{fig:ZdotBifDiagrams_FrictionComp_D}
    \end{subfigure}
   \caption{Bifurcation diagrams of impact velocities, $\dot{Z}_{k}$, vs $d$ for different values of $\mu_k$, with $r=0.5$, $\beta=\pi/4$, and $A$ increasing over $[3.906,14.4]$. The blue dots correspond to impact velocities on $\Gamma^{+}$ (bottom membrane), while the green dots correspond to impact velocities on $\Gamma^{-}$ (top membrane). In (d) the gray dots correspond to impact velocities in panel (a) with $\mu_k=0$ for comparison. The red circles and squares indicate PD bifurcations, while the black arrows point to grazing on $\Sigma$ and the red arrows indicate transitions via grazing on $\Gamma$. These occur at $d=d_{\Sigma}$ and $d=d_{\Gamma^{+}}$, respectively, discussed in \cref{sec:ExamplesOfPeriodicMotions}.}
    \label{fig:ZdotBifDiagrams_FrictionComp}
\end{figure}

In \cref{fig:ZdotBifDiagrams_FrictionComp_A}, we see that a period doubling (PD) bifurcation to a 1:1/$2T$ solution occurs at $d = 0.2707$ (red square). The period doubling cascade leads to chaotic behavior for $d\in(0.21143,0.2707)$ after which a grazing bifurcation on $\Gamma^{+}$ leads to the loss  of the 1:1 chaotic (1:1$/C$) behavior and the gain of the stable 2:1 solution (red arrow in \cref{fig:ZdotBifDiagrams_FrictionComp_A}). Such a grazing bifurcation is characterized by a vanishing impact velocity, $\dot{Z}_{k}=0$, on $\Gamma^{\pm}$. For values of $d$ below the grazing bifurcation, the additional branch in the 2:1 periodic 
has small impact velocity $\dot{Z}_{k}$.

The bifurcation diagram is qualitatively similar for $\mu_k=0.05$ (\cref{fig:ZdotBifDiagrams_FrictionComp_B}). Compared to the case without friction ($\mu_k=0$), for $\mu_k=0.05$ the period doubling cascade occurs over a larger range of $d$ values ($0.205<d<0.2762$ with the region indicated by the red square and arrow), after which the transition from the 1:1$/C$ to 2:1 regime follows from a grazing bifurcation on $\Gamma^{+}$.
Furthermore, the resulting impact velocities are smaller. For $\mu_k=0.05$, we observe a grazing-sliding bifurcation \cite{DiBernardo_Review2008,Jeffrey_Hogan2011} on $\Sigma$ occurring at $d=0.266$ (black arrow), which introduces $2T$ - solutions that slide along $\Sigma$ during one of the forcing periods $T$ (see bottom left panel \cref{fig:Schematic_DifferentMotions}). The sliding motion is replaced by crossing on $\Sigma$, as $d$ decreases (see top right panel \cref{fig:Schematic_DifferentMotions}).

For larger values of the dry friction coefficient, $\mu_k$, we observe qualitative differences in the bifurcation diagrams of $\dot{Z}_k$. For $\mu_k=0.15$ (\cref{fig:ZdotBifDiagrams_FrictionComp_C}) a PD bifurcation occurs at $d=0.2768$ (red square). This PD bifurcation is followed by a grazing-sliding bifurcation on $\Sigma$ (black arrow) at $d=0.266$ that introduces segments over which the $2T$-periodic solution slides on or crosses $\Sigma$ ($\dot{Z}=0$). However, we see that a stable period $T$ 1:1-type solution reemerges at $d=0.2506$ (black arrow) via a grazing bifurcation on $\Sigma$, which is subsequently, 
followed by another PD bifurcation at $d=0.246$ (red circle). The grazing bifurcation introduces sliding and crossing motions on $\Sigma$, similar to the cases shown in the top right and bottom panels of \cref{fig:Schematic_DifferentMotions} and \cref{fig:NonStickLoop_Timetrace_AbsFrame,fig:NonStickLoop_Timetrace,fig:NonStickLoop_PhasePlane,fig:SlidingStick_Timetrace_AbsFrame,fig:SlidingStick_Timetrace,fig:SlidingStick_PhasePlane}. The end of the 1:1/$C$ branch follows from a grazing bifurcation on $\Gamma^{+}$ (red arrow). A similar sequence of bifurcations is observed for $\mu_k=0.5$ (\cref{fig:ZdotBifDiagrams_FrictionComp_D}). A 1:1/$2T$ solution occurs for $d\in(0.246,0.2943)$ (red square indicate PD bifurcation and black arrows indicate grazing-sliding bifurcation on $\Sigma$). The first grazing-sliding bifurcation on $\Sigma$ ($d=0.2633$) introduces sliding motions and crossings at $\Sigma$, while the second one ($d=0.2435$) leads to the loss of the $2T$-periodic solution. The second occurrence of a 1:1/$2T$ solution is observed for $d\in(0.1989,0.2183)$ (red circle indicates PD bifurcation and red arrows indicate grazing bifurcation on $\Gamma^{+}$). In addition, the 1:1/$2T$ solution remains stable for this entire interval of $d$-values and solutions of higher periods of the form $2pT$ for $p>1$ are not observed. Thus, the addition of friction can limit the appearance of chaotic regimes in the VI-EH device.

To characterize the different properties of 1:1$/pT$ periodic solutions as well as the influence of the dry friction on the VI-EH system dynamics we develop a mathematical framework as described in the following sections.

\subsection{Derivation of discrete nonlinear maps}
\label{sec:AnalysisIntro} 
\begin{figure}[hbtp!]
    \centering
    \includegraphics[scale=0.41]{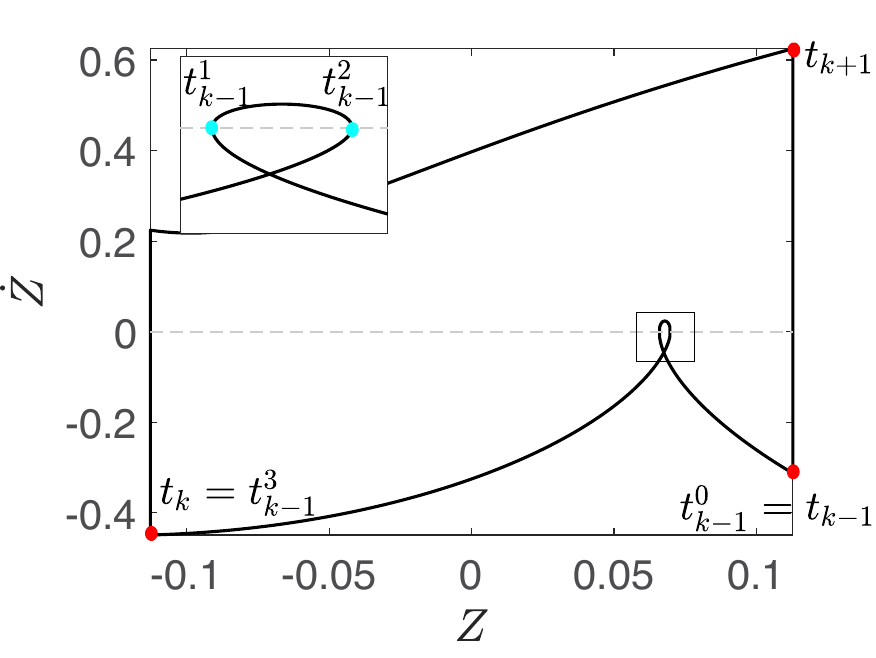}
    \caption{Example of a 1:1$^{\Sigma}_{c}$ periodic solution, crossing $\Sigma$ (gray dashed line) while evolving from $\Gamma^{+}$ to  $\Gamma^{-}$, with no crossing or sliding while evolving from $\Gamma^{-}$ to $\Gamma^{+}$. \Cref{fig:NonStickLoop_PhasePlane} is used here to illustrate the notation in \eqref{eq:EqMotionVel}-\eqref{eq:EqMotionPosition}. The impacts (red circles) occur on $\Gamma^{+}$ at $t=t_{k-1}$ and $t=t_{k+1}$, and on $\Gamma^{-}$ at $t=t_{k}$. The inset focuses on the loop that occurs due to the two crossings (cyan circles) at $\Sigma$ at $t=t^{1}_{k-1}$ and $t=t^{2}_{k-1}$.}
    \label{fig:NotationExplanation}
\end{figure}

We derive discrete maps that capture the dynamics of these solutions between impacts by integrating the equations of motion \cref{eq:ZDifEq} (or \cref{eq:PWSsystem}) between impacts that occur at times $t_{k-1}$ and $t_{k}$ (in \eqref{eq:EqMotionVel},\eqref{eq:EqMotionPosition}), applying the impact conditions \cref{eq:ImpactCondsRelFrame}, and tracking intermediate events on $\Sigma$, to get the the impact times, $t_{k}$, and velocities, $\dot{Z}_{k}$. This integration yields:
\begin{align}
\label{eq:EqMotionVel}
   \dot{Z}_{k} & = -r\dot{Z}_{k-1} + \sum_{j=1}^{N} S(j)\bigg[ \int_{t^{j-1}_{k-1}}^{t^{j}_{k-1}} f(a)\, da  -  L^{\lambda(j)}\Delta t^{j-1}_{k-1}\bigg]\\
   \begin{split}\label{eq:EqMotionPosition}
          Z_k & = Z_{k-1} -r\dot{Z}_{k-1}\Delta t^{0}_{k-1} + \sum_{j=1}^{N} S(j)\bigg[\int_{t^{j-1}_{k-1}}^{t^{j}_{k-1}} F_1(a)\, da   \\
     & + F_1(t^{j-1}_{k-1}) \Delta t^{j-1}_{k-1} -\frac{L^{\lambda(j)}}{2}(\Delta t^{j-1}_{k-1})^2 \bigg]\,,  
   \end{split}
\end{align}
where 
\begin{itemize}
    \item  $k$: index denoting an impact event on $Z= \Gamma$, such that $Z_k=\pm d/2$ (red points in \cref{fig:NotationExplanation}). Note that we have dropped the superscript $-$, so that $\dot{Z}_{k}=\dot{Z}^{-}_{k}$, since $\dot{Z}^{\pm}_{k}$ are continuous between impacts and $\dot{Z}^{+}_{k}$ can be expressed in terms of $\dot{Z}^{-}_{k}$ using the impact condition \eqref{eq:ImpactCondsRelFrame}.
    \item $j$: index denoting intermediate events on $\Sigma$ between impacts occurring at times $t^j_{k-1}$ (cyan points in \cref{fig:NotationExplanation}).
    \item $N$: number of intermediate subintervals between impact times. Note that $t^{0}_{k-1} =t_{k-1}$ and $t^{N}_{k-1}  =t_{k}$. If a solution does not reach $\Sigma$ between impacts, then $N=1$ and $t^{1}_{k-1} = t_{k}$.
    \item $\Delta t^{j-1}_{k-1}$: length of the time interval $[t^{j-1}_{k-1},t^{j}_{k-1}]$, namely $\Delta t^{j-1}_{k-1}=t^{j}_{k-1}-t^{j-1}_{k-1}$, for $j=1,\dots, N$. If  $N=1$, then $\Delta t^{0}_{k-1} = \Delta t_{k-1}=t_{k}-t_{k-1}$.
    \item $S(j)$: indicator variable that flags sliding on $\Sigma$: If sliding occurs on the interval $[t^{j-1}_{k-1},t^{j}_{k-1}]$, then $S(j)=0$, otherwise $S(j)=1$. In both cases the trajectory reaches $\Sigma$ at $t=t^{j-1}_{k-1}$. In the sliding case ($S(j)=0$), $t^{j-1}_{k-1}$ satisfies \eqref{eq:SlidingOnset2a}-\eqref{eq:SlidingOnset2c} and $t^{j}_{k-1}$ is determined by \eqref{eq:NonPassToSemiPassCondToOmega2_1}-\eqref{eq:NonPassToSemiPassCondToOmega2_3}. Otherwise ($S(j)=1$), $t^{j-1}_{k-1}$ satisfies \eqref{eq:CrossingConditionsMinusToPlus_1}-\eqref{eq:CrossingConditionsMinusToPlus_2}. Note that for $j=1$, $S(j) \neq 0$, since we do not consider the case of side-stick motion here. If $S(j) =0$ on an interval $[t^{j-1}_{k-1},t^{j}_{k}]$, for $j>1$, $S(j) =0$ ensures that the relative velocity, $\dot{Z}_k$, is zero and the relative displacement, $Z_k$, is constant on that interval. The absolute velocity and displacement continue to vary with $t$ according to \eqref{eq:CapsuleDifEq}.
    \item $\lambda(j)=+$ or $-$, if the trajectory evolves in $\Sigma^{+}$ or $\Sigma^{-}$, respectively, on the interval $[t^{j-1}_{k-1},t^{j}_{k-1}]$.
\end{itemize}
It is also useful to define the relative velocities, $\dot{Z}^{j}_{k-1}$, and positions, ${Z}^{j}_{k-1}$, at times $t^{j}_{k-1}$, respectively:
\begin{align}
\label{eq:IntermediateVelocity}
   \dot{Z}^{j}_{k-1} & = -r\dot{Z}_{k-1} + \sum_{i=1}^{j} S(i)\bigg[ \int_{t^{i-1}_{k-1}}^{t^{i}_{k-1}} f(a)\, da  -  L^{\lambda(i)}\Delta t^{i-1}_{k-1}\bigg]\\
    \begin{split}\label{eq:IntermediatePosition}
          Z^{j}_{k-1} & = Z_{k-1} -r\dot{Z}_{k-1}\Delta t^{0}_{k-1} + \sum_{i=1}^{j} S(i)\bigg[\int_{t^{i-1}_{k-1}}^{t^{i}_{k-1}} F_1(a)\, da   \\
     & + F_1(t^{i-1}_{k-1}) \Delta t^{i-1}_{k-1} -\frac{L^{\lambda(i)}}{2}(\Delta t^{i-1}_{k-1})^2 \bigg]\,,\ j=1,\dots,N\,.
    \end{split}
\end{align}
The system of equations \eqref{eq:EqMotionVel}-\eqref{eq:EqMotionPosition} describes solutions that may or may not include sliding or crossing segments. To determine the occurrence of sliding or crossing segments on $\Sigma$, we can use these together with the conditions summarized in \cref{apx:Conditions}.

\cref{fig:NotationExplanation} depicts a 1:1$^{\Sigma}_{c}$ periodic solution to illustrate the notation in \eqref{eq:EqMotionVel} and \eqref{eq:EqMotionPosition}. Here, the motion from $\Gamma^{+}$ to $\Gamma^{-}$ is broken into three subintervals ($N=3$), namely $[t^{0}_{k-1},t^{1}_{k-1}]$, $[t^{1}_{k-1},t^{2}_{k-1}]$ and $[t^{2}_{k-1},t^{3}_{k-1}]$, where $t^{3}_{k-1}=t_k$, while the motion from $\Gamma^{-}$ to $\Gamma^{+}$ occurs on the interval $[t_{k},t_{k+1}]$. In addition, we have $S(1)=S(2)=S(3)=1$, since no sliding occurs, and $\lambda(1)=\lambda(3)=-$, while $\lambda(2)=+$. 

In (i)-(ix) below we define the maps between surfaces, which can be composed to capture the steps represented in \eqref{eq:EqMotionVel}-\eqref{eq:EqMotionPosition}. For example, the composition $P_{\Gamma^{-}\Gamma^{+}}\circ P_{\Sigma\Gamma^{-}}\circ P^{+}_{\Sigma\Sigma}\circ P_{\Gamma^{+}\Sigma}$ captures the segments of the solution shown in \cref{fig:NotationExplanation}. 
\begin{enumerate}[label=(\roman*)]
    \item $P_{\Gamma^{+}\Gamma^{-}}:(t_{k},Z_{k}\in\Gamma^{+},\dot{Z}_{k})\mapsto (t_{k+1},Z_{k+1}\in\Gamma^{-},\dot{Z}_{k+1})$ 
    \item $P_{\Gamma^{-}\Gamma^{+}}:(t_{k},Z_{k}\in\Gamma^{-},\dot{Z}_{k})\mapsto (t_{k+1},Z_{k+1}\in\Gamma^{+},\dot{Z}_{k+1})$\\
    The maps (i)-(ii) correspond to transitions between the impact surfaces with no events on $\Sigma$.
    \item $P_{\Gamma^{+}\Sigma}:(t_{k},Z_{k}\in\Gamma^{+},\dot{Z}_{k})\mapsto (t^{1}_{k},Z^{1}_{k},0)$
    \item $P_{\Gamma^{-}\Sigma}:(t_{k},Z_{k}\in\Gamma^{-},\dot{Z}_{k})\mapsto (t^{1}_{k},Z^{1}_{k},0)$
     \item $P_{\Sigma \Gamma^{+}}:(t^{N-1}_{k},Z^{N-1}_{k},0)\mapsto (t_{k+1},Z_{k+1}\in\Gamma^{+},\dot{Z}_{k+1})$
    \item $P_{\Sigma \Gamma^{-}}: (t^{N-1}_{k},Z^{N-1}_{k},0)\mapsto (t_{k+1},Z_{k+1}\in\Gamma^{-},\dot{Z}_{k+1})$\\
    The maps (iii)-(vi) correspond to transitions between the impact surfaces $\Gamma^{\pm}$ and $\Sigma$ through $\Sigma^{\pm}$.
    \item $P_{\Sigma\Sigma}^{+}: (t^{j}_{k},Z^{j}_{k},0)\mapsto (t^{j+1}_{k},Z^{j+1}_{k},0)$ (through $\Sigma^{+}=\{\dot{Z}>0\}$)
    \item $P_{\Sigma\Sigma}^{-}: (t^{j}_{k},Z^{j}_{k},0)\mapsto (t^{j+1}_{k},Z^{j+1}_{k},0)$ (through $\Sigma^{-}=\{\dot{Z}<0\}$).
    Maps (vii) and (viii) capture transitions from $\Sigma$ to $\Sigma$ through $\Sigma^{\pm}$, so that a crossing event takes place on $\Sigma$ at $t=t^{j}_{k}$.
     \item $P_{\Sigma\Sigma}^{s}: (t^{j}_{k},Z^{j}_{k},0)\mapsto (t^{j+1}_{k},Z^{j}_{k},0)$, where $t^{j+1}_{k}=t^{j}_{k}+\Delta S$, and $\Delta S$ is the duration of the sliding motion for this map.
\end{enumerate}
 The maps from (to) $\Gamma^{\pm}$ correspond to $Z_{k}=\pm d/2$ in \eqref{eq:EqMotionPosition}. Additionally, the maps initiated or ending at $\Sigma$ ($\dot{Z}^{j}_{k}=0$) correspond to individual terms of the sums in \eqref{eq:EqMotionVel}, \eqref{eq:EqMotionPosition} for properly assigned values of the indicator function, $S(j)$.
 In \cref{sec:StabilityAnalysis} we use \eqref{eq:EqMotionVel}-\eqref{eq:EqMotionPosition} to obtain expressions associated with the Jacobians of these maps that are necessary to analyze the linear stability.
\subsection{Analytical construction of key periodic solutions}\label{sec:Framework}
Together with additional periodicity conditions, the maps 
$P_{{\bullet \bullet}}$ in (i)-(ix) above provide a systematic way for deriving an analytical form of the periodic motions in terms of several quantities: the
impact velocity $\dot{Z}_{k-1}$, the time intervals between impacts on $\Gamma^{\pm}$, i.e. $\Delta t^{j-1}_{k-1} = t^{j}_{k}- t^{j-1}_{k-1}$, and the phase difference at the initial impact $\varphi_{k-1} = \mod(\pi t_{k-1} + \phi,2\pi)$. We may assume that $t_{k-1}=0$ and $\varphi_{k-1} = \phi$, for $\phi\in[0,2\pi)$. 

As mentioned in \cref{sec:AnalysisIntro}, the time intervals between impact times may be further divided into subintervals, during which sliding or non-stick (crossing) motions on $\Sigma$ occur. For the four types of periodic solutions described in \cref{sec:PhasePlane}, the motion from $\Gamma^{+}$ to $\Gamma^{-}$ is divided into $N$ subintervals of duration $\Delta t^{j}_{k-1}$, $j=0,\dots,N-1$, while the motion from $\Gamma^{-}$ to $\Gamma^{+}$ consists of one interval of duration $\Delta t_{k} =t_{k+1}-t_{k}$.
Each type of 1:1 solution is defined by a system of $2N+3$ equations, resulting from determining $Z$ and $\dot{Z}$ at the $N+1$ event times of interest (two impacts on $\Gamma^{\pm}$ and $N-1$ events at $\Sigma$) and from requiring that the solution is periodic. This system can be reduced to an associated system consisting of $N +2$ equations. Then the solution  
 \begin{align}\label{eq:Soltuple}
     (\dot{Z}_{k-1},\phi,\Delta t^{0}_{k-1}, \dots,\Delta t^{N-1}_{k-1})
 \end{align}

of the reduced system determines the different types of 1:1 solutions. Namely, we have:
\begin{itemize}
    \item Two equations capturing $\dot{Z}_k$ at each of the impacts at $\Gamma^{\pm}$ \cref{eq:EqMotionPosition}.
    \item $N-1$ equations that determine the intermediate times of events on $\Sigma$.
    \item One equation representing the periodicity condition, i.e. $\dot{Z}_{k+1}=\dot{Z}_{k-1}$.
\end{itemize}

\noindent We obtain the reduced system of equations in the variables \eqref{eq:Soltuple} as follows, assuming that $Z_{k-1}=d/2$ for simplicity: 
\begin{enumerate}[label=(\arabic*)]
\item Substitute all the intermediate displacements, $Z^{j}_{k-1}$ \eqref{eq:IntermediatePosition}, $j=1,\dots,N-1$ into the equation $Z_k=Z^{N}_{k-1} = -d/2$ successively. Similarly, substitute the intermediate displacement, $Z_{k}$, into \cref{eq:EqMotionPosition}, $Z_{k+1}=Z^{1}_{k} = +d/2$.
\item Substitute $\dot{Z}_k$ using its associated expression in terms of the variables in \eqref{eq:Soltuple} into \cref{eq:EqMotionVel} for $\dot{Z}_{k+1}$ . 
\item Use the periodicity conditions $\dot{Z}_{k+1}=\dot{Z}_{k-1}$, $F_1(t_{k+1})=F_1(t_{k-1})$, and $F_2(t_{k+1})=  F_2(t_{k-1})$ and eliminate any remaining quantities based on (1) and (2).
\item We find a solution  $(\dot{Z}_{k-1},\phi,\Delta t^{0}_{k-1}, \dots,\Delta t^{N-1}_{k-1})$ to the reduced system obtained from steps (1)-(3) numerically, using the Matlab solvers, fsolve and vpasolve. Note that this does not include inequality constraints, such as conditions when reaching $\Sigma$ or $\Gamma^{\pm}$, \eqref{eq:CrossingConditionsMinusToPlus_1} and
\eqref{eq:CrossingConditionsMinusToPlus_2} in \cref{sec:NonStickMotionCondsApp} and when sliding on $\Sigma$, \cref{eq:SlidingOnset2a,eq:SlidingOnset2b,eq:SlidingOnset2c} in \cref{sec:SlidingStick}, ensuring $\Delta t_{k-1}^{j}>0$, and so on. Therefore, appropriate constraints must be applied to remove solutions that violate these constraints, which we term as either \textit{infeasible}, e.g. those that yield $\Delta t_k^j<0$, or \textit{unphysical}, those that violate the model description, e.g. not applying appropriate equations of motion when crossing or sliding on $\Sigma$  or impacting $\Gamma^{\pm}$. Significantly, constraints for eliminating unphysical solutions correspond to grazing bifurcations on $\Sigma$ or $\Gamma$ as discussed below for \cref{eq:SwitchingDynamicsOnsetCondition} or \cref{eq:GrazingCondition}, respectively.
\end{enumerate}

\begin{remark}
The framework described here can be generalized to include periodic solutions that reach $\Sigma$ during the motion from $\Gamma^{-}$ to $\Gamma^{+}$. In this case, the reduced system of equations would consist of $N_u+N_d+1$ equations, where $N_d$ and $N_u$ are the number of subintervals composing the motion from $\Gamma^{-}$ to $\Gamma^{+}$ and from $\Gamma^{+}$ to $\Gamma^{-}$, respectively.
\end{remark}
\section{ Bifurcations of 1:1-type periodic solutions}\label{sec:ExamplesOfPeriodicMotions}

We analytically determine different types of periodic motions with alternating impacts on $\Gamma^{\pm}$ from reduced systems of equations, following steps (1)-(3) in \cref{sec:Framework}. Here we provide full details for 1:1 and 1:1$^{\Sigma}_{s}$ periodic solutions described in \cref{sec:PhasePlane}, with the details for others in \cref{sec:SlidingSol_MS} and \cref{sec:LoopingSlidingSolAnalysis}. We also provide conditions for grazing bifurcations associated with the switching manifold $\Sigma$ and capsule ends $\Gamma^{\pm}$. 
\subsection{1:1 periodic solution}\label{sec:SimpleSolAnalysis}

Below we describe the simple 1:1 periodic motion represented by the composition $P_{\Gamma^{-}\Gamma^{+}}\circ P_{\Gamma^{+}\Gamma^{-}}$, similar to
  $P_1 \circ P_2$ in \cite{Serdukova2019StabilityAB}.
  Here, the distinct vector fields in $\Sigma^{+}$ and $\Sigma^{-}$ when $\mu_k\neq0$ influence the analytical expressions capturing this motion. Recall that 1:1 periodic solutions do not cross or exhibit sliding on $\Sigma$.

\subsubsection*{Analytical expressions for the 1:1 periodic solution}

First we outline the key time points for 1:1 periodic motion using the notation defined in \eqref{eq:EqMotionVel},\eqref{eq:EqMotionPosition}:
\begin{itemize}
    \item First impact on $\Gamma^{+}$ $\left(Z=\frac{d}{2}\right)$ occurs at $t=t_{k-1}$.
    \item Second impact on $\Gamma^{-}$ $\left(Z=-\frac{d}{2}\right)$ occurs at $t=t_{k}$.
     \item Third impact on $\Gamma^{+}$ $\left(Z=\frac{d}{2}\right)$ occurs at $t=t_{k+1} = t_{k-1}+T$.
     \item The duration of the upward motion from $\Gamma^{+}$ to $\Gamma^{-}$ is given by $\Delta t_{k-1}$ (as defined in \eqref{eq:EqMotionVel},\eqref{eq:EqMotionPosition}). Similarly, the duration of the downward motion from $\Gamma^{-}$ to $\Gamma^{+}$ is $\Delta t_{k}$. Then, the period is $T=\Delta t_{k-1}+\Delta t_{k}$.
\end{itemize}
Since there are no intermediate intervals of interest (i.e. $N=1$ between impacts), we drop the superscript $j$ on the variables in \cref{eq:EqMotionVel,eq:EqMotionPosition}.
 These key time points are then used in the  full system of equations for a 1:1 periodic orbit following \cref{eq:EqMotionVel,eq:EqMotionPosition}.
\subsubsection*{From first to second impact, $P_{\Gamma^{+}\Gamma^{-}}: (t_{k-1},d/2,\dot{Z}_{k-1})\to(t_{k},-d/2,\dot{Z}_{k})$}
\begin{align}
\label{eq:FirstImpactZdot}
     \dot{Z}_k &= -r\dot{Z}_{k-1}+F_1(t_k)-F_1(t_{k-1})-L^{-}\Delta t_{k-1}\\\label{eq:FirstImpactZ}
    -d &= -r\dot{Z}_{k-1}\Delta t_{k-1}+F_2(t_k)-F_2(t_{k-1})-F_1(t_{k-1})\Delta t_{k-1}-L^{-}\frac{(\Delta t_{k-1})^2}{2}
\end{align}
\subsubsection*{From second to third impact, $P_{\Gamma^{-}\Gamma^{+}}: (t_{k},-d/2,\dot{Z}_{k})\to(t_{k+1},d/2,\dot{Z}_{k+1})$}
\begin{align}
\label{eq:SecondImpactZdot}
     \dot{Z}_{k+1} &= -r\dot{Z}_{k}+F_1(t_{k+1})-F_1(t_{k})-L^{+}\Delta t_{k}\\\label{eq:SecondImpactZ}
    d &= -r\dot{Z}_{k}\Delta t_{k}+F_2(t_{k+1})-F_2(t_{k})-F_1(t_{k})\Delta t_{k}-L^{+}\frac{(\Delta t_{k})^2}{2}
\end{align}
\subsubsection*{Periodicity conditions}
\begin{align}\label{eq:PeriodicityCond_SimpleSol}
   \dot{Z}_{k+1}&=\dot{Z}_{k-1}
\end{align}
In \cref{sec:ZeroFrictionComparison} we compare our calculations to the expressions obtained in \cite{Serdukova2019StabilityAB}, where we highlight the modifications in the solution due to friction.

\subsection*{Reduced system of equations for the 1:1 periodic solution}

We can reduce \eqref{eq:FirstImpactZdot}-\eqref{eq:PeriodicityCond_SimpleSol} and solve the resulting  subsystem to obtain a triple 
$(\dot{Z}_{k-1},\phi, \Delta t_{k-1})$, for a 1:1 periodic solution using steps (1)-(3) from \cref{sec:Framework}. 
\begin{align} \label{eq:RSE_SimpleSol}
\begin{split}
   \begin{cases}
     -d &= -r\dot{Z}_{k-1}\Delta t_{k-1}+F_2(t_k)-F_2(t_{k-1})-F_1(t_{k-1})\Delta t_{k-1}-L^{-}\frac{(\Delta t_{k-1})^2}{2}\\
  d&  =r^2\dot{Z}_{k-1}\Delta t_{k}-(r+1)F_1(t_k)\Delta t_{k}+rF_1(t_{k-1})\Delta t_{k}+F_2(t_{k-1})-F_2(t_{k})\\
  &+rL^{-}\Delta t_{k-1}\Delta t_{k}-L^{+}\frac{(\Delta t_{k})^2}{2}\\
  \dot{Z}_{k-1} &= \dfrac{1}{1-r^2}\bigg[-(r+1)(F_1(t_k)-F_1(t_{k-1}))+rL^{-}\Delta t_{k-1}-L^{+}\Delta t_{k}\bigg]
   \end{cases} 
   \end{split}
\end{align}
Results obtained by solving system \eqref{eq:RSE_SimpleSol} for certain parameter combinations are shown in \cref{fig:ZdotBifDiagram_muk=0p5_AnalyticalSols_A}, \cref{fig:LSA_OneOneBranch}, and \cref{fig:ParSet2_BifDiagrams} (black circles).
\subsubsection{Grazing-sliding bifurcation and unphysical 1:1 periodic solutions}\label{sec:SlidGrazingBif}
The numerical solutions of \eqref{eq:RSE_SimpleSol} may yield triples for solutions that reach and cross $\Sigma$, since the reduced system \eqref{eq:RSE_SimpleSol} does not include conditions to avoid this behavior, which violates the assumptions for 1:1 periodic solutions. 
Specifically, there may be a critical value of $d$ at which a 1:1 solution reaches $\Sigma$, with $\dot{Z} = 0$ and $\ddot{Z}=0$, so that the model trajectory exhibits a tangency at $\Sigma$ in the $Z-\dot{Z}$ phase plane. These conditions correspond to a grazing bifurcation on $\Sigma$ 
 at $d=d_{\Sigma}$, called a grazing-sliding bifurcation in \cite{DiBernardo_Review2008,Jeffrey_Hogan2011}. The critical $d_{\Sigma}$ is defined as follows: 
\begin{align}\label{eq:SwitchingDynamicsOnsetCondition}
\begin{split}
    d_{\Sigma} &= [d \,|\, \exists\, t_{\Sigma}\in(t_{k-1},t_{k})  \text{ s.t. }  \dot{Z} (t_{\Sigma})   = -r\dot{Z}_{k-1} + F_1(t_{\Sigma})-F_1(t_{k-1})-L^{-}(t_{\Sigma}-t_{k-1})=0\,,\\
    &\ddot{Z}(t_{\Sigma})=\cos(\pi t_{\Sigma}+\phi)-L^{-}=0,\, Z(t_{\Sigma})\in(-d/2,d/2)]\,.\\
\end{split} 
\end{align}

The grazing-sliding bifurcation on $\Sigma$ at $d=d_{\Sigma}$ marks the onset of a 1:1$^{\Sigma}_{s}$ solution whose sliding interval increases as $d$ decreases ($A$ increases). At $d=d_{\Sigma}$, the onset and exit times for sliding coincide, represented by $t_{\Sigma}$ which satisfies conditions \cref{eq:NonPassToSemiPassCondToOmega2_1,eq:NonPassToSemiPassCondToOmega2_2,eq:NonPassToSemiPassCondToOmega2_3} (where $t_{\Sigma}=t_{se}=t_{so}$).  Typically, triples obtained from system \eqref{eq:RSE_SimpleSol} for $d<d_{\Sigma}$ are flagged as unphysical 1:1 solutions.

\subsection{1:1$^{\Sigma}_{s}$ periodic solution}
Below we describe the 1:1$^{\Sigma}_{s}$ periodic motion represented by the composition $P_{\Gamma^{-}\Gamma^{+}}\circ P_{\Sigma\Gamma^{-}}\circ P^{s}_{\Sigma\Sigma}\circ P_{\Gamma^{+}\Sigma}$. The composition  $P_{\Sigma\Gamma^{-}}\circ P^{s}_{\Sigma\Sigma}\circ P_{\Gamma^{+}\Sigma}$
corresponds to the following motion: while the bullet is traveling from $\Gamma^{+}$ to $\Gamma^{-}$, the trajectory reaches $\Sigma$ ($\dot{Z}=0$), where its absolute velocity is equal to the absolute velocity of the capsule ($P_{\Gamma^{+}\Sigma}$). Then due to the dry friction, there is a sliding motion ($P^{s}_{\Sigma\Sigma}$). At the end of the sliding interval, the bullet continues to travel towards $\Gamma^{-}$ at a velocity larger than the velocity of the capsule ($\dot{Z}<0$, $P_{\Sigma\Gamma^{-}}$). As before, the map $P_{\Gamma^{-}\Gamma^{+}}$ 
corresponds to the simple motion during which the bullet travels from $\Gamma^{-}$ to $\Gamma^{+}$ without reaching $\dot{Z}=0$ from $\Sigma^{+}$, that is, its absolute velocity never becomes equal to the absolute velocity of the capsule. 

\subsubsection*{Analytical expressions for the 1:1$^{\Sigma}_{s}$ periodic solution}

First, we introduce notation corresponding to the key time points in 1:1$^{\Sigma}_{s}$ periodic motion:
\begin{itemize}
    \item First impact on $\Gamma^{+}$ $\left(Z=\frac{d}{2}\right)$ occurs at $t=t^{0}_{k-1}$, where $t^{0}_{k-1}=t_{k-1}$.
    \item Intersecting $\Sigma$ $(\dot{Z}=0)$ at $t= t^{1}_{k-1}$ (onset of sliding).
    \item Exiting sliding motion at $t=t^{2}_{k-1}$. 
    \item The duration of the sliding motion is $\Delta S =\Delta t^{1}_{k-1}=t^{2}_{k-1}-t^{1}_{k-1}$.
     \item Second impact on $\Gamma^{-}$ $\left(Z=-\frac{d}{2}\right)$ occurs at $t=t^{3}_{k-1}=t_{k}$.
     \item Third impact on $\Gamma^{+}$ $\left(Z=\frac{d}{2}\right)$ occurs at $t=t_{k+1}=t_{k}+\Delta t_{k}$.
\end{itemize}
Note that the onset time, $t^{1}_{k-1}$, and exit time, $t^{2}_{k-1}$, of sliding satisfy the following:
\begin{align}\label{eq:SlidingConditionInMS_1}
    \mod(\pi t^{1}_{k-1}+\phi, 2\pi)\in[\arccos(L^{+}),\arccos(L^{-}))\\\label{eq:SlidingConditionInMS_2}
     \mod(\pi t^{2}_{k-1}+\phi, 2\pi)=\arccos(L^{-})\,,
\end{align}
as determined by conditions \eqref{eq:SlidingOnset2a}-\eqref{eq:SlidingOnset2c} and \eqref{eq:NonPassToSemiPassCondToOmega2_1}-\eqref{eq:NonPassToSemiPassCondToOmega2_3}.
These key time points then appear in the full system of equations for the 1:1$^{\Sigma}_{s}$ periodic orbit, following  \cref{eq:EqMotionVel,eq:EqMotionPosition}. Note that $\dot{Z}^{0}_{k-1}=\dot{Z}_{k-1}$. The superscript $0$ in \eqref{eq:SlidingOnsetEquationZdotFromOmega2}-\eqref{eq:PeriodicityCond_SlidingSol} is included to highlight the   initial condition at the impact on $\Gamma^{+}$, as formulated in 
\cref{eq:EqMotionVel,eq:EqMotionPosition}. We drop it in \eqref{eq:RSESlidingSol}, since $\dot{Z}_{k-1}$ is the only impact velocity in the reduced system.

\subsubsection*{From first impact to sliding onset,  $P_{\Gamma^{+}\Sigma}: (t_{k-1},d/2,\dot{Z}_{k-1})\to(t^{1}_{k-1},Z^{1}_{k-1},0)$}
\begin{align}
\label{eq:SlidingOnsetEquationZdotFromOmega2}
     0&=\dot{Z}^{1}_{k-1} = -r\dot{Z}^{0}_{k-1}+F_1(t^{1}_{k-1})-F_1(t^{0}_{k-1})-L^{-}\Delta t^{0}_{k-1}\\\begin{split}\label{eq:SlidingOnsetEquationZFromOmega2}
    Z^{1}_{k-1} &= \frac{d}{2}-r\dot{Z}^{0}_{k-1}\Delta t^{0}_{k-1}+F_2(t^{1}_{k-1})-F_2(t^{0}_{k-1})-F_1(t^{0}_{k-1})\Delta t^{0}_{k-1}-L^{-}\frac{(\Delta t^{0}_{k-1})^2}{2}
      \end{split}
\end{align}
\subsubsection*{From sliding onset to sliding exit time, $P^{s}_{\Sigma\Sigma}: (t^{1}_{k-1},Z^{1}_{k-1},0)\to (t^{2}_{k-1},Z^{2}_{k-1},0)$}
\begin{align}\label{eq:ZdotEqSlidingExit_SlidingSol}
\dot{Z}^{2}_{k-1}&=\dot{Z}^{1}_{k-1}=0\\
\begin{split}\label{eq:ZEqSlidingExit_SlidingSol}
       Z^{2}_{k-1} & =Z^{1}_{k-1} 
    \end{split}
\end{align}
\subsubsection*{From exit from sliding to second impact,  $P_{\Sigma\Gamma^{-}}: (t^{2}_{k-1},Z^{2}_{k-1},0)\to (t^{3}_{k-1},-d/2,\dot{Z}^{3}_{k-1})$}
\begin{align}
\label{eq:SlidingtoSecondImpactZdotFromOmega2}
     \dot{Z}^{3}_{k-1} &= F_1(t^{3}_{k-1})-F_1(t^{2}_{k-1})-L^{-}\Delta t^{2}_{k-1}\\\label{eq:SlidingtoSecondImpactZdotFromOmega2}
    -\frac{d}{2} &= Z^{2}_{k-1}+F_2(t^{3}_{k-1})-F_2(t^{2}_{k-1})-F_1(t^{2}_{k-1})\Delta t^{2}_{k-1}-L^{-}\frac{(\Delta t^{2}_{k-1})^2}{2}
\end{align}
\subsubsection*{From second to third impact, $P_{\Gamma^{-}\Gamma^{+}}: (t_{k},-d/2,\dot{Z}_{k})\to(t_{k+1},d/2,\dot{Z}_{k+1})$}
We use \eqref{eq:SecondImpactZdot} and \eqref{eq:SecondImpactZ} for $t_{k} = t^{3}_{k-1}$, $\dot{Z}_{k}= \dot{Z}^{3}_{k-1}$, and $t_{k+1} = t^{1}_{k}$.
\begin{align}
\label{eq:SecondImpactZdotinSlidingSolution}
     \dot{Z}_{k+1} &= -r\dot{Z}^{3}_{k-1}+F_1(t^{1}_{k})-F_1(t^{3}_{k-1})-L^{+}\Delta t_{k}\\\label{eq:SecondImpactZinSlidingSolution}
    d&= -r\dot{Z}^{3}_{k-1}\Delta t_{k}+F_2(t^{1}_{k})-F_2(t^{3}_{k-1})-F_1(t^{3}_{k-1})\Delta t_{k}-L^{+}\frac{(\Delta t_{k})^2}{2}
\end{align}
\subsubsection*{Periodicity conditions}
\begin{align}\label{eq:PeriodicityCond_SlidingSol}
    \dot{Z}_{k+1}&=\dot{Z}^{0}_{k-1}
\end{align}
\subsubsection*{Reduced system of equations for 1:1$^{\Sigma}_{s}$ periodic solution}

 Following steps (1)-(3) described in \cref{sec:Framework}, we can reduce system \eqref{eq:SlidingOnsetEquationZdotFromOmega2}-\eqref{eq:PeriodicityCond_SlidingSol} to the following subsystem \eqref{eq:RSESlidingSol} and obtain a solution 
$(\dot{Z}_{k-1},\phi, \Delta t^{0}_{k-1},\Delta t^{1}_{k-1},\Delta t^{2}_{k-1})$, corresponding to a 
1:1$^{\Sigma}_{s}$ periodic solution. Recall that $\Delta t^{j-1}_{k-1} = t^{j}_{k-1}-t^{j-1}_{k-1}$, $j=1,2,3$,  with onset and exit of sliding determined from \cref{eq:SlidingConditionInMS_1,eq:SlidingConditionInMS_2}.
\begin{subequations}\label{eq:RSESlidingSol}
\begin{align}
  \label{eq:RSESlidingSol_1}
      0= -r\dot{Z}_{k-1}+F_1(t^{1}_{k-1})-F_1(t^{0}_{k-1})-L^{-}\Delta t^{0}_{k-1}\\\label{eq:RSESlidingSol_2}
      \mod(\pi t^{2}_{k-1}+\phi,2\pi)=\arccos{(L^{-})}\\\label{eq:RSESlidingSol_3}
          -d =-r\dot{Z}_{k-1}\Delta t^{0}_{k-1}+F_2(t^{1}_{k-1})-F_2(t^{0}_{k-1})-F_1(t^{0}_{k-1})\Delta t^{0}_{k-1}-L^{-}\frac{(\Delta t^{0}_{k-1})^2}{2}\\ \nonumber
       +F_2(t^{3}_{k-1})-F_2(t^{2}_{k-1})-F_1(t^{2}_{k-1})\Delta t^{2}_{k-1}-L^{-}\frac{(\Delta t^{2}_{k-1})^2}{2}\\\label{eq:RSESlidingSol_4}
  d=-(r+1)\Delta t_{k}F_1(t^{3}_{k-1})+r\Delta t_{k}F_1(t^{2}_{k-1})+rL^{-}\Delta t_{k}\Delta t^{2}_{k-1}\\\nonumber
  +F_2(t^{0}_{k-1})-F_2(t^{3}_{k-1})-L^{+}\frac{(\Delta t_{k})^2}{2}\\\label{eq:RSESlidingSol_5}
        \dot{Z}_{k-1} = -(r+1)F_1(t^{3}_{k-1})+rF_1(t^{2}_{k-1})+rL^{-}\Delta t^{2}_{k-1}+F_1(t^{0}_{k-1})-L^{+}\Delta t_{k}  
\end{align}
\end{subequations}

Note that once again the necessary variables are $(\dot{Z}_{k-1},\phi,\Delta t^{0}_{k-1})$, since $\Delta t^{1}_{k-1}$ and $\Delta t^{2}_{k-1}$ are implicit functions of $\dot{Z}_{k-1},\phi$ and $\Delta t^{0}_{k-1}$. The details of the calculations that yield system \eqref{eq:RSESlidingSol} are given in \cref{sec:RSESliding_Calculations}.
\begin{remark} The duration of the sliding interval $[t^{1}_{k-1},t^{2}_{k-1}]$ for the 1:1$^{\Sigma}_{s}$ solution is given by $\Delta S  =t^{2}_{k-1}-t^{1}_{k-1}$, with $\mod(\pi t^{2}_{k-1}+\phi,2\pi) = \arccos(L^{-})$. The maximal length of this sliding interval is $\Delta S= (\arccos{(L^{-})-\arccos{(L^{+})}})/\pi$ with $\mod(\pi t^{1}_{k-1}+\phi,2\pi) = \arccos(L^{+})\,.$
\end{remark}
\begin{remark}
 Note that $L^{+}-L^{-} = 2M\mu_kg\cos{\beta}/A$ is the length of the image of the harmonic forcing, $f(t)$, on the maximal sliding interval $[t^{1}_{k-1},t^{2}_{k-1}]$. Therefore, a larger value of the friction coefficient, $\mu_k$, leads to a longer maximal length of the sliding interval. 
\end{remark}

\subsubsection{Switching-sliding and crossing-sliding bifurcations on $\Sigma$ and unphysical 1:1$^{\Sigma}_{s}$ periodic solutions}\label{sec:SwithcingSlidAndCrossingSlidBifs}

The transition from sliding  to crossing typically occurs through a sequence of two bifurcations, that is, through the sequence 1:1$_s^\Sigma \to $ 1:1$_{cs}^\Sigma \to$ 1:1$_{c}^{\Sigma}$ for a sequence of decreasing $d$ in this application. 

The first is a switching-sliding bifurcation at 
$d=d_{ss}$ \cite{DiBernardo_Review2008,Jeffrey_Hogan2011} for 1:1$_{s}^\Sigma \to$ 1:1$_{cs}^\Sigma$, for which the time $t^{1}_{k-1}$ of the first intersection with $\Sigma$, corresponding to the sliding start time for 1:1$_{s}^{\Sigma}$ in \eqref{eq:RSESlidingSol}, reaches the boundary of the sliding interval defined by \eqref{eq:SlidingOnset2a}-\eqref{eq:SlidingOnset2c}. Then, $t^{1}_{k-1}$ corresponds to crossing at $d=d_{ss}$ instead of sliding. That is,
\begin{align}\label{eq:SwitchingSlidingBifOnsetCondition}   d_{ss} &= [d \,|\, \exists\, t_{ss}\in(t_{k-1},t^{2}_{k-1})\subset (t_{k-1},t_{k}) \text{ s.t. }
    \dot{Z} (t_{ss}) = -r\dot{Z}_{k-1} + F_1(t_{ss})-F_1(t_{k-1})-L^{-}(t_{ss}-t_{k-1})=0\,,\\\nonumber
   & \ddot{Z}(t_{ss})=\cos(\pi t_{ss}+\phi)-L^{+}=0\,,\dddot{Z}(t_{ss})=-\pi\sin(\pi t_{ss}+\phi)<0\,, Z(t_{ss})\in(-d/2,d/2)]\,.
\end{align}
To analytically obtain a 
1:1$_{cs}^{\Sigma}$ periodic solution captured by the map $P_{\Gamma^{-}\Gamma^{+}} \circ P_{\Sigma\Gamma^{-}} \circ P^{s}_{\Sigma\Sigma}\circ P^{+}_{\Sigma\Sigma}\circ P_{\Gamma^{+}\Sigma}$, appropriate intermediate times $t_{k-1}^{j}$ are included in the system 
\eqref{eq:EqMotionVel}-\eqref{eq:EqMotionPosition} to account for the crossing of $\Sigma$ at $t_{k-1}^{1}$ and subsequent dynamics in $\Sigma^{+}$ until the sliding interval $[t_{k-1}^{2},t_{k-1}^{3}]$ on $\Sigma$ defined by \eqref{eq:SlidingOnset2a}-\eqref{eq:SlidingOnset2c}.

A crossing-sliding bifurcation \cite{DiBernardo_Review2008, Jeffrey_Hogan2011} leads to the transition from 1:1$^{\Sigma}_{cs}$ to 1:1$^{\Sigma}_{c}$ periodic solutions. At the critical value  $d = {d_c}$, the trajectory from $\Sigma^{+}$ intersects $\Sigma$ at $t_c = t_{k-1}^2$ at the end of the (potentially) sliding region defined by \eqref{eq:SlidingOnset2a}-\eqref{eq:SlidingOnset2c}. Then, the 1:1$^{\Sigma}_{c}$ solution has two crossings of $\Sigma$ at $t^{1}_{k-1}$ and $t^{2}_{k-1}$ rather than sliding, defined respectively by \eqref{eq:CrossingConditionsMinusToPlus_1}-\eqref{eq:CrossingConditionsMinusToPlus_2} and \eqref{eq:CrossingConditionsPlusToMinus_1}-\eqref{eq:CrossingConditionsPlusToMinus_2}, and captured by the map $P_{\Gamma^{-}\Gamma^{+}} \circ P_{\Sigma\Gamma^{-}}\circ P^{+}_{\Sigma\Sigma}\circ P_{\Gamma^{+}\Sigma}$. In particular, 
  \begin{align}\label{eq:CrossingSlidingBifOnsetCondition}
\begin{split}
    d_{c} &= [d \,|\, \exists\, t_{c}\in(t^{1}_{k-1},t_{k})\subset (t_{k-1},t_{k}) \text{ s.t. }\dot{Z} (t_{c})   = -r\dot{Z}_{k-1} + F_1(t_{c})-F_1(t^{1}_{k-1})-L^{+}(t_{c}-t^{1}_{k-1})=0\,,\\
    &\ddot{Z}(t_{c})=\cos(\pi t_{c}+\phi)-L^{-}=0,\,\dddot{Z}(t_{c})=-\pi\sin(\pi t_{c}+\phi)<0\,,  Z(t_{c})\in(-d/2,d/2)]\,.\\
\end{split} 
\end{align}
Details for obtaining 1:1$^{\Sigma}_{c}$ and 1:1$^{\Sigma}_{cs}$ periodic solutions are given in \cref{sec:LoopingSolAnalysis} and \cref{sec:LoopingSlidingSolAnalysis}. Solutions obtained from system \eqref{eq:RSESlidingSol} (\eqref{eq:RSELoopingSolution}) for $d<d_{ss}$ ($d<d_{c}$) are flagged as unphysical 1:1$^{\Sigma}_{s}$ (1:1$^{\Sigma}_{cs}$) solutions.
 \subsection{Grazing bifurcations on $\Gamma^{+}$ and unphysical 1:1$^{\Sigma}_{c}$ solutions}\label{sec:GrazingOnGamma}
Here, we provide conditions for a 
grazing bifurcation denoted $d_{\Gamma^{+}}$ corresponding to intersections of trajectories with $\Gamma^{+}$. We determine the value $d_{\Gamma^{+}}$, as follows:
\begin{align}\label{eq:GrazingCondition}
  d_{\Gamma^{+}} = [d \,|\,\exists\, t_{\Gamma^{+}} \text{ s.t. } \dot{Z} (t_{\Gamma^{+}})   = F_1(t_{\Gamma^{+}})-F_1(t^{1}_{k-1})-L^{+}(t_{\Gamma^{+}}-t^{1}_{k-1})=0\,, \\\nonumber
  Z(t_{\Gamma^{+}})=d/2\,,\text{ where }t^{1}_{k-1}\in(t_{k-1},t_{\Gamma^{+}}),\text{ so that } \dot{Z}(t^1_{k-1})=0 \text{ and } Z(t^1_{k-1})\in(-d/2,d/2)]\,.
\end{align}
As seen in \cref{sec:ComparisonAndStabAnalysis}, for $d<d_{\Gamma^{+}}$, the 1:1$^{\Sigma}_{c}$ solutions obtained via \eqref{eq:RSELoopingSolution} are classified as unphysical.

\begin{remark}\label{rem:AnalyticsFor2TSolutions}
In \cref{fig:ZdotBifDiagrams_FrictionComp} there are several instances of
1:1/$2T$-type periodic solutions, also compared with analytical results, e.g. in \cref{fig:ZdotBifDiagram_muk=0p5_AnalyticalSols_A} and following. The steps in (1)-(3) in \cref{sec:Framework} can be adapted to obtain reduced systems of equations for these types of periodic solutions, providing an augmented system for $\dot{Z}_{k-1},\dots,\dot{Z}_{k+3}$. Specifically, the equations for the displacement $Z^{j}_{k-1}$ at times of intersection with $\Sigma$, $t^{j}_{k-1}$ in \eqref{eq:EqMotionPosition} need to be substituted into the equation for $Z_{k+3}$. The equations \eqref{eq:EqMotionVel} for the impact velocities, $\dot{Z}_{k+\ell}$, $\ell=0,1,2$ need to be substituted into the equation for $\dot{Z}_{k+3}$ and the periodicity conditions become $\dot{Z}_{k+3}=\dot{Z}_{k-1}$, $F_1(t_{k+3})=F_1(t_{k-1})$ and $F_2(t_{k+3})=F_2(t_{k-1})$. 
Similar adaptations are used for other types of 1:1/$2T$ solutions that cross or slide on $\Sigma$, discussed below.
\end{remark}
\begin{remark}\label{rem:Notationd2T}
 We obtain the critical values of $d$ at the grazing and sliding bifurcations of 1:1$/2T$-type solutions and denote them as $d^{2T}_{\bullet}$, using the appropriate subscripts ($\Sigma$, $c$, $ss$) given above. We summarize the notation in the table below:
 \begin{table}[H]
     \centering
     \begin{tabular}{|c|c|}
     \hline
        \text{Notation}  & \text{Description}  \\\hline
      $d_{\Gamma^{+}}$    &  \text{Grazing of a 1:1 periodic orbit on $\Gamma^{+}$ characterized by $\dot{Z}_k=0$,} \\
      & leading to an additional low-velocity impact.\\\hline
      $d_{\Sigma}$    &  Grazing of a 1:1 periodic orbit on $\Sigma$ characterized by $\dot{Z}(t_k^j)= 0$\\
      & indicating the start or end of sliding.\\\hline
         $d_{ss}$    &  Switching-sliding bifurcation of a 1:1 periodic orbit that indicates \\&the onset of crossing on $\Sigma$ following sliding.\\\hline
         $d_{c}$    &  Crossing-sliding bifurcation of a 1:1 periodic orbit that indicates \\& the onset of solutions that cross $\Sigma$ twice.\\\hline
             &  Grazing of a 1:1$/2T$ periodic orbit on $\Gamma^{+}$ characterized by $\dot{Z}_k=0$,\\
      $d_{\Gamma^{+}}^{2T}$ & leading to an additional low-velocity impact. The intersection with $\Gamma^{+}$ occurs\\
      & only during one $\Gamma^{+}\to\Gamma^{-}$ transition in the period doubled solution.\\\hline
             &  Grazing of a 1:1$/2T$ periodic orbit on $\Sigma$ characterized by $\dot{Z}(t_k^j)=0$,\\
      $d_{\Sigma}^{2T}$ & indicating the start or end of sliding. The intersection may occur\\
      & only during one $\Gamma^{+}\to\Gamma^{-}$ transition in the period doubled solution.\\\hline
      &  Switching-sliding bifurcation of a 1:1$/2T$ periodic orbit indicating \\
       $d_{ss}^{2T}$ & the start of crossing on $\Sigma$ following sliding on $\Sigma$. The intersection may occur\\
      & only during one $\Gamma^{+}\to\Gamma^{-}$ transition in the period doubled solution\\\hline
     \end{tabular}
     \caption{Summary of notation for the bifurcation parameter $d$ at non-smooth bifurcations.}
     \label{tab:Notation_d}
 \end{table}
\end{remark}
\subsection{Sequences of stable solutions via nonlinear maps and numerical simulations} \label{sec:StableSolComparison}

\begin{figure}[hbtp!]
    \begin{subfigure}[t]{0.45\textwidth}
\centering
    \includegraphics[scale=0.42]{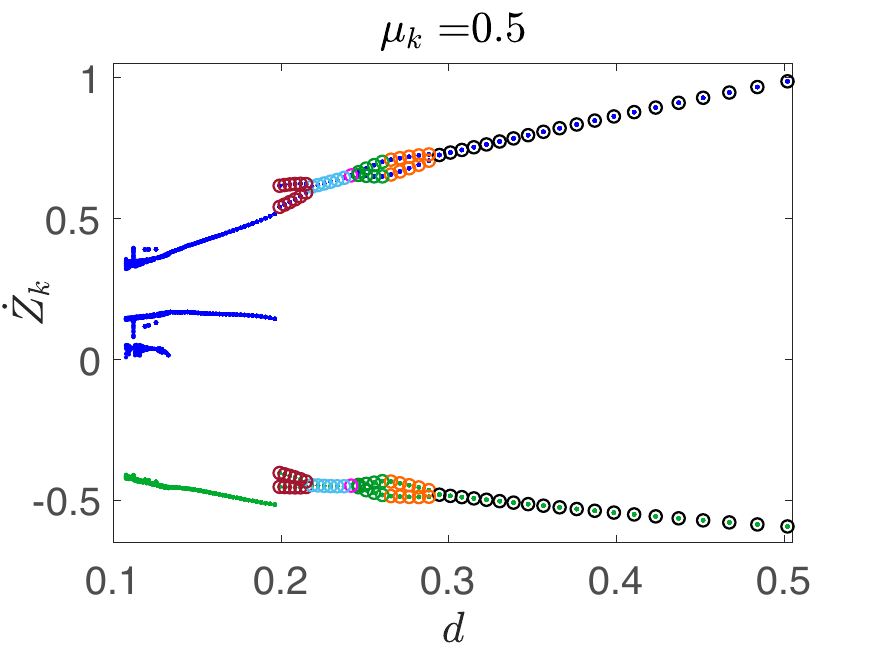}
    \caption{}
\label{fig:ZdotBifDiagram_muk=0p5_AnalyticalSols_A}
\end{subfigure}\hfill
    \begin{subfigure}[t]{0.45\textwidth}
\centering
    \includegraphics[scale=0.42]{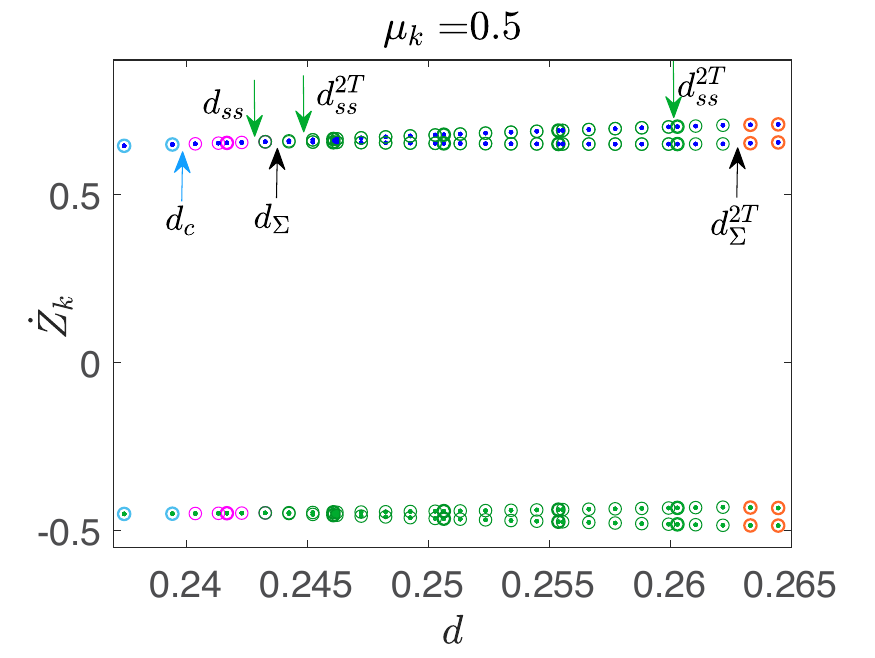}
    \caption{}
\label{fig:ZdotBifDiagram_muk=0p5_AnalyticalSols_A2}
\end{subfigure}
    \begin{subfigure}[t]{0.45\textwidth}
 \centering
    \includegraphics[scale=0.42]{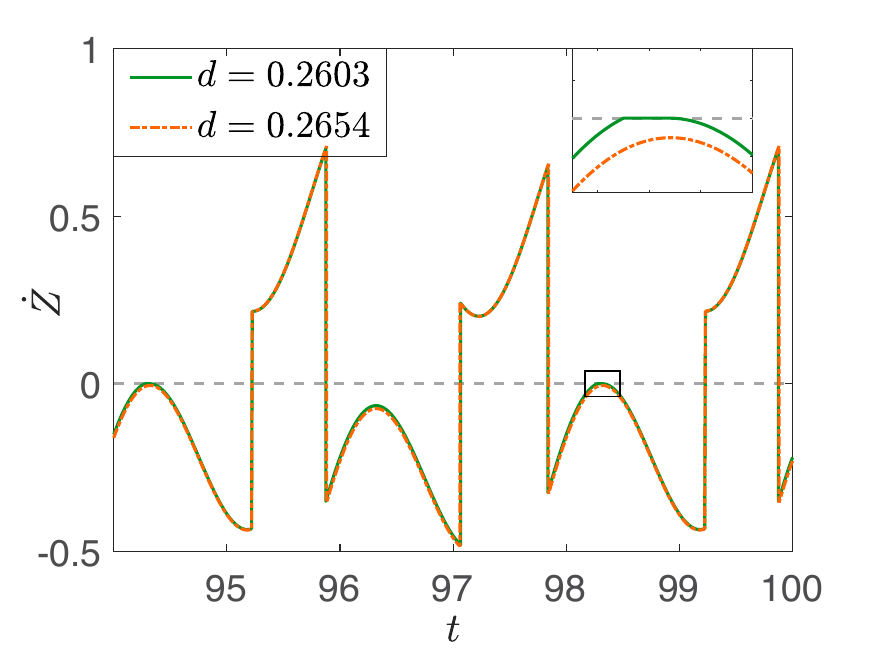}
    \caption{}
    \label{fig:ZdotBifDiagram_muk=0p5_AnalyticalSols_B}
\end{subfigure}\hfill
\begin{subfigure}[t]{0.45\textwidth}
 \centering
    \includegraphics[scale=0.42]{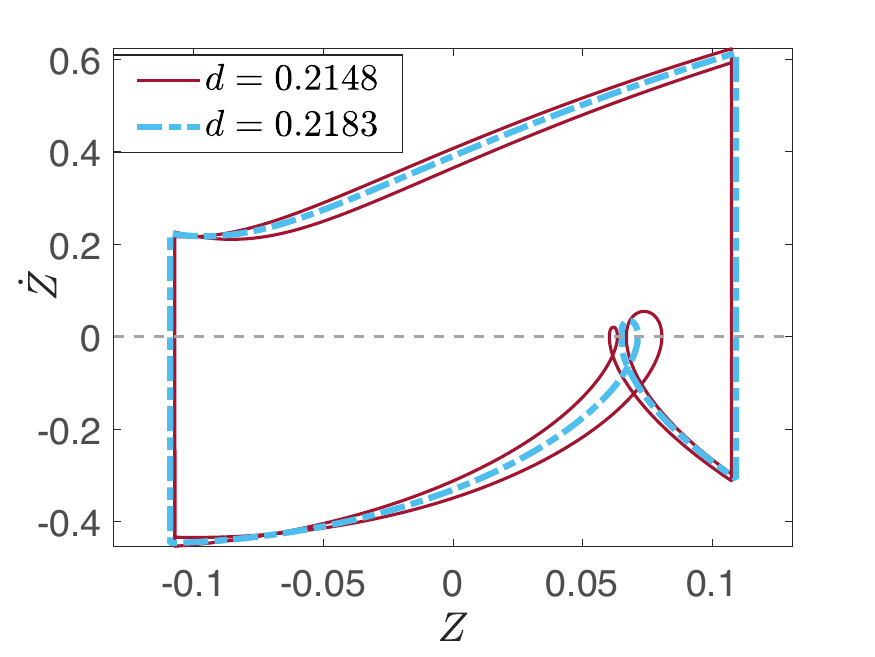}
    \caption{}
    \label{fig:ZdotBifDiagram_muk=0p5_AnalyticalSols_C}
\end{subfigure}
     \caption{(a),(b) Bifurcation diagram of impact velocities, $\dot{Z}_{k}$, for $\mu_k=0.5$. The remaining parameters are: $\beta=\pi/4$, $r=0.5$, $A\in[3.1,14.5]$. Panel (b) zooms in on $d\in[0.235,0.265]$. The blue and green dots correspond to impact velocities on $\Gamma^{+}$ and $\Gamma^{-}$ (at the bottom and top membrane), respectively, obtained from numerical simulations of \eqref{eq:PWSsystem}-\eqref{eq:ImpactCondsRelFrame}. The colored circles correspond to impact velocities of periodic solutions that were obtained analytically using \eqref{eq:EqMotionVel}-\eqref{eq:EqMotionPosition}. (a) The branch of 1:1-type solutions consists of the following periodic solutions: 1:1 (black circles), 1:1/$2T$ (orange circles), the sequence of solutions 1:1-1:1$^{\Sigma}_{s}/2T\to$1:1-1:1$^{\Sigma}_{cs}/2T\to$1:1-1:1$^{\Sigma}_{s}/2T\to$1:1$^{\Sigma}_{s}$ (green circles),  1:1$^{\Sigma}_{cs}$ (magenta circles), 1:1$^{\Sigma}_{c}$ (light blue circles), and 1:1$^{\Sigma}_{c}/2T$ (dark red circles). (b) The critical values of $d\in[0.235,0.265]$ at grazing-sliding ($d^{\bullet}_{\Sigma}$, \cref{sec:SlidGrazingBif}), switching-sliding ($d^{\bullet}_{ss}$, \cref{sec:SwithcingSlidAndCrossingSlidBifs}) and crossing-sliding ($d^{\bullet}_{c}$, \cref{sec:SwithcingSlidAndCrossingSlidBifs}) bifurcations generating the sequence of solutions $\textrm{1:1}/2T\rightarrow \textrm{1:1-1:1}^{\Sigma}_{s}/2T\rightarrow\textrm{1:1-1:1}^{\Sigma}_{cs}/2T\rightarrow\textrm{1:1-1:1}^{\Sigma}_{s}/2T\rightarrow\textrm{1:1}^{\Sigma}_{s}\rightarrow\textrm{1:1}^{\Sigma}_{c}$ are designated by black, green and light blue arrows, respectively. (c) Comparison of $\dot{{Z}}$ time series for a 1:1$/2T$ periodic solution at $d=0.2654$ (dashed-dotted orange curve) and a 1:1-1:1$^{\Sigma}_{s}/2T$ periodic solution at $d=0.2603\,(<d^{2T}_{\Sigma}=0.2633)$ (solid green curve). (d) Comparison of a 1:1$^{\Sigma}_{c}$ periodic solution at $d=0.2183$ (dashed-dotted light blue curve) and a 1:1$^{\Sigma}_{c}/2T$ periodic solution at $d=0.2148$ (solid dark red curve) in the $Z-\dot{{Z}}$ phase plane. Dashed gray lines in (c) and (d) correspond to $\dot{Z}=0$.}  \label{fig:ZdotBifDiagram_muk=0p5_AnalyticalSols}
\end{figure}

 We apply the analytical approach discussed above to capture stable and unstable solutions and to detect bifurcations. 
\Cref{fig:ZdotBifDiagram_muk=0p5_AnalyticalSols_A,fig:ZdotBifDiagram_muk=0p5_AnalyticalSols_A2} show impact velocities of 1:1/$pT$ periodic solutions obtained analytically (colored circles) that perfectly overlap with impact velocities obtained from numerical simulations of the model (green and blue dots). Note that \cref{fig:ZdotBifDiagram_muk=0p5_AnalyticalSols_A,fig:ZdotBifDiagram_muk=0p5_AnalyticalSols_A2} show the same bifurcation diagram as \cref{fig:ZdotBifDiagrams_FrictionComp_D} over a slightly larger interval of $d$. 

 The distinct colors of the circles differentiate the distinct types of 1:1 solutions, as well as the sequence of bifurcations for decreasing $d$. The black circles correspond to stable 1:1 solutions in which the velocity, $\dot{Z}$, does not reach $\Sigma$. The orange circles correspond to 1:1$/2T$ solutions, following a PD bifurcation at $d=0.2943$, as discussed in \cref{sec:FrameworkSection} and in more detail in \cref{sec:ComparisonAndStabAnalysis}. The green circles show a sequence of 1:1-$2T$-type solutions that exhibit sliding and crossing segments. In particular, at $d^{2T}_{\Sigma}= 0.2633$ (rightmost black arrow in \cref{fig:ZdotBifDiagram_muk=0p5_AnalyticalSols_A2}) a grazing-sliding bifurcation results in the gain of a $2T$-periodic solution that behaves like a combination of 1:1 solution (its velocity $\dot{Z}$ does not reach $\Sigma$) during the first transition from $\Gamma^{+}$ to $\Gamma^{-}$, while it behaves like a 1:1$^{\Sigma}_{s}$ solution (its velocity $\dot{Z}$ slides on $\Sigma$) during its second transition from $\Gamma^{+}$ to $\Gamma^{-}$. We denote this solution by 1:1-1:1$^{\Sigma}_{s}/2T$. At $d^{2T}_{ss}=0.2603$ (rightmost green arrow in \cref{fig:ZdotBifDiagram_muk=0p5_AnalyticalSols_A2}) a switching-sliding bifurcation leads to a 1:1-1:1$^{\Sigma}_{cs}/2T$ periodic solution that behaves like a 1:1 and a 1:1$^{\Sigma}_{cs}$ solution during alternating transitions $\Gamma^{+}\to\Gamma^{-}\to \Gamma^{+}$. At $d^{2T}_{ss}=0.24516$ (middle green arrow in \cref{fig:ZdotBifDiagram_muk=0p5_AnalyticalSols_A2}), a second switching-sliding bifurcation leads to a 1:1-1:1$^{\Sigma}_{s}/2T$ periodic solution. 

 At $d_{\Sigma}=0.2435$ (leftmost black arrow in \cref{fig:ZdotBifDiagram_muk=0p5_AnalyticalSols_A2}), a grazing-sliding bifurcation on $\Sigma$ leads to the transition 1:1-1:1$^{\Sigma}_{s}/2T\to$ 1:1$^{\Sigma}_{s}$. A switching-sliding bifurcation leads to a stable 1:1$^{\Sigma}_{cs}/T$ at $d_{ss}=0.24251$ (leftmost green arrow in \cref{fig:ZdotBifDiagram_muk=0p5_AnalyticalSols_A2}) introduces the stable 1:1$^{\Sigma}_{cs}/T$ solutions shown in magenta circles. At $d_{c}=0.2395$ (light blue arrow in \cref{fig:ZdotBifDiagram_muk=0p5_AnalyticalSols_A2}), a crossing-sliding bifurcation leads to the transition 1:1$^{\Sigma}_{cs}\to$ 1:1$^{\Sigma}_{c}$. In the 1:1$^{\Sigma}_{c}$ (light blue circles) periodic solutions $\dot{Z}$ crosses $\Sigma$. At $d=0.2183$ a PD bifurcation leads to the last sequence of 1:1-type solutions (dark red circles in \cref{fig:ZdotBifDiagram_muk=0p5_AnalyticalSols_A}), which is characterized by 1:1$^{\Sigma}_{c}/2T$ solutions in which the associated trajectories cross from $\Sigma^{-}$ to $\Sigma^{+}$ and then back to $\Sigma^{-}$ in each transition from $\Gamma^{+}$ to $\Gamma^{-}$. At $d^{2T}_{\Gamma^{+}}=0.1989$ a grazing bifurcation on $\Gamma^{+}$ marks the loss of stability of the 1:1$^{\Sigma}_{c}/2T$ solution and subsequent gain of the 2:1 periodic solution. The full sequence of stable solutions for decreasing $d$ (increasing $A$) can then be summarized as follows: 
\begin{eqnarray}
\label{eq:BifSequenceParSet1}
\begin{split}
    \textrm{1:1}\rightarrow&\textrm{1:1}/2T\rightarrow \textrm{1:1-1:1}^{\Sigma}_{s}/2T\rightarrow\textrm{1:1-1:1}^{\Sigma}_{cs}/2T\rightarrow\textrm{1:1-1:1}^{\Sigma}_{s}/2T\rightarrow\\&\textrm{1:1}^{\Sigma}_{s}\rightarrow\textrm{1:1}^{\Sigma}_{cs}\rightarrow\textrm{1:1}^{\Sigma}_{c}\rightarrow\textrm{1:1}^{\Sigma}_{c}/2T\rightarrow\textrm{2:1}\,.
    \end{split}
\end{eqnarray}
We note that in this parameter regime, we increase the excitation amplitude, $A$, thus, providing more energy to the VI-EH system. Even though higher friction slows down the device and leads to sliding, by increasing $A$, the system can overcome its dissipative influence. The consequent increase of energy output, discussed in \cref{sec:Energy}, is consistent with increasing $A$. At the same time,  the dissipation from friction  manifests itself through the quenching of the period-doubling regime via crossing and sliding on $\Sigma$, that allows for the reappearance of 1:1-type solutions.
 \begin{figure}[hbtp!]
    \begin{subfigure}[t]{0.5\textwidth}
    \centering
         \includegraphics[scale=0.45]{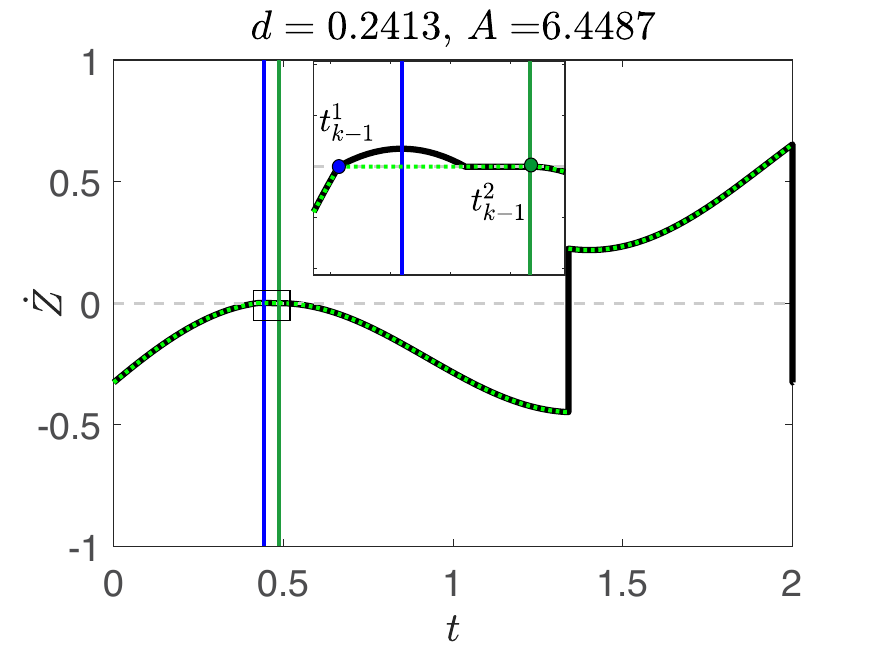}
         \caption{}
    \label{fig:UnphysicalSolutions_Simple}
    \end{subfigure}\hspace{-0.5cm}
      \begin{subfigure}[t]{0.5\textwidth}
      \centering
        \includegraphics[scale=0.45]{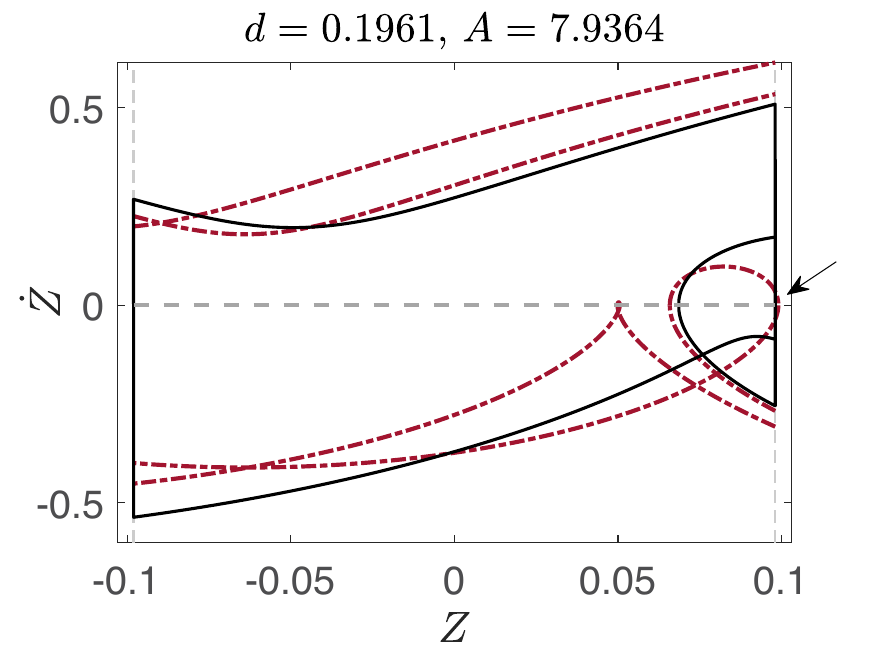}
    \caption{}
    \label{fig:UnphysicalSolutions_Looping}
    \end{subfigure}
 \caption{Illustration of unphysical 1:1-type solutions for different values of $d$ using the parameters in \cref{fig:ZdotBifDiagram_muk=0p5_AnalyticalSols}. (a) For $d=0.24133$ ($A=6.4487$, $d<d_{ss}<d_{\Sigma}$), the light green dotted trajectory shows the unphysical 1:1$^{\Sigma}_{s}$ periodic solution, as $t^{1}_{k-1}$(blue circle), the first time of intersection at $\Sigma$ (gray dashed line), violates condition \eqref{eq:SlidingConditionInMS_1} and lies outside the sliding interval (solid vertical blue and dark green lines). The black solid trajectory shows one period of the stable 1:1$^{\Sigma}_{cs}$ periodic solution. (b) For $d=0.1961$ ($A=7.9364$, $d<d^{2T}_{\Gamma^{+}}$), the dark red dashed trajectory corresponds to the unphysical (exceeds $\Gamma^{+}$, see black arrow) 1:1$^{\Sigma}_{c}/2T$ periodic solution. The black trajectory shows the stable 2:1 periodic solution.}
    \label{fig:UnphysicalSolutions}
\end{figure}
\Cref{fig:ZdotBifDiagram_muk=0p5_AnalyticalSols_B} shows a time series for the relative velocity, $\dot{Z}$ of a 1:1/$2T$ solution at $d=0.2654$ (orange curve) and a 1:1-1:1$^{\Sigma}_{s}/2T$ solution at $d=0.2603$ (green curve). \Cref{fig:ZdotBifDiagram_muk=0p5_AnalyticalSols_C} shows the phase plane for a 1:1$^{\Sigma}_{c}$ solution at $d=0.2183$ (light blue trajectory) and a 1:1$^{\Sigma}_{c}/2T$ solution at $d=0.2148$ (dark red trajectory).

\cref{fig:UnphysicalSolutions_Simple} illustrates an unphysical 1:1$^{\Sigma}_{s}$ periodic solution for $d<d_{ss} = 0.24251$, generated using
the quintuple obtained from \cref{eq:RSESlidingSol} at $d=0.2413$ yields a solution that slides along $\Sigma$ on $[t^{1}_{k-1}, t^{2}_{k-1}]$ indicated by the blue and green circles. This solution is unphysical since $t^{1}_{k-1}$, the first time of intersection with $\Sigma$, violates condition \eqref{eq:SlidingConditionInMS_1}, as it lies outside the actual sliding interval designated by the vertical blue and dark green lines, \eqref{eq:SlidingOnset2a}-\eqref{eq:SlidingOnset2c} and \eqref{eq:NonPassToSemiPassCondToOmega2_1}-\eqref{eq:NonPassToSemiPassCondToOmega2_3}. The black solid trajectory shows one period of the 1:1$^{\Sigma}_{cs}$ solution obtained from the full model \eqref{eq:PWSsystem}-\eqref{eq:ImpactCondsRelFrame}, using the same initial conditions as for the light green trajectory, and exhibiting non-stick and sliding motions at times of intersection with $\Sigma$ that respect the switching dynamics of the system.

 \cref{fig:UnphysicalSolutions_Looping} shows an unphysical 1:1$^{\Sigma}_{c}$ solution for 
$d< d^{2T}_{\Gamma^{+}}=0.1989$. At $d=0.1961$, the solution obtained from an adaptation for period-doubling of \eqref{eq:RSELoopingSolution}, violates the impacting dynamics by crossing $\Gamma^{+}$ (dark red trajectory). The black trajectory shows the stable 2:1-type solution at $d=0.1961$ obtained using \eqref{eq:PWSsystem}-\eqref{eq:ImpactCondsRelFrame}.

\section{Stability Analysis}\label{sec:StabilityAnalysis} The $T$-periodic motions discussed here are captured by composing maps of the form $P= P_{\Gamma^{-}\Gamma^{+}}\circ P_{\Gamma^{+}\dots \Gamma^{-}}$, where $P_{\Gamma^{+}\dots\Gamma^{-}}$ denotes an appropriate map composition from $\Gamma^{+}$ to $\Gamma^{-}$   
 as defined in (i)-(ix) in \cref{sec:AnalysisIntro}. For example, for 1:1 periodic motion  $P_{\Gamma^{+}\dots\Gamma^{-}} = P_{\Gamma^{+}\Gamma^{-}}$, while for 1:1$^{\Sigma}_{s}$ periodic solutions  $P_{\Gamma^{+}\dots\Gamma^{-}} = P_{\Gamma^{-}\Sigma}\circ P^{s}_{\Sigma\Sigma}\circ P_{\Gamma^{+}\Sigma}$. The map $P_{\Gamma^{+}\dots\Gamma^{-}}$ takes $\mathbf{H}_{k-1}=(t_{k-1},\dot{Z}_{k-1})$ to $\mathbf{H}_{k}=(t_{k},\dot{Z}_{k})$, and the map $P_{\Gamma^{-}\Gamma^{+}}$ takes $\mathbf{H}_{k}=(t_{k},\dot{Z}_{k})$ to $\mathbf{H}_{k+1}=(t_{k+1},\dot{Z}_{k+1})$. To perform a linear stability analysis about a periodic orbit of this map $(\mathbf{H}^{*}_{k-1},\mathbf{H}^{*}_{k})$ we consider a small perturbation, $\delta \mathbf{H}_{k-1}$ about $\mathbf{H}^{*}_{k-1}$ and write:
\begin{align}
 \delta \mathbf{H}_{k+1}&=  DP_{\Gamma^{-}\Gamma^{+}}(\mathbf{H}^{*}_{k})DP_{\Gamma^{+}\dots\Gamma^{-}}(\mathbf{H}^{*}_{k-1})\delta \mathbf{H}_{k-1}\,.
\end{align}
Therefore, we need to evaluate the corresponding Jacobians:
\begin{align}\label{eq:JacobiansGeneral}
DP_{\Gamma^{-}\Gamma^{+}}\cdot DP_{\Gamma^{+}\dots\Gamma^{-}}=\begin{bmatrix}
    \frac{\partial t_{k+1}}{\partial t_{k}} &\frac{\partial t_{k+1}}{\partial \dot{Z}_{k}} \\
     \frac{\partial \dot{Z}_{k+1}}{\partial t_{k}} & \frac{\partial \dot{Z}_{k+1}}{\partial \dot{Z}_{k}}
\end{bmatrix}\cdot\begin{bmatrix}
    \frac{\partial t_{k}}{\partial t_{k-1}} &\frac{\partial t_{k}}{\partial \dot{Z}_{k-1}} \\
     \frac{\partial \dot{Z}_{k}}{\partial t_{k-1}} & \frac{\partial \dot{Z}_{k}}{\partial \dot{Z}_{k-1}}
\end{bmatrix}\,.
\end{align}
Generally, the partial derivatives $\frac{\partial t_{j}}{\partial ()_{j-1}}$ are obtained using the Implicit Function Theorem on \cref{eq:EqMotionPosition}, while the  partial derivatives $\frac{\partial \dot{Z}_{j}}{\partial ()_{j-1}}$ are obtained by direct differentiation of \cref{eq:EqMotionVel} \cite{Luo2013VibroImpactDL,SHAW1983129}. To obtain expressions for the elements of $DP_{\Gamma^{+}\dots\Gamma^{-}}$ we use the Chain Rule due to the dependence of $t_k$ and $\dot{Z}_k$ on the intermediate times and displacements, $t^{j}_{k-1}$ and $Z^{j}_{k-1}$, $j=1,\dots,N-1$.

\subsection{Stability analysis of the 1:1 periodic motion}\label{sec:StabAnalysis_SimpleOneOne}
The 1:1 periodic motion is described by the following composition $P_{\Gamma^{-}\Gamma^{+}}\circ P_{\Gamma^{+}\Gamma^{-}}$. We find expressions for the partial derivatives in $DP_{\Gamma^{+}\Gamma^{-}}$ as follows:
\begin{align}\label{eq:dtkdZdotkm1}
     \frac{\partial t_{k}}{\partial \dot{Z}_{k-1}}  &= \frac{r\Delta t_{k-1}}{-r\dot{Z}_{k-1}+F_1(t_k)-F_1(t_{k-1})-L^{-}\Delta t_{k-1}} \\\label{eq:dtkdtkm1}
    \frac{\partial t_{k}}{\partial t_{k-1}} &=- \frac{r\dot{Z}_{k-1}-f(t_{k-1})\Delta t_{k-1}+L^{-}\Delta t_{k-1}}{-r\dot{Z}_{k-1}+F_1(t_k)-F_1(t_{k-1})-L^{-}\Delta t_{k-1}}\\\label{eq:dZdotkdtkm1}
    \frac{\partial \dot{Z}_{k}}{\partial t_{k-1}}&= \frac{\partial t_{k}}{\partial t_{k-1}}(f(t_{k})-L^{-})-(f(t_{k-1})-L^{-})\\
\label{eq:dZdotkdZdotkm1}
    \frac{\partial \dot{Z}_{k}}{\partial \dot{Z}_{k-1}} & = -r +\frac{\partial t_{k}}{\partial \dot{Z}_{k-1}}(f(t_k)-L^{-})
\end{align}
Similarly, we obtain the partial derivatives in $DP_{\Gamma^{-}\Gamma^{+}}$ that have the same form as \eqref{eq:dtkdZdotkm1}-\eqref{eq:dZdotkdZdotkm1} replacing $\Delta t_{k-1}$ with $\Delta t_{k}$ and $L^{-}$ with $L^{+}$ (see \eqref{eq:dtkp1dtk}-\eqref{eq:dZdotkp1dZdotk} in \cref{sec:JacobianTtoB}).
\subsection{Stability analysis for the 1:1$^{\Sigma}_{s}$ periodic solution}
The 1:1$^{\Sigma}_{s}$ periodic solution is described by the composition $P_{\Gamma^{-}\Gamma^{+}}\circ 
(P_{\Gamma^{-}\Sigma}\circ P^{s}_{\Sigma\Sigma}\circ P_{\Gamma^{+}\Sigma})$. 
Recall that $t^{j}_{k-1}=t^{j-1}_{k-1} + \Delta t^{j-1}_{k-1}$, where the subscript corresponds to the $k-1^{\textrm{st}}$ impact and the superscript $j$ corresponds to the $j^{\textrm{th}}$ event time between the $k-1^{\textrm{st}}$ and the $k^{\textrm{th}}$ impacts. The entries of the corresponding Jacobians are:
\begin{align}
    \frac{\partial t_k}{\partial t_{k-1}}= \frac{\partial t^{3}_{k-1}}{\partial t^{0}_{k-1}}&=  \frac{\partial t^{3}_{k-1}}{\partial Z^{2}_{k-1}} \frac{\partial Z^{2}_{k-1}}{\partial Z^{1}_{k-1}} \bigg(\frac{\partial Z^{1}_{k-1}}{\partial t^{0}_{k-1}}+\frac{\partial Z^{1}_{k-1}}{\partial t^{1}_{k-1}}\frac{\partial t^{1}_{k-1}}{\partial t^{0}_{k-1}}\bigg)\\
  \frac{\partial t_{k}}{\partial \dot{Z}_{k-1}} =   \frac{\partial t^{3}_{k-1}}{\partial \dot{Z}_{k-1}}&=  \frac{\partial t^{3}_{k-1}}{\partial Z^{2}_{k-1}}\frac{\partial Z^{2}_{k-1}}{\partial Z^{1}_{k-1}}\bigg( \frac{\partial Z^{1}_{k-1}}{\partial \dot{Z}_{k-1}}+\frac{\partial Z^{1}_{k-1}}{\partial t^{1}_{k-1}}\frac{\partial t^{1}_{k-1}}{\partial \dot{Z}_{k-1}}\bigg)\\
   \frac{\partial \dot{Z}_{k}}{\partial t_{k-1}} =     \frac{\partial \dot{Z}^{3}_{k-1}}{\partial t^{0}_{k-1}}&= [f(t^{3}_{k-1})-L^{-}]\frac{\partial t^{3}_{k-1}}{\partial t^{0}_{k-1}}\\
     \frac{\partial \dot{Z}_{k}}{\partial \dot{Z}_{k-1}}=    \frac{\partial \dot{Z}^{3}_{k-1}}{\partial \dot{Z}_{k-1}}&= [f(t^{3}_{k-1})-L^{-}]\frac{\partial t^{3}_{k-1}}{\partial \dot{Z}_{k-1}}
\end{align}
We find the partial derivatives at each intermediate point of the motion between $\Gamma^{+}$ and $\Gamma^{-}$:
\begin{enumerate}
    \item From $(t^{0}_{k-1},d/2,\dot{Z}_{k-1})$ to $(t^{1}_{k-1},Z^{1}_{k-1},0)$: 
    \begin{align}
    \frac{\partial t^{1}_{k-1}}{\partial t^{0}_{k-1}} & = \dfrac{f(t^{0}_{k-1})-L^{-}}{f(t^{1}_{k-1})-L^{-}}\\
    \frac{\partial t^{1}_{k-1}}{\partial \dot{Z}_{k-1}} &  = \dfrac{r}{f(t^{1}_{k-1})-L^{-}}\\
    \frac{\partial Z^{1}_{k-1}}{\partial \dot{Z}_{k-1}} & =   - r\Delta t^{0}_{k-1}\\
    \frac{\partial Z^{1}_{k-1}}{\partial t^{0}_{k-1}} 
    &=  F_1(t^{1}_{k-1})-F_1(t^{0}_{k-1}) -  f(t^{0}_{k-1})\Delta t^{0}_{k-1}
\end{align}
\item From $(t^{1}_{k-1},Z^{1}_{k-1},0)$ to $(t^{2}_{k-1},Z^{2}_{k-1},0)$:
    \begin{align}
        \frac{\partial t^{2}_{k-1}}{\partial t^{1}_{k-1}}& =\frac{\partial t^{2}_{k-1}}{\partial Z^{1}_{k-1}} =\frac{\partial Z^{2}_{k-1}}{\partial t^{1}_{k-1}}  =0\\
         \frac{\partial Z^{2}_{k-1}}{\partial Z^{1}_{k-1}} & = 1\\
  \frac{\partial Z^{1}_{k-1}}{\partial t^{1}_{k-1}} &= -r\dot{Z}_{k-1} + F_1(t^{1}_{k-1})-F_1(t^{0}_{k-1})-L^{-}\Delta t^{0}_{k-1} = 0  
\end{align}
\item From $(t^{2}_{k-1},Z^{2}_{k-1},0)$ to $(t^{3}_{k-1},-d/2,\dot{Z}^{3}_{k-1})$:
 \begin{align}
     \frac{\partial t^{3}_{k-1}}{\partial t^{2}_{k-1}} &=\frac{\partial \dot{Z}^{3}_{k-1}}{\partial t^{2}_{k-1}}= \frac{\partial \dot{Z}^{3}_{k-1}}{\partial Z^{2}_{k-1}} = 0\\
    \frac{\partial t^{3}_{k-1}}{\partial Z^{2}_{k-1}} & = -\frac{1}{F_1(t^{3}_{k-1})-F_1(t^{2}_{k-1})-L^{-}\Delta t^{2}_{k-1}}
\end{align}
\end{enumerate}
\begin{remark}
 The Jacobian corresponding to the transition from $\Gamma^{-}$ to $\Gamma^{+}$, $DP_{\Gamma^{-}\Gamma^{+}}$, involves the partial derivatives described by \eqref{eq:dtkp1dtk}-\eqref{eq:dZdotkp1dZdotk} in \cref{sec:JacobianTtoB}.
    \end{remark}
   \begin{remark} Similarly, we obtain expressions for the partial derivatives in the Jacobian matrices needed for the linear stability analysis of the 1:1$^{\Sigma}_{c}$ and 1:1$^{\Sigma}_{cs}$ periodic solutions. The details of such calculations are given in \cref{sec:LoopingSolAnalysis} and \cref{sec:LoopingSlidingSolAnalysis}.
\end{remark}
\section{Comparison of numerical and analytical results} \label{sec:ComparisonAndStabAnalysis}
We demonstrate the use of the analytical approach from previous sections to determine the existence of different types of physical or unphysical solutions and perform  linear stability analyses using the maps described in \cref{sec:StabilityAnalysis}. 
\subsection{Larger $\beta$, $r$, $\omega$} \label{sec:ComparisonAndStabAnalysis_LargePars}
\begin{figure}[hbtp!]
    \centering
    \begin{subfigure}[b]{0.28\textwidth}
    \centering
         \includegraphics[scale=0.38]{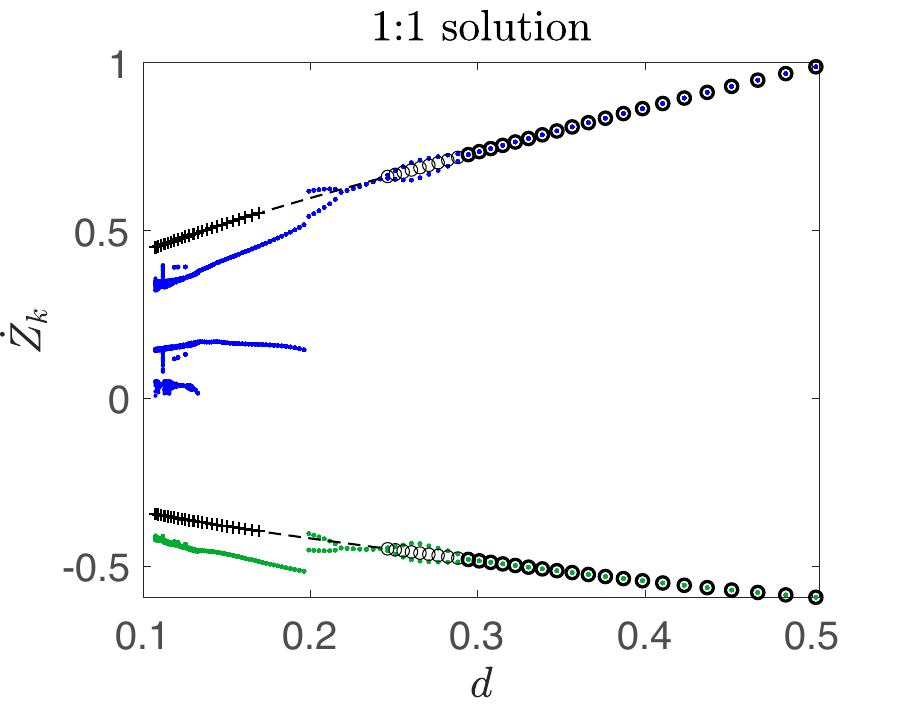}
    \caption{}
    \label{fig:SimpleSolutionBranch}
    \end{subfigure}
    \hfill
    \begin{subfigure}[b]{0.28\textwidth}
    \centering
         \includegraphics[scale=0.38]{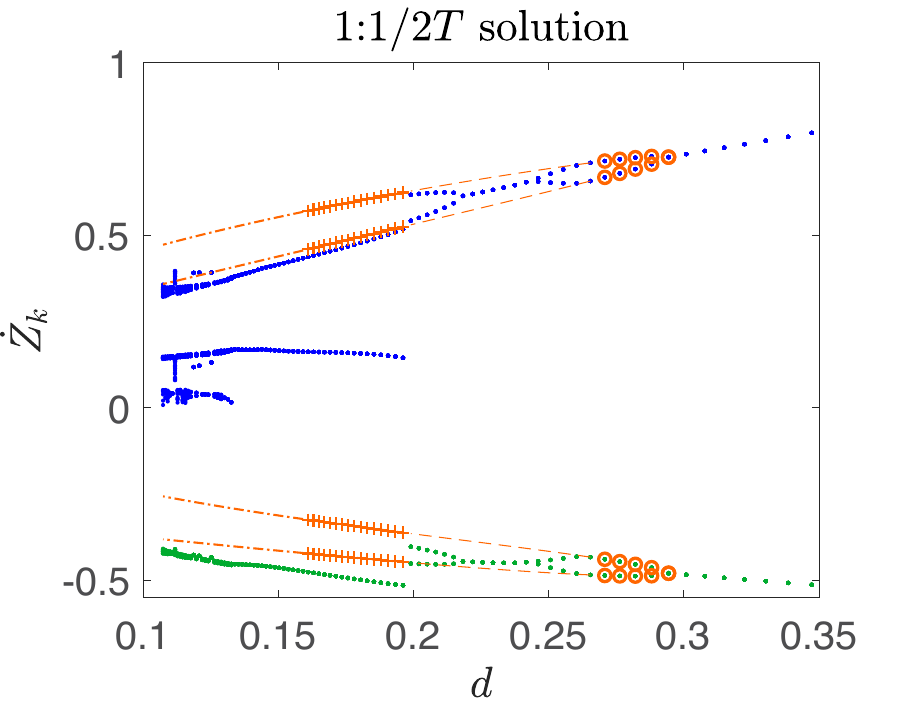}
    \caption{}
    \label{fig:PDSimpleSolutionBranch}
    \end{subfigure}
    \hfill
    \begin{subfigure}[b]{0.28\textwidth}
    \centering
         \includegraphics[scale=0.38]{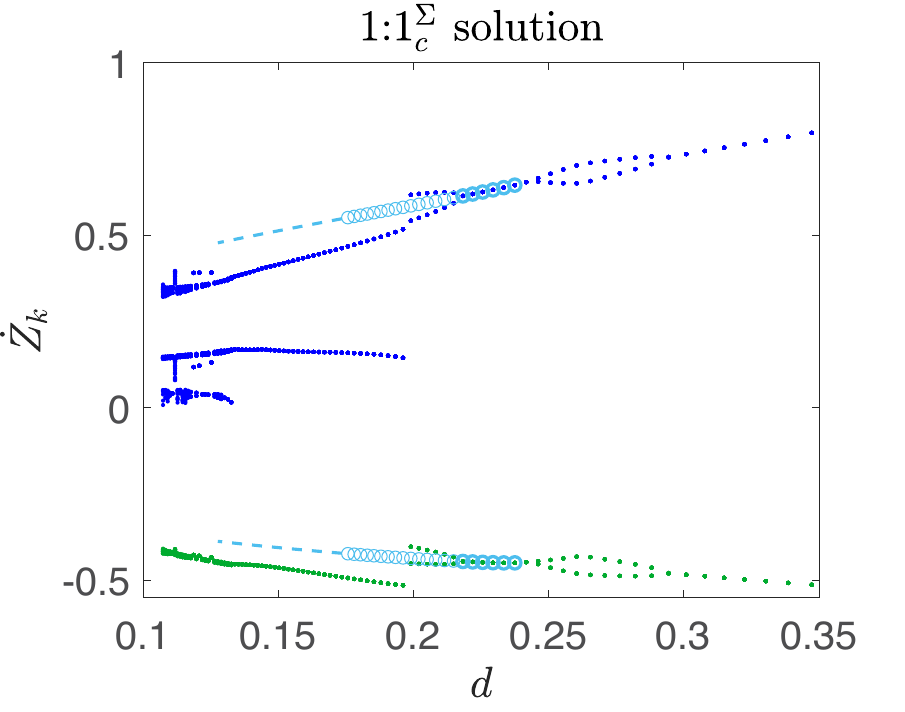}
    \caption{}
    \label{fig:LoopingSolutionBranch}
    \end{subfigure}
    \hfill
    \caption{Analytical results for stable, unstable, and unphysical (a) 1:1 solutions by solving system \eqref{eq:RSE_SimpleSol}, (b) 1:1/$2T$ solutions as discussed in \cref{rem:AnalyticsFor2TSolutions}, and (c) 1:1$^{\Sigma}_{c}$ solutions by solving system \eqref{eq:RSELoopingSolution} plotted over the numerical results (blue and green dots) and using the parameters in \cref{fig:ZdotBifDiagram_muk=0p5_AnalyticalSols_A}, i.e. $\mu_k=0.5$, $\beta = \pi/4$, $r=0.5$, $\omega =5\pi$, $s=0.5$, $A\in[3.1,14.5]$. Throughout the thick/thin circles correspond to a stable/unstable physical periodic solution (see corresponding eigenvalues in \Cref{fig:LSA_OneOneBranch_MS}). The various dashed lines, crosses (black in (a) and orange in (c)) and dashed-dotted lines correspond to unphysical solutions via a grazing bifurcation on $\Sigma$ (grazing-sliding) or $\Gamma^{+}$, since the switching or impacting dynamics of \eqref{eq:PWSsystem}-\eqref{eq:ImpactCondsRelFrame} are violated.} 
\label{fig:LSA_OneOneBranch}
\end{figure}

 In \cref{fig:LSA_OneOneBranch}, we show different 1:1$/T$ (1:1$/2T$ and 1:1$^{\Sigma}_{c}$\,) solutions, obtained numerically from \eqref{eq:PWSsystem}-\eqref{eq:ImpactCondsRelFrame} and compare the analytical results from the reduced systems \eqref{eq:RSE_SimpleSol} for the 1:1$/T$ solutions, the reduced system for 1:1$/2T$ solutions as described in \cref{rem:AnalyticsFor2TSolutions}, and \eqref{eq:RSELoopingSolution} for the 1:1$^{\Sigma}_{c}$ solutions. Different marker types show smooth and non-smooth bifurcations. 
Throughout dots indicate attracting solutions from numerical simulations and circles show solutions detected analytically. We also show unphysical solutions, for values of $d$ beyond certain grazing and sliding bifurcations of stable or unstable solutions. Moreover, as seen in \cref{fig:SimpleSolutionBranch,fig:PDSimpleSolutionBranch,fig:LoopingSolutionBranch}, we obtain unphysical 1:1, 1:1$/2T$, and 1:1$^{\Sigma}_{c}$ solutions over a large interval of $d$ values, confirming that the inclusion of additional constraints to the analytical systems of equations is important for the proper understanding of periodic behavior in the system (also illustrated in \cref{fig:UnphysicalSolutions}).  Additional grazing events in the unphysical branch of the 1:1 solution may also provide a lower bound for the parameter $d$ at bifurcations of physical 1:1$^{\Sigma}_{c}$-type solutions (\cref{fig:SimpleSolutionBranch}). As discussed in \cref{sec:StableSolComparison}, the transition 1:1$\to$1:1$^{\Sigma}_{c}$ occurs through a sequence of sliding bifurcations on $\Sigma$ which postpones the progression of physical solutions to a grazing bifurcation on $\Gamma^{+}$. However, such sliding events are not accounted for in the 1:1 unphysical solutions, so that  a grazing bifurcations on ${\Gamma^{+}}$ for the unphysical solutions is likely to occur at a larger $d$ than for the physical solutions.   
In \cref{fig:SimpleSolutionBranch} the thick black circles indicate stable 1:1$/T$ solutions, while the thin black circles correspond to unstable 1:1$/T$ solutions due to a PD bifurcation as indicated by the eigenvalues in the bottom panel of \Cref{fig:SimpleSolutionBranch_Eigs}. The calculation of eigenvalues follows from the linear stability analysis in \cref{sec:StabilityAnalysis} of the map $P_{\Gamma^{-}\Gamma^{+}}\circ P_{\Gamma^{+}\Gamma^{-}}$. {At $d<d_{\Sigma}=0.2435$ (\cref{sec:SlidGrazingBif}) the 1:1 solution undergoes a grazing-sliding bifurcation. Then dashed lines show unphysical solutions for $d<d_{\Sigma}$, obtained via Matlab from solutions to the \eqref{eq:RSE_SimpleSol}, but violating the switching conditions \eqref{eq:GeneralSlidingCondition} at times of intersection with $\Sigma$. } 

Similarly, in \cref{fig:PDSimpleSolutionBranch} the thick orange circles correspond to stable 1:1/$2T$ motion for $d\in (0.2633,0.2934]$, as indicated by the eigenvalues of the map $P_{\Gamma^{-}\Gamma^{+}}\circ P_{\Gamma^{+}\Gamma^{-}}\circ P_{\Gamma^{-}\Gamma^{+}}\circ P_{\Gamma^{+}\Gamma^{-}}$ shown in  \Cref{fig:PDSimpleSolutionBranch_Eigs} (bottom panel). At
$d\leq d^{2T}_{\Sigma} =0.2633$ 
there is a grazing-sliding bifurcation, with the superscript $2T$ indicating the use of a grazing condition similar to \eqref{eq:SwitchingDynamicsOnsetCondition} for the 1:1/$2T$ solution.
Unphysical solutions for $d\leq d^{2T}_{\Sigma} =0.2633$ violating the switching conditions \eqref{eq:GeneralSlidingCondition} at times of intersection with $\Sigma$ are shown as orange dashed lines. 

Finally, we obtain 1:1$^{\Sigma}_{c}$ periodic solutions by solving system \eqref{eq:RSELoopingSolution} in Matlab shown in \cref{fig:LoopingSolutionBranch}. The thick light blue circles correspond to stable 1:1$^{\Sigma}_{c}$ solutions, which lose stability via a PD bifurcation at  $d\approx 0.21479$, with eigenvalues shown in \Cref{fig:LoopingSolutionBranch_Eigs}.
The unstable 1:1$^{\Sigma}_{c}$ solutions are shown by
thin light blue  
circles for smaller $d$. Similarly to the other cases,
the unstable 1:1$^{\Sigma}_{c}$ solutions undergo a grazing bifurcation on $\Gamma^{+}$ and 
become unphysical for $d<d_{\Gamma^{+}} = 0.1736$ as given in \eqref{eq:GrazingCondition}. 

While unstable branches do not play a crucial role in the deterministic case, in stochastic settings we may observe solutions that approach unstable branches and thus, it is important to verify their existence in order to understand their influence on solution dynamics under different noise configurations. Furthermore, bifurcations on the unstable branches may offer bounds for bifurcations of stable solutions of a different type, as already noted in  \cite{Serdukova2022FundamentalCO}. For example, the grazing bifurcation of the unstable 1:1 periodic solution actually coincides with the grazing bifurcation that indicates the onset of a stable 1:1$^{\Sigma}_{s}$ periodic solutions at $d=d_{\Sigma}$ (\cref{fig:SimpleSolutionBranch}). Similarly, the grazing bifurcation on $\Gamma^{+}$ of the unstable 1:1$^{\Sigma}_{c}$ solution offers a lower bound on the value of $d^{2T}_{\Gamma^{+}}$, where the stable 1:1$^{\Sigma}_{c}/2T$ ceases to exist due to a grazing bifurcation on $\Gamma^{+}$ (\cref{fig:LoopingSolutionBranch}).
\subsection{Smaller $\beta$, $r$, larger $\omega$} \label{sec:ComparisonAndStabAnalysis_MixedPars}
\begin{figure}[hbtp!]
\begin{subfigure}[b]{0.45\textwidth}
     \centering
    \includegraphics[scale=0.44]{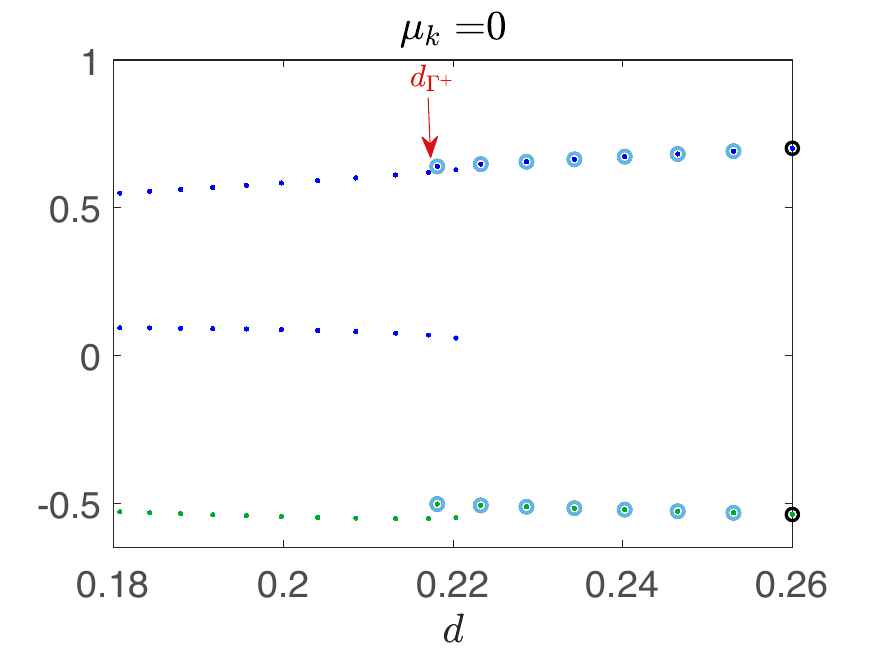}
    \caption{}
    \label{fig:ParSet2_BifDiagram_muk=0}
\end{subfigure}\hfill
\begin{subfigure}[b]{0.45\textwidth}
     \centering
    \includegraphics[scale=0.45]{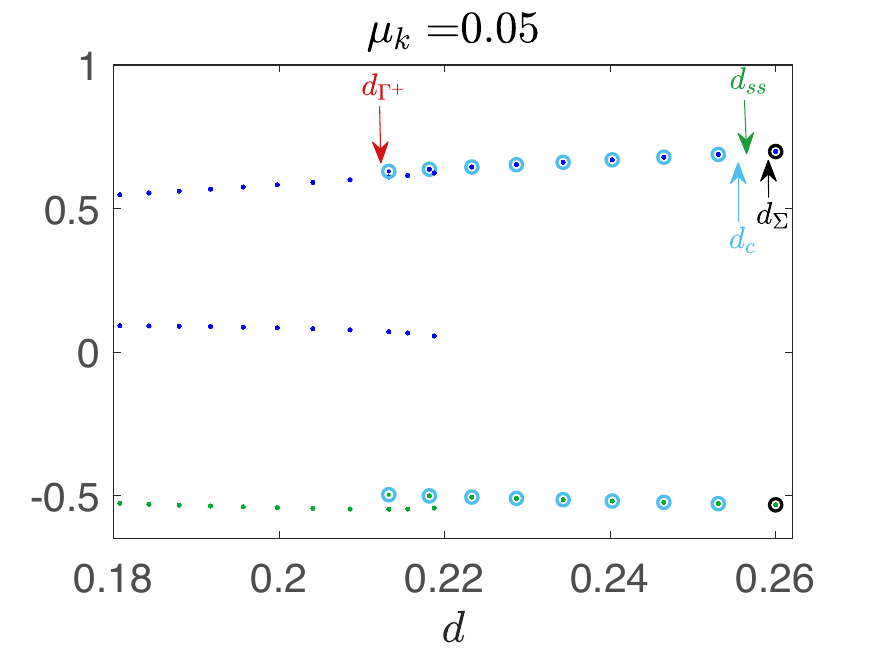}
    \caption{}
 \label{fig:ParSet2_BifDiagram_muk=0p05}
\end{subfigure}
\begin{subfigure}[b]{0.45\textwidth}
     \centering
    \includegraphics[scale=0.45]{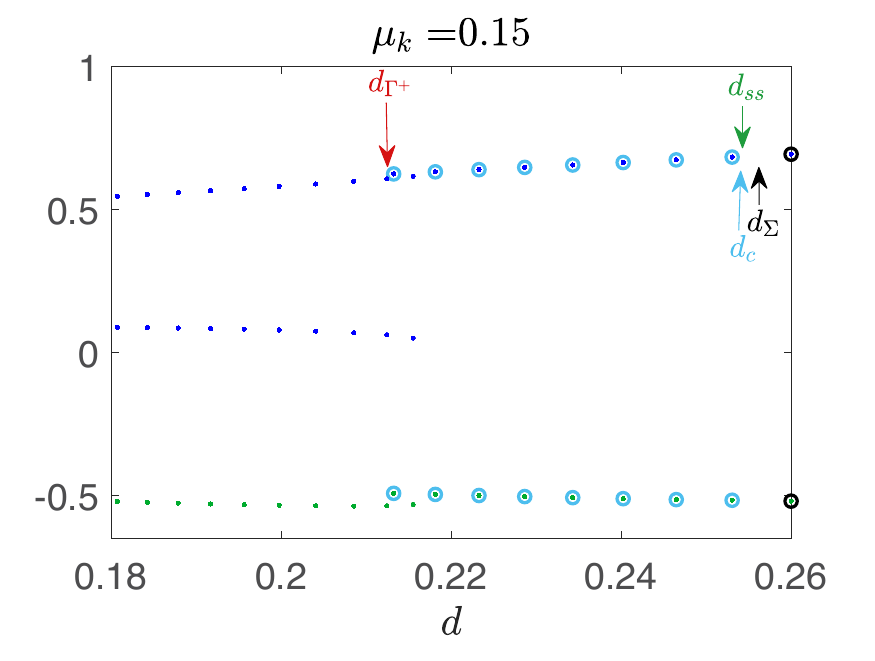}
    \caption{}
    \label{fig:ParSet2_BifDiagram_muk=0p15}
\end{subfigure}\hfill
\begin{subfigure}[b]{0.45\textwidth}
     \centering
    \includegraphics[scale=0.44]{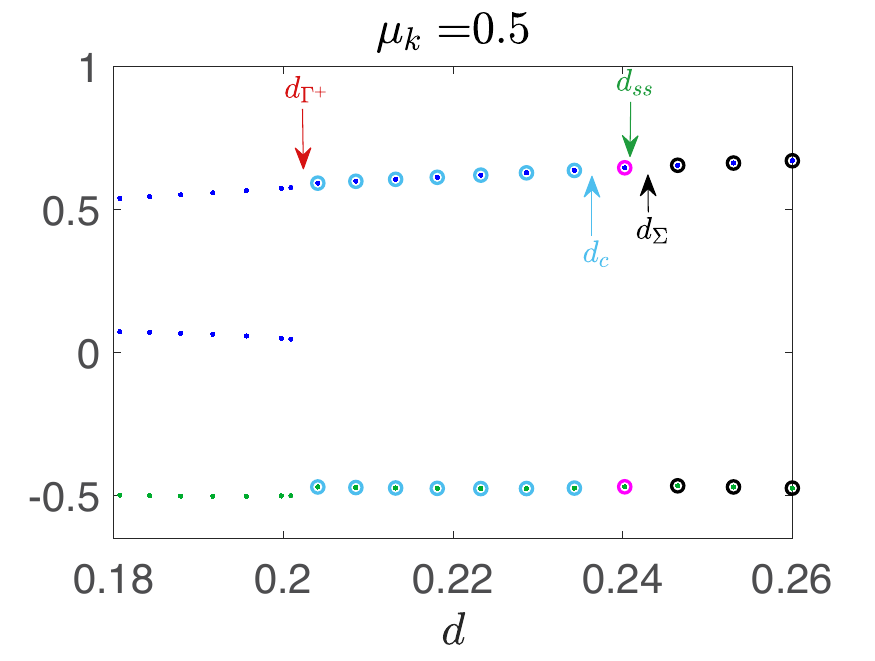}
    \caption{}
    \label{fig:ParSet2_BifDiagram_muk=0p5}
\end{subfigure}  
    \caption{Bifurcation diagrams of impact velocities, $\dot{Z}_{k}$, for (a) $\mu_k=0$, (b) $\mu_k=0.05$, (c) $\mu_k=0.15$, (d) $\mu_k=0.5$. The remaining parameters are: $\beta=\pi/6$, $r=0.25$, $A\in(5.986,8.6113)$, $s=0.5$, $\omega=5\pi$. The blue and green dots correspond to impact velocities $\dot{Z}_{k}$ on $\Gamma^{+}$ and $\Gamma^{-}$, respectively, obtained numerically from \eqref{eq:PWSsystem}-\eqref{eq:ImpactCondsRelFrame}. The thick black, magenta, and light blue circles indicate stable 1:1, 1:1$^{\Sigma}_{cs}$ and 1:1$^{\Sigma}_{c}$ solutions obtained by solving systems \eqref{eq:RSE_SimpleSol}, \eqref{sec:RSE_LoopingSlidingSol} and \eqref{eq:RSELoopingSolution}, respectively. Arrows indicate bifurcations at the critical values of $d$ for grazing-sliding ($d_{\Sigma}$), switching-sliding ($d_{ss}$), crossing-sliding ($d_{c}$) on $\Sigma$, and grazing (on $\Gamma^{+}$, $d_{\Gamma^{+}}$) bifurcations.}
    \label{fig:ParSet2_BifDiagrams}
\end{figure}

In contrast to the bifurcation sequence described in \cref{fig:ZdotBifDiagram_muk=0p5_AnalyticalSols} and \cref{fig:LSA_OneOneBranch}, \cref{fig:ParSet2_BifDiagrams} shows a different route, using small $r$ and $\beta$ ($r=0.25$, $\beta=\pi/6$), and larger $\omega$ ($\omega=5\pi$, same as before), for the transition from 1:1- to 2:1-type solutions for various $\mu_k$. This sequence does not involve a PD bifurcation in the transition from 1:1 to 2:1 solutions. The sequence does involve a grazing-sliding bifurcation (\cref{sec:SlidGrazingBif}), followed by switching-sliding and crossing-sliding bifurcations (\cref{sec:SwithcingSlidAndCrossingSlidBifs}) that lead to the transition from 1:1 to 1:1$^{\Sigma}_{c}$ solution, with a grazing on $\Gamma^{+}$ designating the onset of 2:1 behavior. The sequence of solutions for $\mu_k>0$ is 
\begin{eqnarray}
     \textrm{1:1}\rightarrow \textrm{1:1}^{\Sigma}_{s}\rightarrow \textrm{1:1}^{\Sigma}_{cs}\rightarrow\textrm{1:1}^{\Sigma}_{c}\rightarrow\textrm{2:1}\,.\label{bifseq_r25}
\end{eqnarray}
Once again, using our semi-analytical approach we identify the full sequence of bifurcations on that change the type of the stable periodic solution. The sequence of events on $\Sigma$ is simpler compared to \cref{eq:BifSequenceParSet1} and does not involve $2T$-periodic solutions. This highlights the influence of various parameters in the complexity of solution dynamics and their interaction with the switching manifold $\Sigma$.

Consistent with our previous observations, higher dry friction (larger $\mu_k$) slows down the system enough and shifts the grazing-sliding (on $\Sigma$) to grazing (on $\Gamma^{+}$) 
 sequence of bifurcations to smaller values of $d$. In particular, we solve the analytical system \eqref{eq:RSE_SimpleSol} and obtain 1:1 stable periodic solutions (thick black circles) and detect the grazing-sliding bifurcation for $\mu_k>0$ that flag the corresponding triples as unphysical. Our semi-analytical results show that the values of the critical parameters $d_{\Sigma}$ for the grazing-sliding bifurcation on $\Sigma$ are  $0.25778$, $0.25518$, and $0.24529$ for $\mu_k=0.05$ (\cref{fig:ParSet2_BifDiagram_muk=0p05}), $0.15$ (\cref{fig:ParSet2_BifDiagram_muk=0p15}), and $0.5$ (\cref{fig:ParSet2_BifDiagram_muk=0p5}), respectively, with the switching-sliding and crossing-sliding bifurcations occurring for smaller $d$ as indicated in \cref{fig:ParSet2_BifDiagrams}. For reference, the 1:1 solution reaches $\Sigma$ at $d=0.2591$ for $\mu_k=0$ (\cref{fig:ParSet2_BifDiagram_muk=0}). Similarly, we solve system \eqref{eq:RSELoopingSolution} to verify the 1:1$^\Sigma_{c}$ solutions obtained numerically using \eqref{eq:PWSsystem}-\eqref{eq:ImpactCondsRelFrame} (light blue circles). We also find that the loss of the 1:1$^\Sigma_{c}$ solution due to a grazing bifurcation on $\Gamma^{+}$ occurs at $d_{\Gamma^{+}}=0.21812$ for $\mu_k=0$, $d_{\Gamma^{+}}=0.21278$ for $\mu_k=0.05$, $d_{\Gamma^{+}}=0.21084$ for $\mu_k=0.15$, and $d_{\Gamma^{+}}=0.20162$ for $\mu_k=0.5$. 

Finally, a continuation-type method for both increasing and decreasing $d$ is used in \cref{fig:ParSet2_BifDiagrams}, showing bistability for $\mu_k=0, 0.05, 0.15$, in 
 panels (a)-(c). For $\mu_k=0$ and $d_{\Gamma^{+}}\leq d\leq0.22033$, for $\mu_k=0.05$ and $d_{\Gamma^{+}}\leq d\leq0.218704$, and for $\mu_k=0.05$ and $d_{\Gamma^{+}}\leq d\leq0.215516$ we detect both stable 2:1 solutions numerically, and  1:1$^{\Sigma}_{c}$ solutions obtained numerically and with our semi-analytical approach based on \eqref{eq:RSELoopingSolution}. Interestingly, for the largest value
shown $\mu_k=0.5$, we do not observe bistability, suggesting that dry friction may limit such regions, since it shifts grazing bifurcations on $\Gamma^+$ as discussed in the next section. We also note, that in these parameter regimes, where bistability is prevalent near the grazing bifurcations, distinct solutions have close impact velocities, and therefore, identifying them numerically may be more challenging due to sensitivity in numerical error. The maps guide our choice of initial conditions for continuing existing or finding new branches of the bifurcation diagram and for numerical validation.
\begin{figure}[hbtp!]
   \begin{subfigure}[b]{0.45\textwidth}
   \centering
       \includegraphics[scale=0.32]{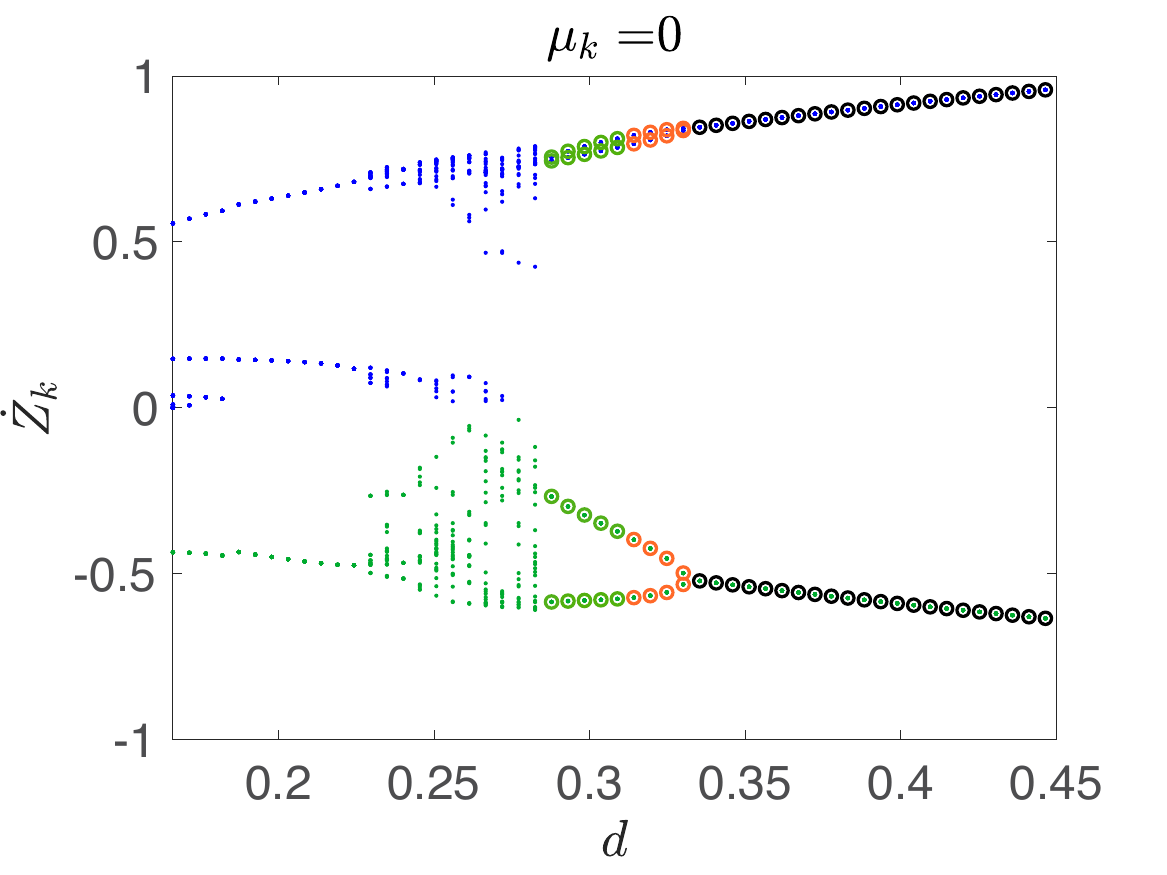}
\caption{}\label{fig:BifDiagram_omega=2pi_muk=0}
   \end{subfigure}\hfill
   \begin{subfigure}[b]{0.45\textwidth}
   \centering
       \includegraphics[scale=0.32]{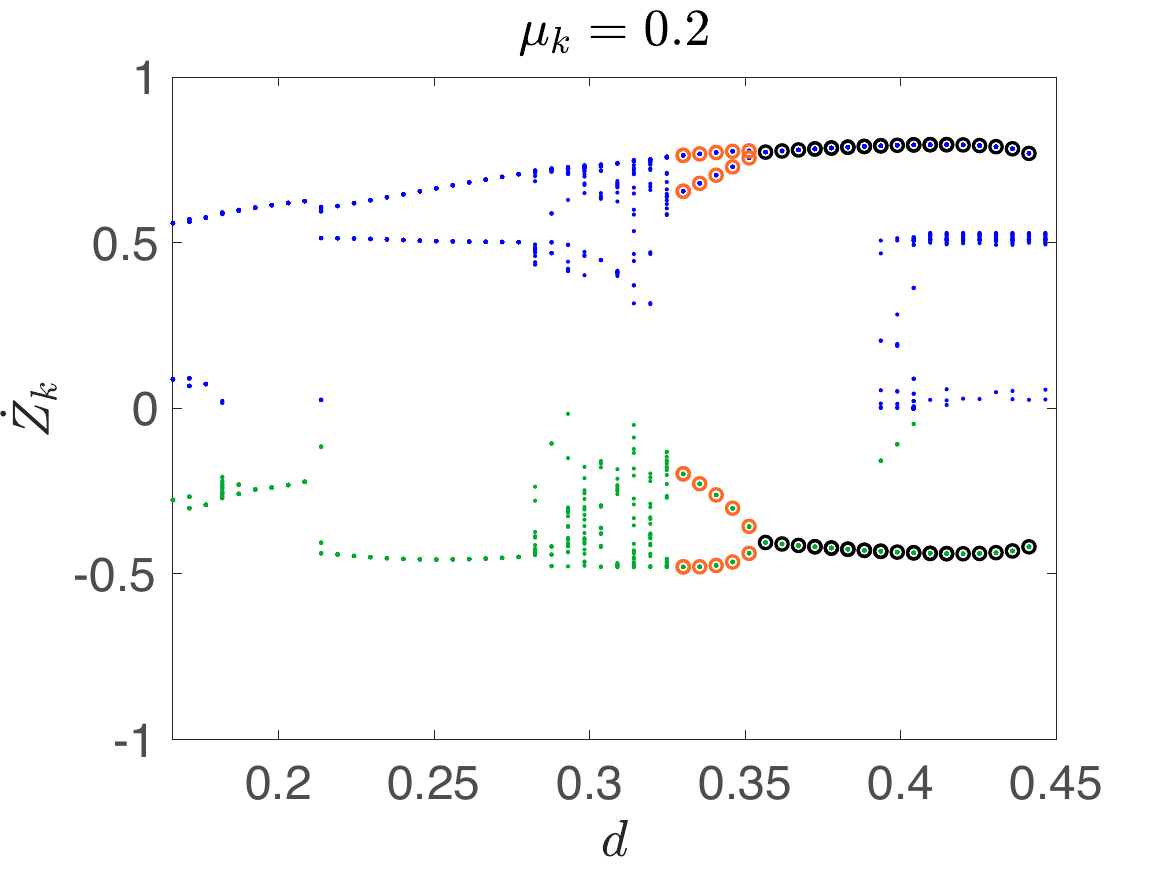}
       \caption{}  \label{fig:BifDiagram_omega=2pi_muk=0p2}
   \end{subfigure}
  \begin{subfigure}[b]{0.45\textwidth}
   \centering
       \includegraphics[scale=0.32]{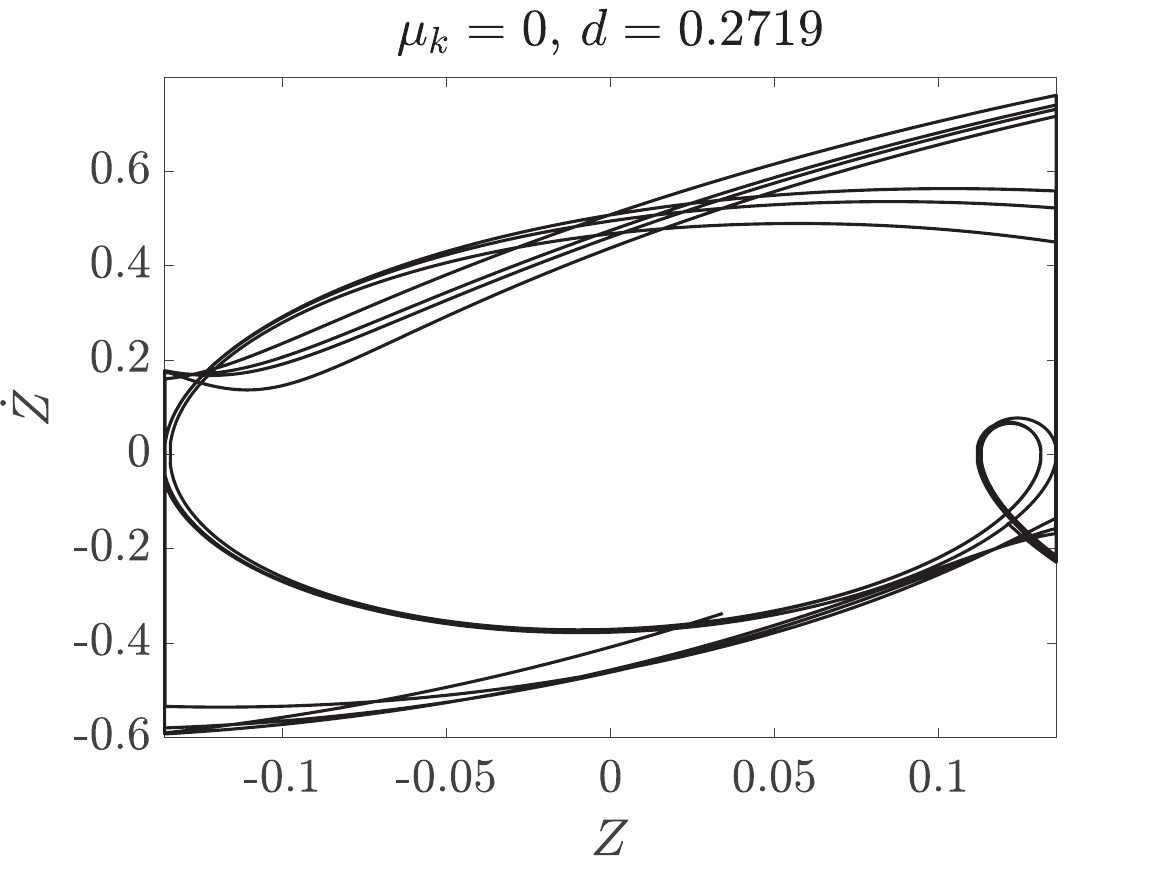}
       \caption{}  \label{fig:PhasePlane_omega=2pi_muk=0_d=0p2719}
   \end{subfigure} \hspace{1.25cm}
   \begin{subfigure}[b]{0.45\textwidth}
   \centering
       \includegraphics[scale=0.32]{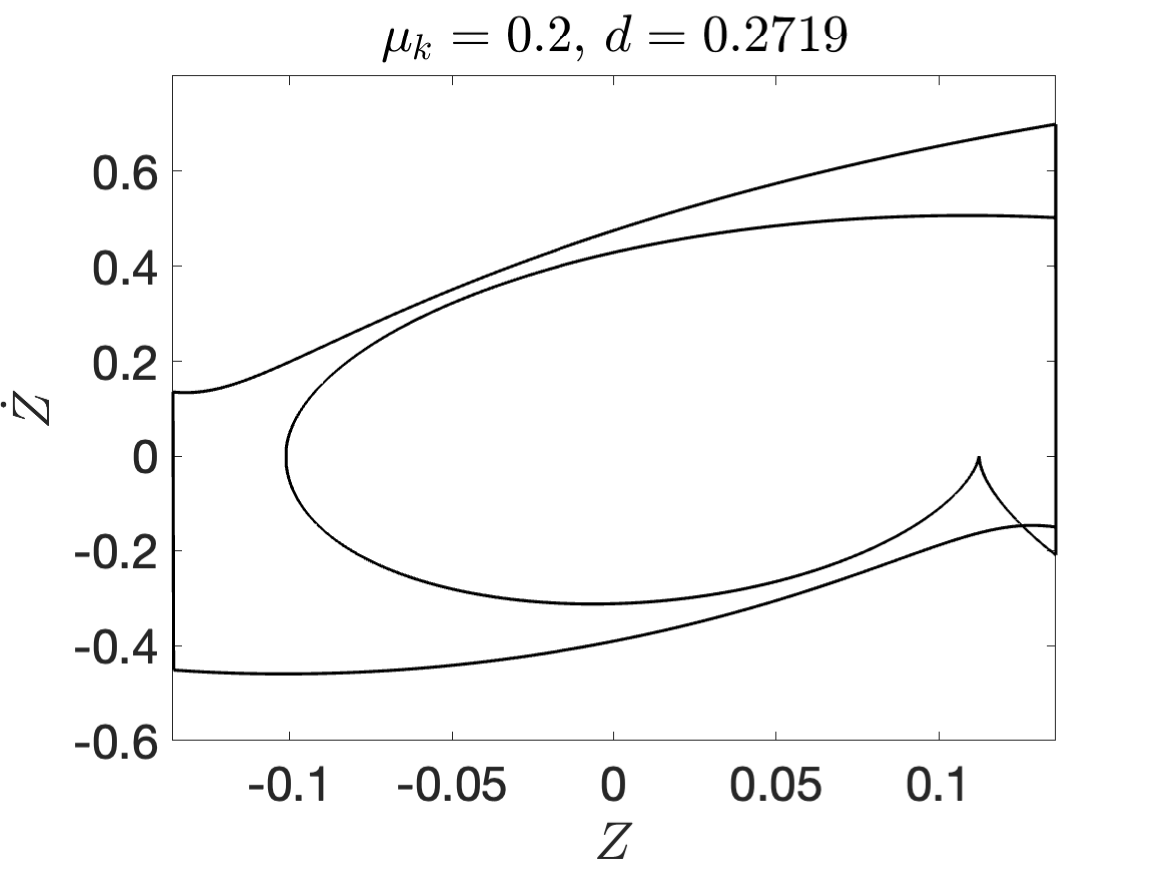}
       \caption{}\label{fig:PhasePlane_omega=2pi_muk=0p2_d=0p2719}
   \end{subfigure}
    \caption{(a),(b): Bifurcation diagrams of impact velocities, $\dot{Z}_{k}$, for (a) $\mu_k=0$, (b) $\mu_k=0.2$. The remaining parameters are: $\omega=2\pi$, $\beta=\pi/12$, $r=0.3$, $A=1.5$, $s\in(0.5,1.35)$. The blue and green dots correspond to impact velocities $\dot{Z}_{k}$ on $\Gamma^{+}$ and $\Gamma^{-}$, respectively, obtained numerically from \eqref{eq:PWSsystem}-\eqref{eq:ImpactCondsRelFrame}. The black, orange and green circles correspond to stable 1:1$/T$, 1:1$/2T$ and 1:1-1:1$^{\Sigma}_{c}/2T$ obtained analytically from \eqref{eq:EqMotionVel}-\eqref{eq:EqMotionPosition}. (c),(d): Phase planes for $d=0.2719$ for (c) $\mu_k=0$, (d) $\mu_k=0.2$. For $\mu_k=0$ the VI-EH system yields a chaotic solution exhibiting occasional small impact velocities on $\Gamma^{\pm}$. For $\mu_k=0.2$, the system yields a 1:1-0:1$^{\Sigma}_{cs}/2T$ periodic solution that impacts once on $\Gamma^{+}$ and half the time on $\Gamma^{-}$ per period $T$. The solution also exhibits sliding during the creation of the loop in the $Z-\dot{Z}$ plane.}
    \label{fig:BifDiagrams_omega=2pi}
\end{figure}
\subsection{Smaller $\beta$, $r$, $\omega$}\label{sec:ComparisonAndStabAnalysis_SmallPars}
A different type of bistability is also observed in smaller frequency regimes when friction is present in the VI-EH system. \Cref{fig:BifDiagram_omega=2pi_muk=0p2} shows the $\dot{Z}_k$ bifurcation diagram for the parameters $\omega=2\pi$, $\beta=\pi/12$, $r=0.3$, $A=1.5$, $s\in(0.5,1.35)$ and $\mu_k=0.2$. Note that for $d>0.393$ we observe solutions that impact multiple times on $\Gamma^{+}$, but never or rarely on $\Gamma^{-}$, resembling the dynamics of a bouncing ball. However, for $d<0.4413$, the system also yields a stable 1:1 solution. This type of bistability is not observed for $\mu_k=0$ in the same $d$-range (\cref{fig:BifDiagram_omega=2pi_muk=0}). In both cases, the 1:1 solution is followed by a PD bifurcation at $d=0.3353$ and $d =0.3565$ for $\mu_k=0$ and $\mu_k=0.2$, respectively. For $\mu_k=0$ and decreasing $d$ the 1:1$/2T$ solution is followed by a chaotic regime for $d<0.2877$ that introduces impacts on $\Gamma^{+}$ with some impact velocities $\dot{Z}_k$ close to $0$ (see for example, \cref{fig:PhasePlane_omega=2pi_muk=0_d=0p2719}). For $d<0.229$ the system yields stable 2:1$/T$ solutions. For $\mu_k=0.2$ the 1:1$/2T$ solution is also followed by a chaotic solution with additional impacts on $\Gamma^{+}$. However, the magnitude of those impact velocities $\dot{Z}_k$ range between 0.3 and 0.5. For $d<0.282$ and $\mu_k=0.2$ the system yields a stable solution that repeats every $2T$ and impacts on $\Gamma^{+}$ once per period and on $\Gamma^{-}$ every other period $T$. We denote this solution by 1:1-0:1$/2T$ and we show an example of it for $d=0.2719$ in \cref{fig:PhasePlane_omega=2pi_muk=0p2_d=0p2719}, noting that it exhibits sliding during the creation of the loop in the $Z-\dot{Z}$ plane. For $d<0.22736$ the VI-EH system yields a 2:1$/T$ stable solution. The maps were employed here to further guide numerical explorations, as they confirmed the existence of the 1:1-type solution in the larger range of $d$ values. The analytical method provides a great computational advantage, as we can evaluate the nonlinear maps at a large set of initial conditions quickly, and obtain insights about potential bistability regions.

    Numerical simulations perfectly agree with the analytically obtained 1:1-type stable solutions depicted by thick colored circles in \cref{fig:BifDiagram_omega=2pi_muk=0,fig:BifDiagram_omega=2pi_muk=0p2}. However, generating bifurcation diagrams numerically may be computationally expensive, as various conditions related to the direction of the flow have to be checked to detect crossing and sliding on $\Sigma$ at each intersection. Our continuation method also contributes to the duration of computations. In fact, in regions of co-existing solutions, continuation has to be implemented in both the increasing and decreasing direction of the bifurcation parameter $d$ from carefully selected initial conditions (see e.g. \cref{fig:ParSet2_BifDiagrams,fig:BifDiagram_omega=2pi_muk=0p2}). In some cases, we found that the semi-analytically obtained solutions provide an alternative and accelerated way for completing the bifurcation diagram and gaining insight about important bifurcations.  For example, in \cref{fig:BifDiagram_omega=2pi_muk=0p2}, the impact velocities corresponding to 1:1 motion for $d\in(0.393,0.4413)$ have been obtained via our semi-analytical maps for the 1:1 solution and then confirmed numerically via a continuation method with increasing $s$ (increasing $d$). The rightmost end of this stable branch corresponds to a fold bifurcation, revealed by linear stability analysis, through which the unstable and stable branches coalesce and disappear as $d$ is increased. The maps help us obtain the stable as well as the unstable branches (not shown here) and the fold bifurcation designates the onset of the 1:1 solution.
\section{Numerical investigation of energy harvesting}\label{sec:Energy} 
In this section we analyze how the addition of dry friction in the system influences the energy output. We calculate two measures of energy, namely the average value of generated voltage per impact, $\overline{U}_I$, and the average value of generated voltage per unit of time using the formulas 
\begin{align}
   \displaystyle \overline{U}_{I} &= \dfrac{\sum_{k=1}^{N} U_k}{N}\,, \qquad
    \overline{U}_{T} = \dfrac{\sum_{k=1}^{N} U_k}{t_f-t_0} 
\end{align}
\noindent where $N$ is the number of impacts over the (non-dimensional) interval $[t_0,t_f]$ (here, we consider the last $30$ forcing periods $T$ of our simulation), $U_k=U^{imp}_k-U_{in}$ is the net output voltage at the $k^{\textrm{th}}$ impact, $U_{in} = 2000$ mV is a constant input voltage
applied to the membranes, and $U^{imp}_k$ is the output voltage
across the deformed dielectric elastomer membrane
at the $k^{\textrm{th}}$ impact. Detailed formulas and parameter values for the calculation of $U^{imp}_k$ are given in \cref{apx:EnergyCalculations}.
\subsection{Larger $\beta$, $r$, $\omega$} \label{sec:Energy_LargePars}
\begin{figure}[hbtp!]
    \centering
    \begin{subfigure}[b]{0.45\textwidth}
    \centering
         \includegraphics[scale=0.32]{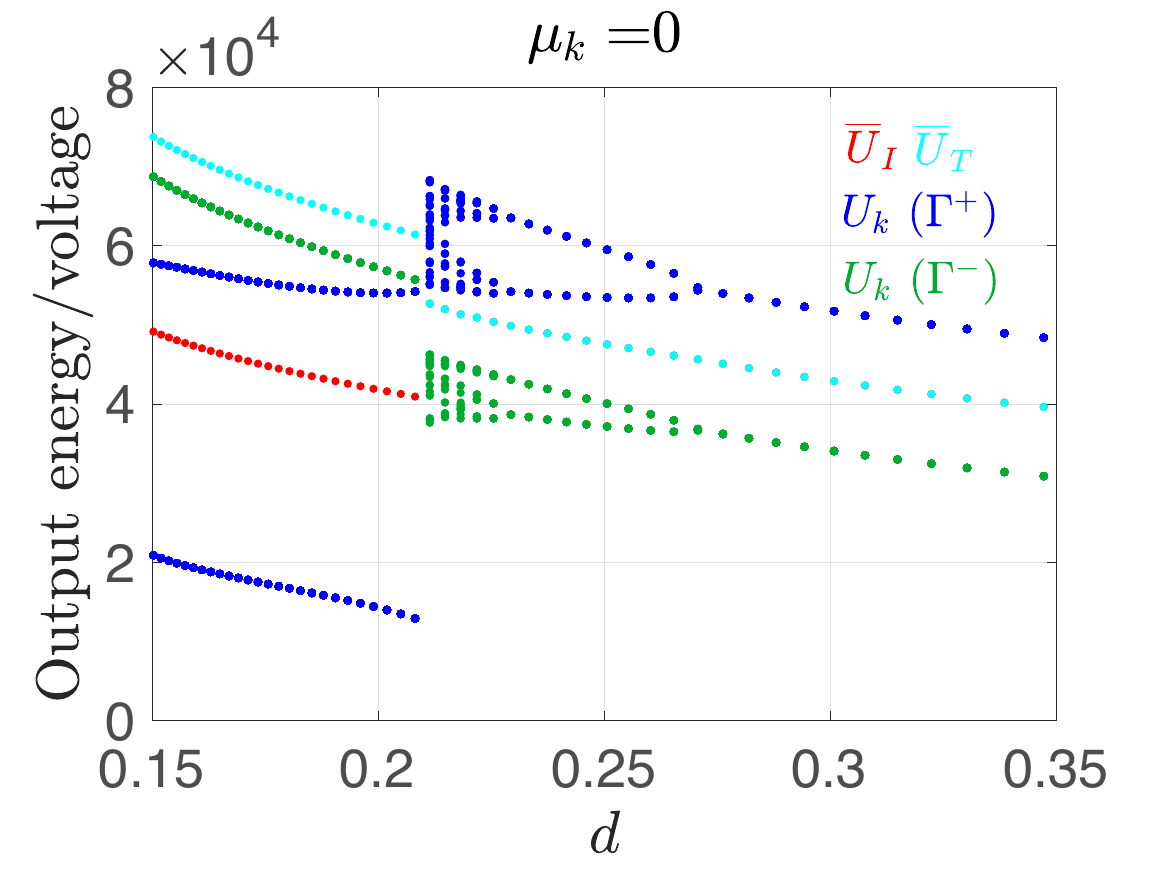}
    \caption{}
    \label{fig:ParSet1_EnergyPlot_muk=0}
    \end{subfigure}
    \hfill
    \begin{subfigure}[b]{0.45\textwidth}
    \centering
         \includegraphics[scale=0.32]{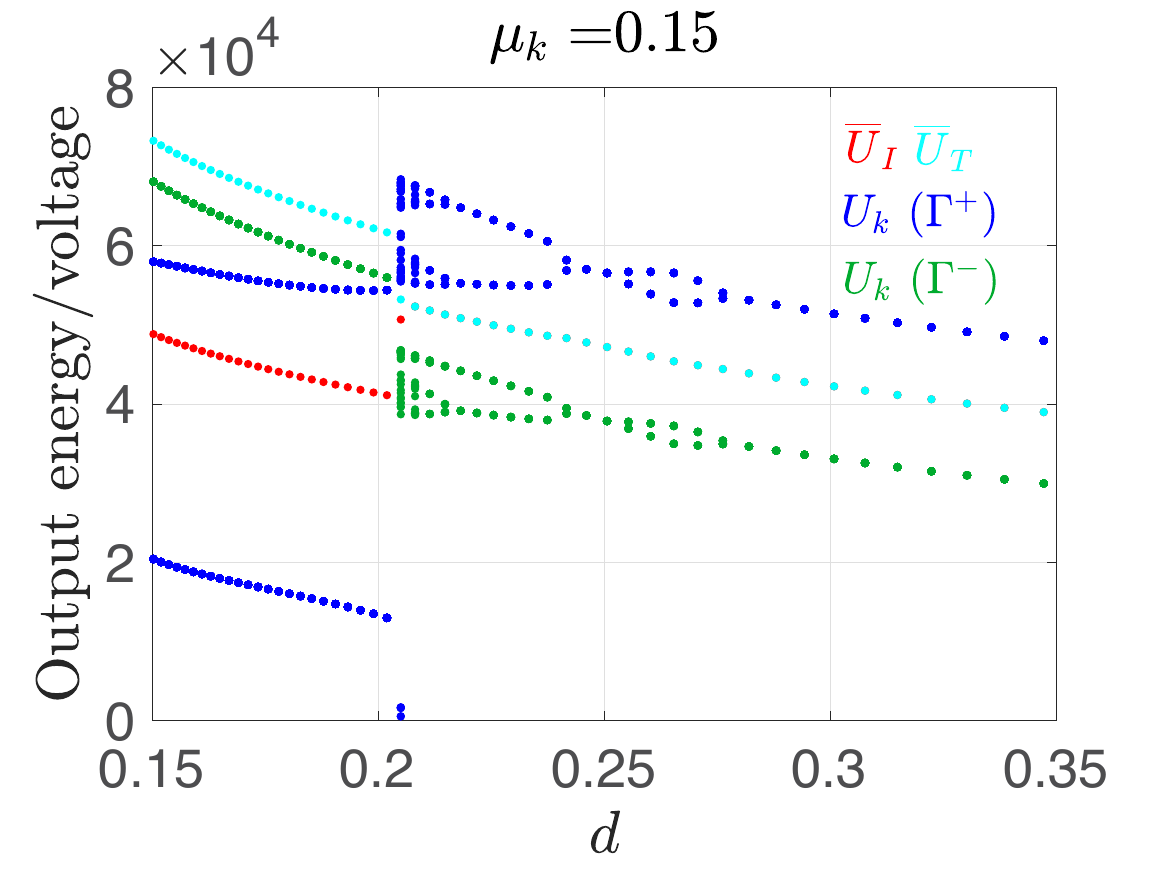}
    \caption{}
    \label{fig:ParSet1_EnergyPlot_muk=0p15}
    \end{subfigure}
    \hfill
   \begin{subfigure}[b]{0.45\textwidth}
   \centering
         \includegraphics[scale=0.32]{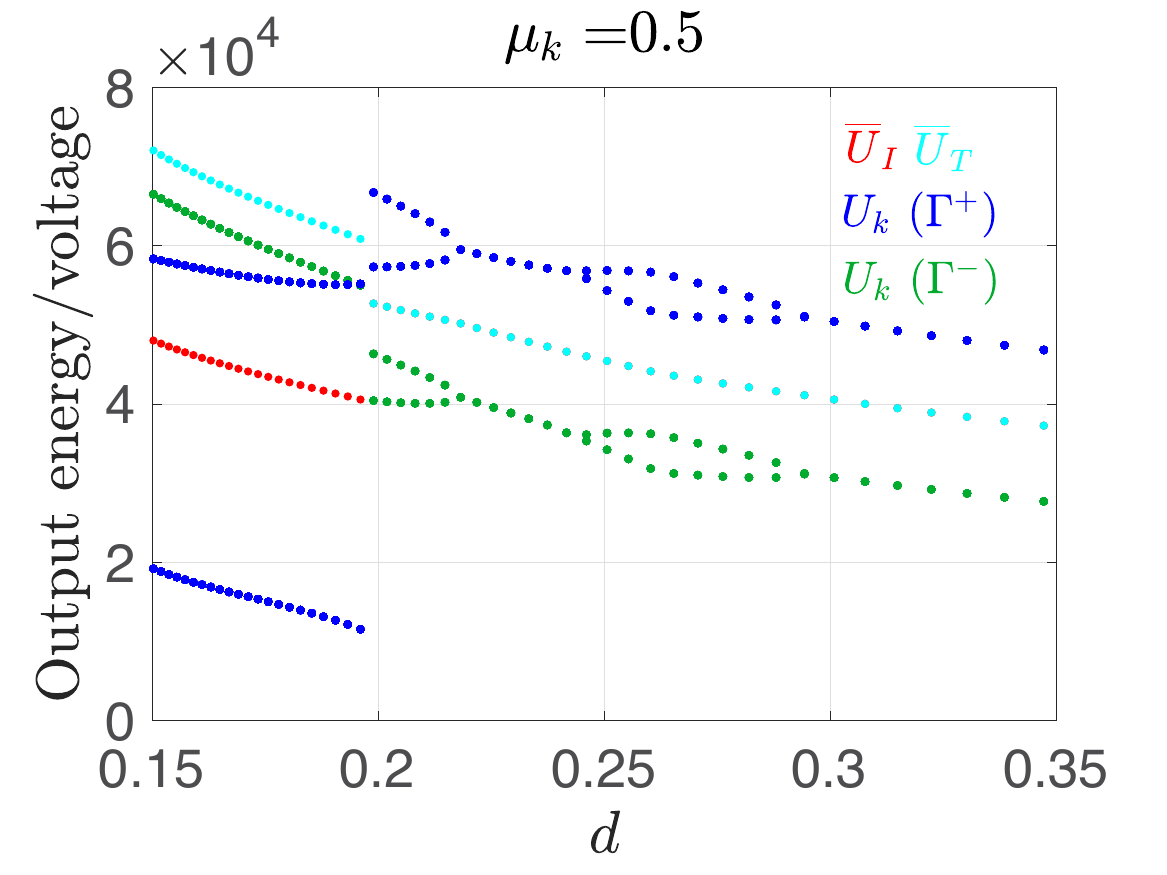}
    \caption{}
    \label{fig:ParSet1_EnergyPlot_muk=0p5}
    \end{subfigure}\hfill
    \begin{subfigure}[b]{0.45\textwidth}
    \centering
         \includegraphics[scale=0.285]{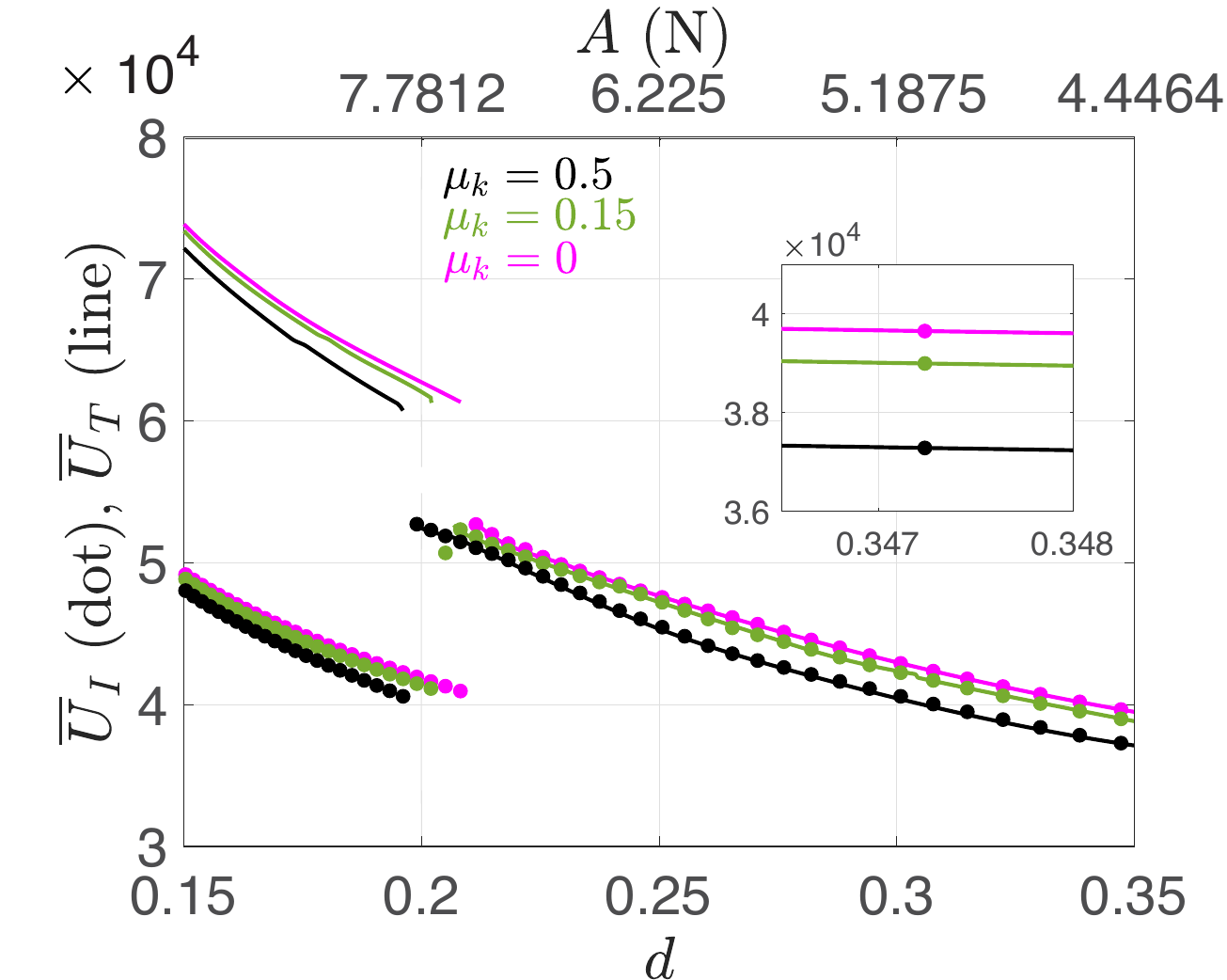}
    \caption{}
    \label{fig:AveVoltageCompA}
    \end{subfigure}
    \caption{(a)-(c): Plots of the voltage $U_k$ generated due to an impact at $\Gamma^{+}$ (blue) and $\Gamma^{-}$ (green), and of the average output voltage per impact and per unit of (non-dimensional) time, $\overline{U}_{I}$ (red) and $\overline{U}_{T}$ (cyan), respectively, for decreasing $d$ (increasing $A$) using the following parameters: $\beta = \pi/4$, $r=0.5$, $\omega =5\pi$, $s=0.5$, $A\in(4.45,14.5]$ (as in \cref{fig:ZdotBifDiagrams_FrictionComp}). The different panels show these four quantities at various values of the dry friction coefficient, namely (a) $\mu_k=0$, (b) $\mu_k=0.15$, and (c) $\mu_k=0.5$. (d) Bifurcation diagram of the average output voltage per impact, $\overline{U}_{I}$ (dots), and average output voltage per time, $\overline{U}_{T}$ (lines), for decreasing $d$, for $\mu_k=0$ (magenta), $\mu_k=0.15$ (green), and $\mu_k=0.5$ (black). Note that $\overline{U}_{I}$ and $\overline{U}_{T}$ overlap in the 1:1-type-solution regime. In all cases, the abrupt drop (increase) in $\overline{U}_{I}$ ($\overline{U}_{T}$) for $d\approx 0.2$ is due to the grazing bifurcation on $\Gamma^{+}$ that introduces 2:1 behavior.}
\label{fig:EnergyComparison}
\end{figure}

\begin{figure}[hbtp!]
     \begin{subfigure}[b]{0.5\textwidth}
    \centering
         \includegraphics[scale=0.3]{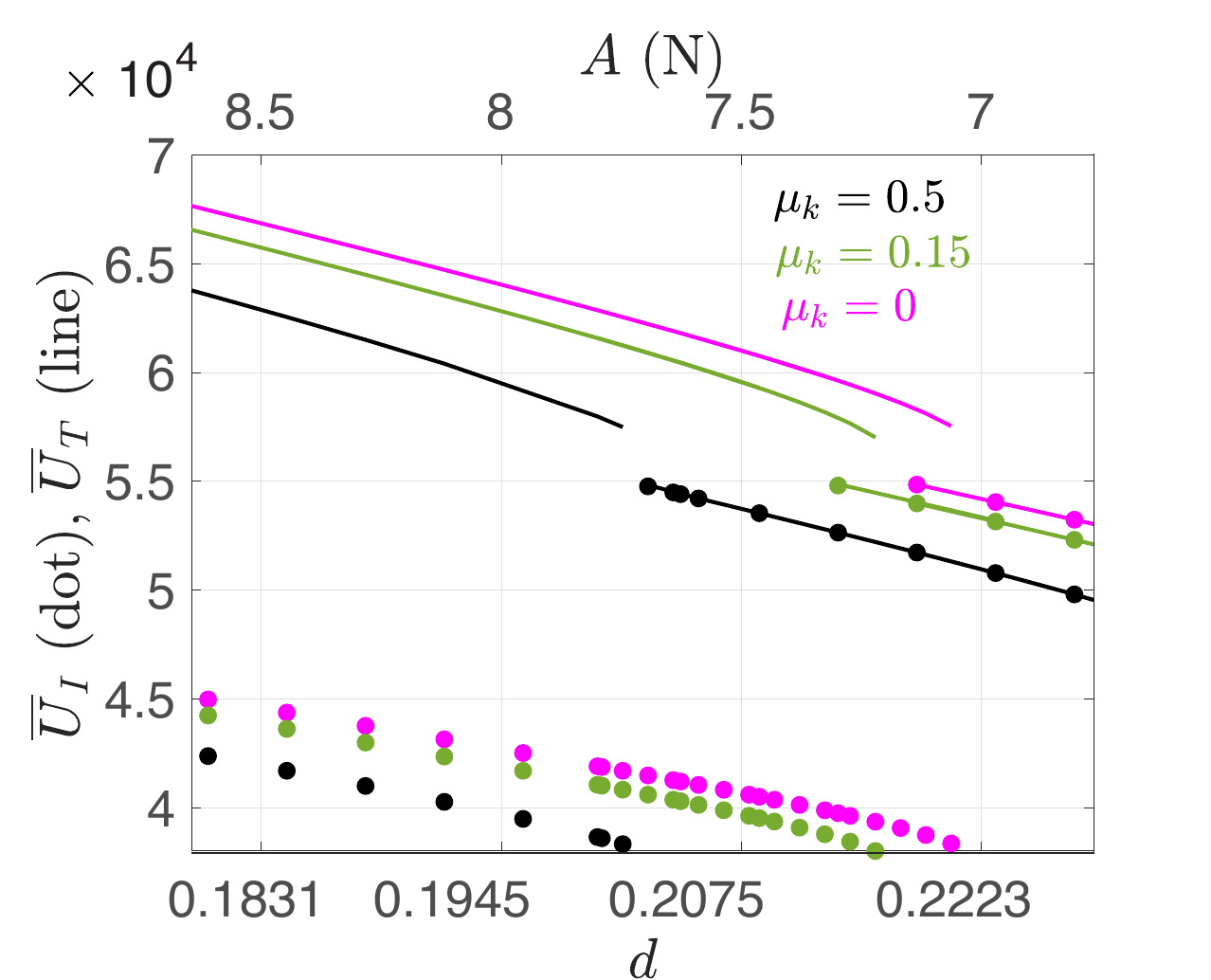}
    \caption{}
    \label{fig:ParSet2_AveVoltageCompA}
    \end{subfigure}\hfill
     \begin{subfigure}[b]{.5\textwidth}
    \centering
         \includegraphics[scale=0.3]{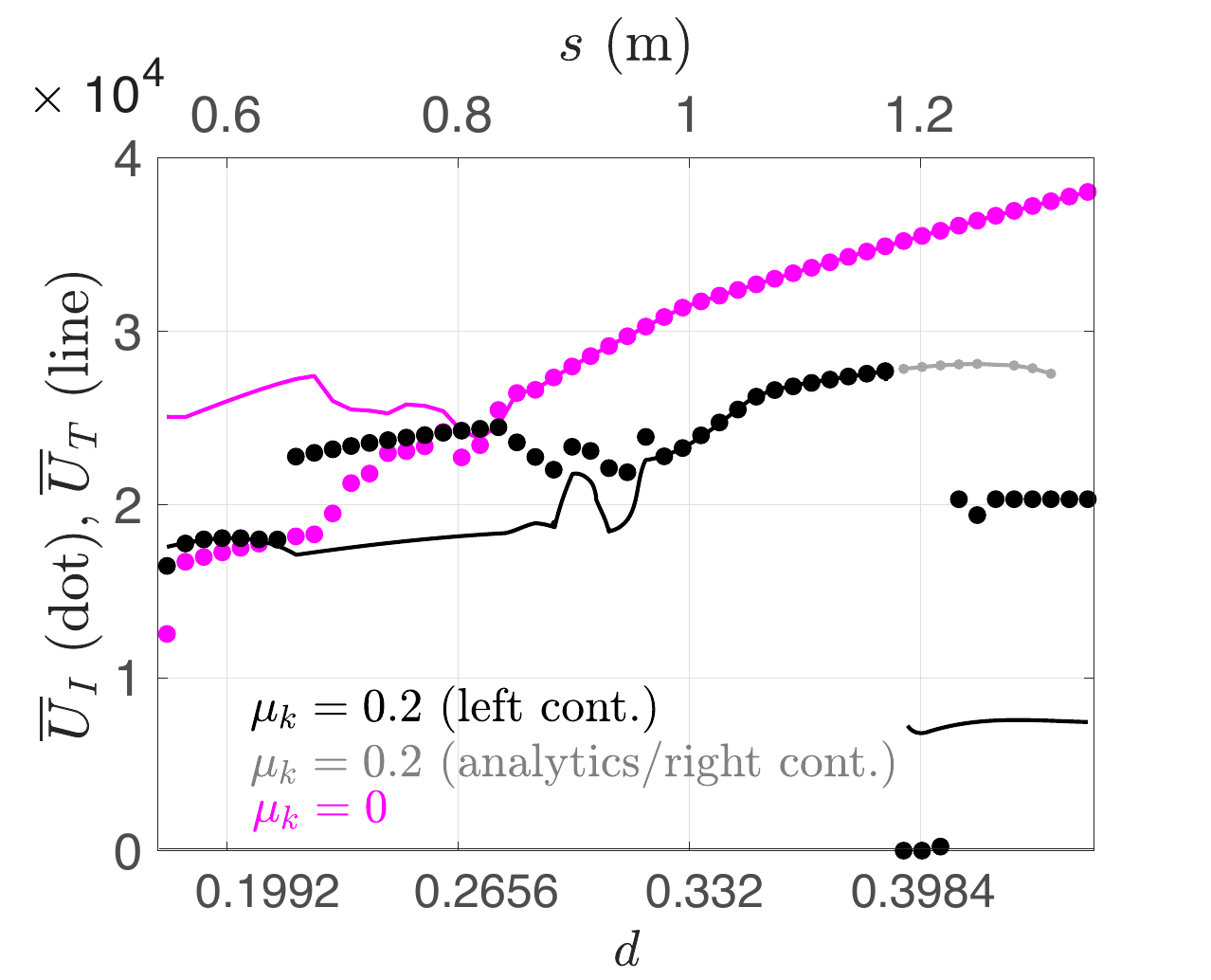}
    \caption{}
    \label{fig:ParSetOmega2pi_AveVoltageComp}
    \end{subfigure}
    \caption{Bifurcation diagrams of the average output voltage per impact, $\overline{U}_{I}$ (dots), and average output voltage per time, $\overline{U}_{T}$ (diamonds), for decreasing $d$ (increasing $A$ or decreasing $s$), for various values of $\mu_k$ and parameters (a): $\omega=5\pi$, $\beta=\pi/6$, $r=0.25$, $A\in(5.986,8.6113)$, $s=0.5$,  (as in \cref{fig:ParSet2_BifDiagrams} (a)-(d)), and (b): $\omega=2\pi$, $\beta=\pi/12$, $r=0.3$, $A=1.5$, $s\in(0.5,1.35)$ (as in \cref{fig:BifDiagrams_omega=2pi}). The different colors correspond to: $\mu_k=0$ (magenta),  $\mu_k=0.15$ (green), and (a): $\mu_k=0.5$ and (b): $\mu_k=0.2$ (black). In panel (b) the gray dots and line indicate the energy outputs for $\mu_k=0.2$ obtained for increasing $d$ (increasing $s$). As $\mu_k$ increases, the energy produced from a 1:1/$pT$ solution decreases. Panel (a) shows that
    higher friction ($\mu_k=0.5$, black markers) maintains the stability of the 1:1$^{\Sigma}_{c}/{pT}$ solution for $d>0.198$, and $d>0.2$ respectively, in contrast to the case for $\mu_k=0$ (magenta markers), with 2:1 solutions and lower (larger) $\overline{U}_{I}$ ($\overline{U}_{T}$) for the same $d$. }
\label{fig:AverageEnergyComparison}
\end{figure}
The bifurcation diagrams of the voltage, $U_k$, generated by an impact on $\Gamma^{+}$ (bottom membrane) and $\Gamma^{-}$ (top membrane) (blue and green dots in \cref{fig:ParSet1_EnergyPlot_muk=0,fig:ParSet1_EnergyPlot_muk=0p15,fig:ParSet1_EnergyPlot_muk=0p5}) are structurally similar to the bifurcation diagrams of the impact velocities (\cref{fig:ZdotBifDiagrams_FrictionComp}). For all values of the parameter $\mu_k$, the average output voltage $\overline{U}_{I}$ ($\overline{U}_{T}$) increases for decreasing dimensionless length $d$  (i.e. increasing $A$) corresponding to 1:1-type motion. We observe an abrupt decrease (increase) in $\overline{U}_{I}$ ($\overline{U}_{T}$) due to a grazing bifurcation $d_{\Gamma^{+}}$ \eqref{eq:GrazingCondition} from a stable 1:1-type motion to a stable 2:1 periodic motion. The energy produced by 2:1 motion also increases as $d$ decreases ($A$ increases) until additional grazing events yield more complex dynamics. Noting that 2:1 periodic motions include more impacts per period, this regime is typically less robust and small perturbations in $A$, and subsequently $d$, can yield large changes in the harvested energy. Then the robustness of the system in the 1:1 regime can be desirable to achieve consistently high levels of harvested energy, obtained for lower forcing of the capsule.

A comparison of \cref{fig:AveVoltageCompA}
 to \cref{fig:ZdotBifDiagrams_FrictionComp} indicates that $\overline{U}_{I}$ and $\overline{U}_{T}$ decrease as $\mu_k$ increases, away from grazing bifurcations on $\Gamma^{+}$ (e.g. to 2:1 motion). This decrease follows from the friction in the VI-EH device increasing the damping, thus reducing the impact velocity, $\dot{Z}_k$, and output voltage,
$U_k$. It is worth noting that the energy drops at most by about 7$\%$ when $\mu_k=0.5$ compared to $\mu_k=0$. However, for a small window of $d$-values, namely $d\in[0.198, 0.208]$, larger $\mu_k$ yields a stable 1:1-type motion, in contrast to smaller $\mu_k$ that yields stable 2:1 motion with lower average output energy, $\overline{U}_I$, compared to the 1:1-type behavior. Thus, near the grazing bifurcation larger $\mu_k$ may yield larger average energy output per impact. In fact, we find that $\overline{U}_I$ is 31$\%$ larger for $\mu_k=0.5$ compared to $\mu_k=0$ for the same $d$ value in this range. Moreover, note that for larger $\mu_k$ the maximum energy output due to a 1:1-type motion (i.e. near the grazing bifurcation) and the relative energy decrease across the grazing bifurcation are comparable to the maximum energy output obtained from a 1:1-type for smaller $\mu_k$ (in the corresponding grazing bifurcation regimes). On the other hand, $\overline{U}_T$ increases near the grazing bifurcation, with lower value of friction yielding larger energy outputs.
\subsection{Smaller $\beta$, $r$, larger $\omega$} \label{sec:Energy_MixedPars}
Similar results are observed in other parameter regimes as shown in \cref{fig:ParSet2_AveVoltageCompA}, where smaller values for $\beta$ and $r$ are incorporated. In particular, $\beta=\pi/6$, $r=0.25$, while $\omega=5\pi$, $A\in(5.986,8.6113)$, and $s=0.5$. Again, larger $\mu_k$ reduces $\dot{Z}_k$ and 
$U_k$ away from the grazing transition to 2:1 motion. Near this bifurcation, larger $\mu_k$ ($\mu_k=0.5$) can shift the grazing to smaller values of $d$, maintaining a stable 1:1$^{\Sigma}_{c}$ motion into the interval $d\in[0.204, 0.213]$. Then, in this interval, larger $\mu_k$ yields a  larger $\overline{U}_I$ as compared with smaller $\mu_k$. In this way, friction can stabilize the regular 1:1-type dynamics and extend the range of higher energy outputs produced by them.
\subsection{Smaller $\beta$, $r$, $\omega$} \label{sec:Energy_SmallerPars}
In \cref{fig:ParSetOmega2pi_AveVoltageComp}, the two energy measures are plotted in terms of $d$ for the parameter set of \cref{fig:BifDiagrams_omega=2pi}, namely $\omega=2\pi$, $\beta=\pi/12$, $r=0.3$, $A=1.5$, $s\in(0.5,1.35)$, for $\mu_k=0$ (magenta) and $\mu_k=0.2$ (black/gray). Note that in contrast to the cases discussed in \cref{sec:Energy_LargePars} and \cref{sec:Energy_MixedPars}, $\omega$, $\beta$, and $r$ are all smaller. Moreover, as $s$ decreases, $d$ also decreases. As shown in \cref{fig:BifDiagram_omega=2pi_muk=0}, for higher $d$ the frictionless system produces 1:1-type behavior and $\overline{U}_{I}=\overline{U}_{T}$, which decrease as $d$ decreases. For $\mu_k=0.2$, the energy produced in the same $d$-range (where both 1:1-type and bouncing-ball-type dynamics may be present) is lower. As discussed in \cref{sec:ComparisonAndStabAnalysis} the 1:1-type solution ceases to exist as $d$ decreases, and chaotic dynamics are present for $\mu_k=0$ (\cref{fig:BifDiagram_omega=2pi_muk=0}) and $\mu_k=0.2$ (\cref{fig:BifDiagram_omega=2pi_muk=0p2}). The chaotic solution in the frictionless case exhibits small impact velocities $\dot{Z}_k$ near $0$, and as $d$ decreases it is followed by a stable 2:1 solution for $d<0.229$. However, for $\mu_k=0.2$ the impact velocities produced in the chaotic regime are bounded away from $0$, and the chaotic dynamics are succeeded by a stable 1:1-0:1$/2T$ solution for $d < 0.282$. Note that in both cases, there is a separation between $\overline{U}_{I}$ and $\overline{U}_T$ when the model dynamics do not correspond to 1:1 behavior. For $\mu_k=0$, the 2:1-type chaotic or stable dynamics ($d < 0.2877$) lead to $\overline{U}_T>\overline{U}_I$ (see magenta diamonds and dots in \cref{fig:ParSetOmega2pi_AveVoltageComp}). For $d<0.33$ and $\mu_k=0.2$ the model dynamics corresponding to the 1:1-0:1$/2T$-type solution lead to a different relationship in the energy measures as $\overline{U}_{T}<\overline{U}_{I}$ (see black diamonds and dots in \cref{fig:ParSetOmega2pi_AveVoltageComp}). Note also that even though $\overline{U}_{T}$ is larger for $\mu_k=0$ compared to $\mu_k=0.2$ (magenta vs black/gray diamonds in \cref{fig:ParSetOmega2pi_AveVoltageComp}), $\overline{U}_{I}$ is consistently larger for $\mu_k=0.2$ compared to $\mu_k=0$ for $d<0.272$ (magenta vs black/gray dots in \cref{fig:ParSetOmega2pi_AveVoltageComp}). In this $d$-range, the VI-EH system yields a stable 1:1-0:1$/2T$-type solution for $\mu_k=0.2$, and a chaotic or stable 2:1-type solution for $\mu_k=0$. This suggests that, in some parameter regimes, friction may provide a mechanism for producing larger output energy (per impact) and more robust solutions in the model.\section{Conclusions}
\label{sec:Conclusions}
 We have investigated the influence of dry friction on the dynamics of a VI-EH system. The system consists of a harmonically excited capsule and a bullet that is freely moving in its interior. The ends of the capsule are covered with membranes that can deform due to impacts and act as variable capacitance capacitors, allowing for an excess in energy to be harvested. Viewed as a dynamical system, the VI-EH system is non-smooth due to the impacts at the bottom (top) of the capsule, $\Gamma^{+(-)}$. Dry friction introduces a second type of non-smoothness, as it defines an additional state-dependent, time-varying switching boundary, $\Sigma$, corresponding to the relative velocity between the capsule and bullet vanishing, $\dot{Z}=0$. There, the dynamical system is discontinuous, and more specifically, of Filippov type, with $\dot{Z}$ changing sign on either side of $\Sigma$. With harmonic forcing, the system exhibits different types of non-smooth periodic behavior that, due to the dry friction, may include sticking intervals where the capsule and bullet move together. Mathematically, this appears as sliding, where the system follows $\Sigma$ with $\dot{Z}=0$ over that interval. 

In this paper we develop a systematic framework based on compositions of discrete nonlinear maps that are combined to capture sliding, non-stick motion, and impacts. We construct various 1:1-type periodic solutions, with alternating $\Gamma^{\pm}$ impacts per period, by solving reduced systems of equations derived from the maps. These yield the initial impact velocity, the impact phase, and the time intervals between impacts on $\Gamma^{\pm}$ and sliding or crossing events on $\Sigma$, which characterize different behaviors.  With friction included in the model, the complexity of the map compositions and 
associated solutions increases with the number of sliding or crossing events on $\Sigma$, in contrast to the frictionless case where only impacts on $\Gamma^{\pm}$ must be tracked.  The semi-analytical solutions are combined with conditions on the flow across $\Sigma$ and the impact conditions to identify sequences of bifurcations that include grazing-sliding, switching-sliding and crossing-sliding, associated with $\Sigma$, and grazing bifurcations on $\Gamma^{\pm}$. These conditions also remove unphysical solutions that solve the reduced system. These non-smooth bifurcations may be interwoven with smooth bifurcations, such as period doubling (PD), that are determined via linear stability analyses. In these analyses the calculation of the Jacobian associated with the map compositions requires careful tracking of the dependencies between the times, velocities, and displacements at the various non-smooth events. 

We focus on 1:1/$pT$ periodic solutions that are more favorable for energy harvesting \cite{DULIN2022106983,Serdukova2022FundamentalCO,Yurchenko_2017}, comparing the influence of the dry friction parameter $\mu_k$ on the  solution features and the related bifurcation structure. For example, higher values of $\mu_k$ lead to longer maximal intervals of sliding motion. Following the bifurcations in terms of decreasing non-dimensional length parameter $d$,
for certain parameter regimes including low or no friction (see \cref{fig:ZdotBifDiagrams_FrictionComp_A} and \cref{fig:ZdotBifDiagrams_FrictionComp_B}), there is a period doubling cascade to chaos for 1:1-type behavior, followed by   a grazing bifurcation on $\Gamma^{+}$ that leads to stable 2:1 motion. In contrast, for larger values of $\mu_k$, (see \cref{fig:ZdotBifDiagrams_FrictionComp_C} and \cref{fig:ZdotBifDiagrams_FrictionComp_D}), sticking yields a sequence of grazing-sliding, switching-sliding, and crossing-sliding bifurcations interleaved with PD bifurcations for  1:1-type solutions. These reduce or replace the transition to chaos and are followed by the impact-type grazing to 2:1 motion.

Throughout, the analysis captures the trend that dry friction slows the system down, reducing the range of relative velocity $\dot{Z}(t)$ and the related impact velocity $\dot{Z}_k$. While this may appear detrimental to energy output, the bifurcation structure provides a more nuanced picture. The behavior is regularized by friction in the VI-EH system, and as a consequence, it can extend the parameter range where the periodic 1:1-type branch persists, thus extending the range where higher average output voltage per impact is realized. Therefore, our results point to parameter regimes where dry friction is beneficial, as it can limit the appearance of chaos and grazing and sustain the larger energy output of 1:1 motion, yielding a relative increase in $\overline{U}_{I}$ of up to $31\%$, even though a small percent decrease due to friction is expected away from grazing (e.g up to $7\%$ in the parameter regime discussed in \cref{fig:EnergyComparison}). Thus, dry friction may provide a design or intrinsic system  mechanism for limiting undesirable and irregular model dynamics. While in some cases this beneficial range may seem small in terms of the non-dimensional parameter $d$, it might lead to large differences in the physical parameters, such as the length of the capsule (see \cref{fig:ParSetOmega2pi_AveVoltageComp}).

The framework presented in this paper is generalizable for studying the dynamics and bifurcation sequences in similar VI systems, such as the impact pair \cite{bapat1988impact,ZhangFu2017} and models of energy transfer \cite{KUMAR2024118131,LI2021104001}. Preliminary experiments \cite{PrivComm_HSDY} motivate future studies to explore the influence of friction in other parameter regimes, including different forcing frequencies and variable restitution coefficients. There, we expect to find new complex periodic behaviors, regions of bistability, and alternate routes to chaos not observed in the absence of friction, as indicated by \cref{fig:BifDiagram_omega=2pi_muk=0p2}. Our semi-analytical approach can be extended to incorporate side-stick motion as well. All of these examples can severely influence energy output. Finally, our current analysis suggests parameter regimes, e.g. near traditional and non-smooth bifurcations, in which the addition of stochastic effects may influence the system behavior, leading to interaction of model solutions with unstable orbits, and consequently the energy output.
\appendix
\section{Conditions for the occurrence of non-stick and sliding motions}
\label{apx:Conditions}
To find conditions for the occurrence of non-stick and sliding motions, we look at the phase plane $Z-\dot{Z}$ and determine how the vector fields across the different switching boundaries influence model solutions. Let $\dot{\mathbf{{Z}}} = (\dot{Z},\ddot{Z})^T$ be the vector field and $\mathbf{n}_{\Sigma}^T = (0,1)^T $ be the normal vector to the switching boundary, $\Sigma$, in relative frame. Then
\begin{align}
    \mathbf{n}_{\Sigma}^T\cdot \dot{\mathbf{Z}} = \ddot{Z} = f(t) -L^{\lambda} \, ,
\end{align}
\noindent where $f(t) = \cos(\pi t +\phi)$ and $\lambda=\pm$ for $\Sigma^{\pm}$. Note that $L^{-}<0$ and $L^{-}<L^{+}$ for $\beta\in[0,\pi/2)$ and $\mu_k\neq 0$. We also note that throughout the manuscript $L^{\pm}\in(-1,1)$.

Below we provide the key conditions for non-stick and sliding motion for trajectories that evolve from $\Sigma^{-}$ towards $\Sigma$. A more detailed justification and general list of those conditions are given after the statements. 
\begin{itemize}
    \item Non-stick motion from $\Sigma^{-}$ to $\Sigma^{+}$ at $t=t_c$: \begin{align}\label{eq:CrossingConditionsMinusToPlus_1}
 f(t_c) - L^{+}>0\Rightarrow\cos(\pi t_c +\phi) > L^{+} \\\label{eq:CrossingConditionsMinusToPlus_2}
    f(t_c) - L^{-}>0\Rightarrow \cos(\pi t_c +\phi) > L^{-}   
\end{align}
\item  Sliding motion from $\Sigma^{-}$ to $\Sigma$ and back to $\Sigma^{-}$ on $[t_{so},t_{se}]$: 
\begin{itemize}
    \item Onset of sliding motion at $t=t_{so}$:   \begin{align}\label{eq:SlidingOnset2a}
   f(t_{so})>L^{-} &\Rightarrow \cos{(\pi t_{so} + \phi)} >L^{-}\\\label{eq:SlidingOnset2b}
   f(t_{so})\leq L^{+} &\Rightarrow \cos{(\pi t_{so} + \phi)}  \leq L^{+}\\\label{eq:SlidingOnset2c}
f'(t_{so})<0 &\Rightarrow -\pi\sin{(\pi t_{so} + \phi)} <0 \Rightarrow \sin{(\pi t_{so} + \phi)} >0
\end{align}
    \item  Exit time of sliding motion at $t=t_{se}$:\begin{align}\label{eq:NonPassToSemiPassCondToOmega2_1}
    f(t_{se})=L^{-} & \Rightarrow \cos(\pi t_{se} + \phi) = L^{-}
\\\label{eq:NonPassToSemiPassCondToOmega2_2}
   f'(t_{se})<0 & \Rightarrow-\pi\sin{(\pi t_{se} + \phi)} < 0 \Rightarrow \sin{(\pi t_{se} + \phi)} > 0 \\\label{eq:NonPassToSemiPassCondToOmega2_3}
    f(t_{se})<L^{+} &\Rightarrow \cos(\pi t_{se} + \phi) < L^{+}
\end{align}
\end{itemize}
\end{itemize}
\subsection{Non-stick motion}\label{sec:NonStickMotionCondsApp}
During the non-stick motion, the vector fields in $\Sigma^{+}$ and $\Sigma^{-}$ point in the same direction and therefore, allow model solutions to cross $\Sigma$. Let $t=t_c$ be the time at which a model solution crosses $\Sigma$. Then, the following condition holds:
\begin{align}
    (f(t_c) - L^{+})(f(t_c) - L^{-})&> 0 \,,
\end{align}
implying that either $f(t_c)-L^{\pm}>0$ or $f(t_c)-L^{\pm}<0$.

For the flow to pass from $\Sigma^{+}$ to $\Sigma^{-}$ at $t=t_c$ we require that:
\begin{align}\label{eq:CrossingConditionsPlusToMinus_1}
   f(t_c) - L^{+}<0\Rightarrow\cos(\pi t_c +\phi) < L^{+} \\\label{eq:CrossingConditionsPlusToMinus_2}
    f(t_c)- L^{-}<0\Rightarrow \cos(\pi t_c +\phi) < L^{-} 
\end{align}
For these to be satisfied simultaneously:
\begin{itemize}
    \item if $L^{-}\geq -1$, then $\mod{(\pi t_c+\phi,2\pi)}\in[\arccos{(L^{-})},2\pi-\arccos{(L^{-})}]$, 
    \item if $L^{-}<-1$ then there is no solution, namely the flow does not pass from $\Sigma^{+}$ to $\Sigma^{-}$.
\end{itemize}
For the flow to pass from $\Sigma^{-}$ to $\Sigma^{+}$ at $t=t_c$ we require that conditions \cref{eq:CrossingConditionsMinusToPlus_1,eq:CrossingConditionsMinusToPlus_2} are met, which are repeated below for the reader's convenience:
\begin{align*}
 f(t_c) - L^{+}>0\Rightarrow\cos(\pi t_c +\phi) > L^{+} \\
    f(t_c) - L^{-}>0\Rightarrow \cos(\pi t_c +\phi) > L^{-}   
\end{align*}
For these to be satisfied simultaneously:
\begin{itemize}
    \item if $L^{+}\in(-1,1)$ (regardless of $L^{-}$), then $\mod{(\pi t_c+\phi,2\pi)}\in[0,\arccos{(L^{+})})\cup(2\pi-\arccos{(L^{+})},2\pi]$,
    \item if $L^{+}\leq -1$ then $\mod{(\pi t_c+\phi,2\pi)}\in[0,2\pi]$, 
    \item if $L^{+}>1$ then there is no solution, namely the flow does not pass from $\Sigma^{-}$ to $\Sigma^{+}$.
\end{itemize}

\subsection{Sliding motion}\label{sec:SlidingStick}
Sliding (or sliding-stick) motion takes place on an interval $[t_{so}, t_{se}]$ when
\begin{align}\label{eq:GeneralSlidingCondition}
    (f(t) - L^{+})(f(t) - L^{-})&\leq 0 , \text{ for } t\in[t_{so}, t_{se}]\,.
\end{align}
This condition indicates that the vector fields across $\Sigma$ both point towards or away from $\Sigma$ during the sliding. Then, we can define the Filippov sliding vector field, $\dot{\mathbf{{Z}}_{s}}$, to be a convex combination of the vector fields in $\Sigma^{\pm}$, such that 
\begin{align}\label{eq:GeneralSlidingCondition}
    \mathbf{n}_{\Sigma}^T\cdot \dot{\mathbf{Z}}_{s} &= 0 \,.
\end{align}
Therefore, during the sliding motion, the capsule and the bullet move together, the relative displacement $Z$ is constant, and the relative velocity $\dot{Z}=0$ (see \cref{fig:SlidingStick_Timetrace_AbsFrame,fig:SlidingStick_Timetrace,fig:SlidingStick_PhasePlane}).\vspace{0.15cm}

\noindent\underline{\textit{Explicit conditions for the existence of sliding motion}} 
Condition \cref{eq:GeneralSlidingCondition} implies that for $t\in[t_{so},t_{se})$ either 
\begin{align*}
    f(t)-L^{+}<0 \text{ and } f(t)-L^{-}>0 \text{ or }\\
    f(t)-L^{+}>0 \text{ and } f(t)-L^{-}<0
\end{align*}
However, since $L^{-}<L^{+}$ only the first set of inequalities can be satisfied. Given that the forcing function is $f(t) = \cos(\pi t+ \phi)$ we obtain the following for $t$ in the sliding interval $(t_{so},t_{se})$,
\begin{itemize}
    \item if $L^{+},L^{-}\in(-1,1)$ then $$\mod{(\pi t +\phi,2\pi)}\in(\arccos{(L^{+})},\arccos{(L^{-})})\cup(2\pi-\arccos{(L^{-})},2\pi-\arccos{(L^{+})})$$
    \item if $L^{-}\leq -1$ and $L^{+}\in(-1,1)$ then $\mod{(\pi t+\phi,2\pi)}\in(\arccos{(L^{+})},2\pi-\arccos{(L^{+})})$
    \item if $L^{+}\geq 1$ and $L^{+}\in(-1,1)$ then $$\mod{(\pi t +\phi,2\pi)}\in(0,\arccos{(L^{-})})\cup(2\pi-\arccos{(L^{-})},2\pi)$$
    \item if $L^{-}\leq -1$ and $L^{+}\geq 1$ then $\mod{(\pi t +\phi,2\pi)}\in [0,2\pi]\,.$
\end{itemize}
Below we provide conditions for the onset and exit time of sliding for $L^{\pm}\in(-1,1)$.\vspace{0.15cm}

\begin{figure}[hbtp!]
    \centering
    \includegraphics[scale=0.42]{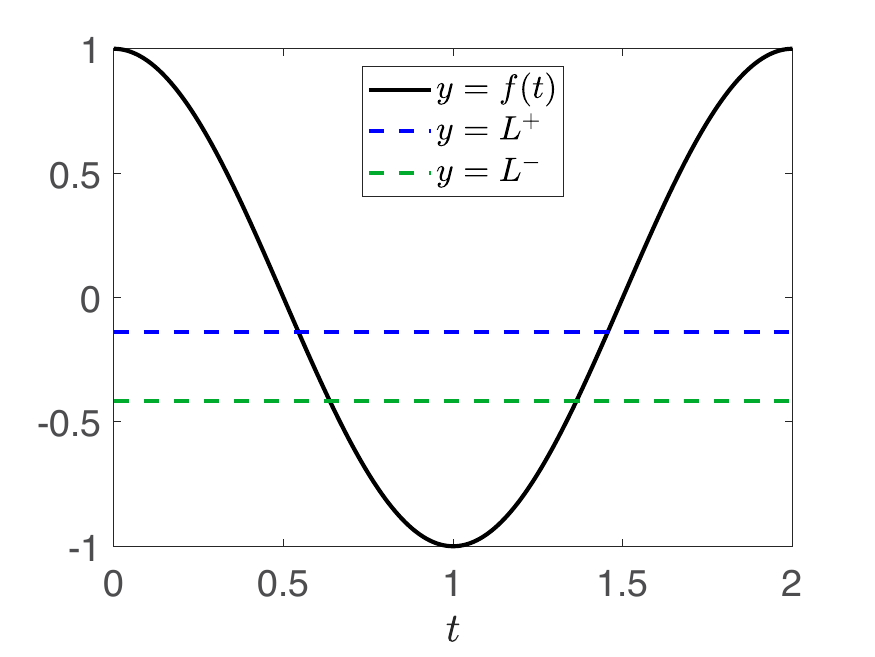}
    \caption{Time series of the forcing term $f(t)=\cos(\pi t+\phi)$ (black curve), the line $y=L^{+}$ (blue curve), and the line $y=L^{-}$ (green curve) such that $L^{\pm}\in(-1,1)$.}
    \label{fig:ForcingTimetrace}
\end{figure}

\noindent\underline{\textit{Onset of sliding motion}} 
For $L^{\pm}\in(-1,1)$, the inequality in \eqref{eq:GeneralSlidingCondition} corresponding to sliding holds for $t$ such that \begin{align}\label{eq:PotentialSlidingIntervals}
 \text{mod}{(\pi t +\phi,2\pi)}\in(\arccos{(L^{+})},\arccos{(L^{-})})\cup(2\pi-\arccos{(L^{-})},2\pi-\arccos{(L^{+})})\,.
\end{align}
The onset time of sliding, $t_{so}$, can be any point in the intervals in \eqref{eq:PotentialSlidingIntervals}, with the maximal duration of sliding being $\Delta S = (\arccos{(L^{-})}-\arccos{(L^{+})})/\pi$.  Then the onset time, $t_{so}$ satisfies conditions \eqref{eq:SlidingOnset2a}, the equality in \eqref{eq:SlidingOnset2b}, and \eqref{eq:SlidingOnset2c}, or \eqref{eq:SlidingOnset1a}, the equality in \eqref{eq:SlidingOnset1b}, and \eqref{eq:SlidingOnset1c}. 
\begin{itemize}    
    \item Onset of sliding motion from $\Sigma^{-}$ to $\Sigma$: Our choice of $\beta>0$ promotes the occurrence of sliding motion on $\Sigma$, when trajectories transition from $\Gamma^{+}$ to $\Gamma^{-}$, and therefore, approach $\Sigma$ from $\Sigma^{-}$. The sliding motion starts at $t=t_{so}$ when conditions \cref{eq:SlidingOnset2a,eq:SlidingOnset2b,eq:SlidingOnset2c} are met. Here, we repeat them  for the reader's convenience:
    \begin{align*}
   f(t_{so})>L^{-} &\Rightarrow \cos{(\pi t_{so} + \phi)} >L^{-}\\
   f(t_{so})\leq L^{+} &\Rightarrow \cos{(\pi t_{so} + \phi)} \leq L^{+}\\
f'(t_{so})<0&\Rightarrow -\pi\sin{(\pi t_{so} + \phi)} <0 \Rightarrow \sin{(\pi t_{so} + \phi)} >0
\end{align*}
 Using \eqref{eq:SlidingOnset2c}, \eqref{eq:SlidingOnset2a}  gives an upper bound on $\mod{(\pi t_{so}+\phi,2\pi)}$ and \eqref{eq:SlidingOnset2b} gives a lower bound. Subsequently, $\mod{(\pi t_{so}+\phi,2\pi)}\in[\arccos(L^{+}),\arccos(L^{-}))$.

To better understand the conditions above, we consider the behavior of the vector fields before and during the sliding motion. A trajectory evolves from $\Sigma^{-}$ to $\Sigma$ and therefore, $\ddot{Z}_{\Sigma^{-}}=f(t)-L^{-}>0$ for $t<t_{so}$ (see \cref{fig:ForcingTimetrace}). 
Therefore, the earliest onset time for sliding satisfies $\ddot{Z}_{\Sigma^{+}} = f(t_{so})-L^{+}=0$, when the trajectory reaches $\Sigma$, so that $\ddot{Z}_{\Sigma^{+}}$ changes sign at $t=t_{so}$.  Therefore, for $t<t_{so}$, $\ddot{Z}_{\Sigma^{+}}>0$, for $t=t_{so}$ $\ddot{Z}_{\Sigma^{+}}=0$, and for $t>t_{so}$, namely after the sliding onset, $\ddot{Z}_{\Sigma^{+}}<0$ while $\ddot{Z}_{\Sigma^{-}}$ is still positive. Since $L^{+}>L^{-}$, this occurs during the decreasing phase of $f(t)$, and thus, $f'(t_{so})<0$ (or $\dddot{Z}_{\Sigma^{+}}$<0 at $t=t_{so}$). This implies that the maximal sliding interval is $
    \{t\ |\mod{(\pi t +\phi,2\pi)}\in[\arccos{(L^{+})},\arccos{(L^{-})}]\}\,.$
    \end{itemize}
     \vspace{-0.5cm}
\begin{remark}
 Note that a trajectory that reaches $\Sigma$ from $\Sigma^{-}$ at $t=t_{c}$ such that $f(t_{c})>L^{+}>L^{-}$ and $f'(t_c)<0$ will first pass to $\Sigma^{+}$, since conditions \cref{eq:CrossingConditionsMinusToPlus_1,eq:CrossingConditionsMinusToPlus_2} for non-stick motion hold. However, since $f'(t)<0$, the flow will approach $\Sigma$ again as it evolves in $\Sigma^{+}$. In this case, the model solution may reach $\Sigma$ again at $t=t_{so}$, such that $\mod{(\pi t_{so} +\phi,2\pi)}\in(\arccos{(L^{+})},\arccos{(L^{-})})$, and exhibit sliding with duration $\Delta S\leq (\arccos{(L^{-})}-\arccos{(L^{+})})/\pi$. This results in a sequence of non-stick and sliding motions as described in 1:1$^{\Sigma}_{cs}$ periodic solutions.    
\end{remark}
\begin{itemize}
  \item Onset of sliding motion from $\Sigma^{+}$ to $\Sigma$:
    \begin{align}\label{eq:SlidingOnset1a}
    f(t_{so})<L^{+} &\Rightarrow \cos(\pi t_{so}+\phi)<L^{+} \\\label{eq:SlidingOnset1b}
 f(t_{so})\geq L^{-} & \Rightarrow \cos(\pi t_{so}+\phi)\geq L^{-}\\\label{eq:SlidingOnset1c}
f'(t_{so}) >0 &\Rightarrow  -\pi\sin{(\pi t_{so} + \phi)} >0\Rightarrow \sin{(\pi t_{so} + \phi)} <0
\end{align}
 Using \eqref{eq:SlidingOnset1c}, \eqref{eq:SlidingOnset1a}  gives an upper bound on $\mod{(\pi t_{so}+\phi,2\pi)}$ and \eqref{eq:SlidingOnset1b} gives a lower bound. Subsequently, $\mod{(\pi t_{so}+\phi,2\pi)}\in[2\pi-\arccos(L^{-}), 2\pi-\arccos(L^{+}))$.\vspace{0.15cm}
\end{itemize}
\noindent\underline{\textit{End of sliding motion}} 
The sliding motion ends at $t=t_{se}$, when the flows in $\Sigma^{\pm}$ allow crossing of the switching boundary, $\Sigma$. From the previous analysis, it follows that the necessary conditions are \eqref{eq:NonPassToSemiPassCondToOmega2_1}-\eqref{eq:NonPassToSemiPassCondToOmega2_3}, which are repeated below for the reader's convenience:
\begin{itemize}
    \item From the boundary $\Sigma$ to $\Sigma^{-}$:
    \begin{align*}
    f(t_{se})=L^{-} & \Rightarrow \cos(\pi t_{se} + \phi) = L^{-}
    \\
   f'(t_{se})<0 & \Rightarrow-\pi\sin{(\pi t_{se} + \phi)} < 0 \Rightarrow \sin{(\pi t_{se} + \phi)} > 0 \\
    f(t_{se})<L^{+} &\Rightarrow \cos(\pi t_{se} + \phi) < L^{+}
\end{align*}

Using \eqref{eq:NonPassToSemiPassCondToOmega2_2}, \eqref{eq:NonPassToSemiPassCondToOmega2_3}  gives an upper bound on $\mod{(\pi t_{se}+\phi,2\pi)}$ and \eqref{eq:NonPassToSemiPassCondToOmega2_1} gives an exact value.

When sliding occurs during the upward motion of the system from $\Gamma^{+}$ to $\Gamma^{-}$, $L^{-}<f(t)<L^{+}$, while $f'(t)<0$. Therefore, the end of sliding occurs at $t=t_{se}$, such that $f(t_{se})-L^{-}=0$. For $t>t_{se}$, $f(t)<L^{-}<L^{+}$, and thus, the model solution enters $\Sigma^{-}$, since conditions \cref{eq:CrossingConditionsPlusToMinus_1,eq:CrossingConditionsPlusToMinus_2} for the occurrence of non-stick motion hold.
    \item From the boundary $\Sigma$ to $\Sigma^{+}$:\begin{align}\label{eq:NonPassToSemiPassCondToOmega1_1}
f(t_{se})=L^{+} &\Rightarrow \cos(\pi t_{se} + \phi) = L^{+}\\\label{eq:NonPassToSemiPassCondToOmega1_2}
  f'(t_{se})>0\Rightarrow  -\pi\sin{(\pi t_{se} + \phi)} &>0 \Rightarrow \sin{(\pi t_{se} + \phi)} <0 \\\label{eq:NonPassToSemiPassCondToOmega1_3}
  f(t_{se})>L^{-} &\Rightarrow \cos(\pi t_{se} + \phi) > L^{-}
\end{align}
Similarly, using \eqref{eq:NonPassToSemiPassCondToOmega1_2}, \eqref{eq:NonPassToSemiPassCondToOmega1_3}  gives a lower bound on $\mod{(\pi t_{se}+\phi,2\pi)}$ and \eqref{eq:NonPassToSemiPassCondToOmega1_1} gives an exact value.
\end{itemize}

\section{The Jacobian $DP_{\Gamma^{-}\Gamma^{+}}$}\label{sec:JacobianTtoB}

The entries of the Jacobian matrix $DP_{\Gamma^{-}\Gamma^{+}}$ for all the periodic motions considered in this work are as follows: 
\begin{align}\label{eq:dtkp1dtk}
    \frac{\partial t_{k+1}}{\partial t_{k}}
   & = \frac{-r\dot{Z}_{k}+(f(t_{k})-L^{+})\Delta t_{k}}{-r\dot{Z}_{k}+F_1(t_{k+1})-F_1(t_{k})-L^{+}\Delta t_{k}}\\\label{eq:dtkp1dZdotk}
    \frac{\partial t_{k+1}}{\partial \dot{Z}_{k}} &=  \frac{r\Delta t_{k}}{-r\dot{Z}_{k}+F_1(t_{k+1})-F_1(t_{k})-L^{+}\Delta t_{k}} \\\label{eq:dZdotkp1dtk}
    \frac{\partial \dot{Z}_{k+1}}{\partial t_{k}}  &= \frac{\partial t_{k+1}}{\partial t_{k}}(f(t_{k+1})-L^{+})-(f(t_{k})-L^{+})\\\label{eq:dZdotkp1dZdotk}
    \frac{\partial \dot{Z}_{k+1}}{\partial \dot{Z}_{k}} 
    & = -r +\frac{\partial t_{k+1}}{\partial \dot{Z}_{k}}(f(t_{k+1})-L^{+})
\end{align}
\section{Calculation of harvested energy}\label{apx:EnergyCalculations}
\begin{table}[H]
    \centering
    \begin{tabular}{|p{0.4\textwidth}|p{0.52\textwidth}|}
    \hline
    \textbf{Formulas} & \textbf{Parameter description}\\\hline\hline
 &    $U^{imp}_k$ - the voltage generated by the membrane deformation at the $k^{\textrm{th}}$ impact\\
$U^{imp}_k = \bigg[\dfrac{A_k}{\pi R^2_{c}}\bigg]^2\cdot U_{in}$ &$U_{\text{in}} = 2000 \, \text{mV}$ - constant input voltage applied to the membranes\\
&$R_c = 6.3 \, \text{mm}$ - the radius of the undeformed membrane\\
\hline
   $ A_k = 2\pi R^2_b(1 - \cos\alpha_k) + \pi R_c^2 - \pi(R_b\sin\alpha_k)^2\cos\alpha_k$
 & $A_k$ - the area of the membrane at the deformed state\\
 &$R_b = 5 \, \text{mm}$ - the radius of the bullet\\\hline
$\cos\alpha_k = \dfrac{-2R_b(\delta_k - R_b) + 2R_c}{\sqrt{R_c^2 + (\delta_k - R_b)^2} - 2\delta_kR_b}$ & $\alpha_k$ - angle at the $k^{\textrm{th}}$ impact defined by the largest deflection $\delta_k$ of the membrane\\
\hline
&$ K = 4.0847\cdot 10^5$ and $\nu = 2.6$ - parameters of the elastic force of the membrane,\\
$\delta_k = \dfrac{\nu + 1}{2K}mV_k^2$ & $V_k$ - the relative dimensional velocity at the $k^{\textrm{th}}$ impact, proportional to $\dot{Z}_k$.\\
\hline
    \end{tabular}
    \caption{The output voltage calculation of the VI-EH device through the membrane deformation upon impact
from \cite{YURCHENKO2017456}.}
    \label{tab:EnergyCalcDetails}
\end{table}

\section{Comparison of expressions in \cite{Serdukova2019StabilityAB} for the 1:1 periodic solution}\label{sec:ZeroFrictionComparison}

In this section, we follow the steps employed in \cite{Serdukova2019StabilityAB} to obtain a 1:1 periodic solution and compare the resulting expressions when dry friction is present in the system, i.e., $\mu_k>0$. 
\begin{enumerate}
    \item Adding \eqref{eq:FirstImpactZdot} and \eqref{eq:SecondImpactZdot} gives:
    \begin{align*}
  \dot{Z}_k + \dot{Z}_{k+1} &= -r(\dot{Z}_{k-1}+ \dot{Z}_{k})+  F_1(t_{k+1})-F_1(t_{k-1})-L^{-}\Delta t_{k-1}-L^{+}\Delta t_{k}  
    \end{align*}
Applying the periodicity conditions and solving for $\dot{Z}_{k}$ gives:
\begin{align}\label{eq:SecondImpactZdotInTemsOfFirst}
   \dot{Z}_k &= -\dot{Z}_{k-1}+\frac{-L^{-}\Delta t_{k-1}-L^{+}\Delta t_{k}}{1+r}
\end{align}
Note: If we write $\Delta t_{k-1}=T-\Delta t_{k}$, then \eqref{eq:SecondImpactZdotInTemsOfFirst} can be rewritten as:
\begin{align}\label{eq:SecondImpactZdotInTemsOfFirstB}
     \dot{Z}_k &= -\dot{Z}_{k-1}+\frac{-L^{-}T+(L^{-}-L^{+})\Delta t_{k}}{1+r}
\end{align}
\item Substitute \eqref{eq:SecondImpactZdotInTemsOfFirstB} in \eqref{eq:SecondImpactZdot}:
\begin{align}\label{eq:DifOfF1betweenFirstTwoImpacts}
    F_1(t_{k})-F_1(t_{k-1}) &= (r-1)\dot{Z}_{k-1}+\frac{-L^{+}\Delta t_{k}+rL^{-}\Delta t_{k-1}}{1+r}
\end{align}
\item Adding \eqref{eq:FirstImpactZ} and \eqref{eq:SecondImpactZ} and using the periodicity conditions \eqref{eq:PeriodicityCond_SimpleSol}, we get:
\begin{align}\label{eq:Zdotkm1InTermsOfZdotk}
     r\dot{Z}_{k-1}\Delta t_{k-1}& = -r\dot{Z}_{k}\Delta t_{k}-F_1(t_{k-1})\Delta t_{k-1} -F_1(t_{k})\Delta t_{k}-L^{-}\frac{(\Delta t_{k-1})^2}{2}-L^{+}\frac{(\Delta t_{k})^2}{2}
\end{align}
Substituting \eqref{eq:SecondImpactZdotInTemsOfFirstB} and \eqref{eq:DifOfF1betweenFirstTwoImpacts} into \cref{eq:Zdotkm1InTermsOfZdotk}, and solving for $\dot{Z}_{k-1}$ yields:
\begin{align}\label{eq:InitialZdotInTermsOft0}
    \dot{Z}_{k-1} & = \frac{1}{r\Delta t_{k-1}-\Delta t_{k}}\bigg(-TF_1(t_{k-1})-L^{-}\frac{(\Delta t_{k-1})^2}{2}+L^{+}\frac{(\Delta t_{k})^2}{2}\bigg)
\end{align}
\item We then eliminate $F_1(t_{k-1})$ in \eqref{eq:FirstImpactZ} using \eqref{eq:InitialZdotInTermsOft0} and obtain an expression for the difference $F_2(t_{k})- F_2(t_{k-1})$:
\begin{align}\label{eq:F2DifferenceBetweenImpacts}
   F_2(t_k)-F_2(t_{k-1}) & =  -d + \frac{\Delta t_{k-1}\Delta t_{k}}{2T}\bigg(L^{+}\Delta t_{k}+L^{-}\Delta t_{k-1}+2(1+r)\dot{Z}_{k-1} \bigg) 
\end{align}
\item Square and add \eqref{eq:DifOfF1betweenFirstTwoImpacts} and \eqref{eq:F2DifferenceBetweenImpacts}
\begin{align}\label{eq:squaredSumOfF1F2Difs}
\begin{split}
      [F_1(t_{k})-F_1(t_{k-1})]^2 &+  [F_2(t_{k})-F_2(t_{k-1})]^2 = \bigg[(r-1)\dot{Z}_{k-1}+\frac{-L^{+}\Delta t_{k}+rL^{-}\Delta t_{k-1}}{1+r}\bigg]^2\\
     &+ \bigg[-d + \frac{\Delta t_{k-1}\Delta t_{k}}{2T}\bigg(L^{+}\Delta t_{k}+L^{-}\Delta t_{k-1}+2(1+r)\dot{Z}_{k-1}\bigg)\bigg]^2
\end{split}
\end{align}
Now, for a 1:1 periodic motion, we define $0\leq q\leq 1$ to be the fraction of the period that corresponds to the motion from $\Gamma^{-}$ to $\Gamma^{+}$, such that $\Delta t_{k-1} = qT$, and $\Delta t_{k}=T(1-q)$, where $T=2$. Also, we let $t_{k-1}=0$, and therefore, $\varphi_{k-1}=\mod(\phi,2\pi)=\phi\in[0,2\pi)$. In the following we obtain a reduced system of equations as formulated in \cite{Serdukova2019StabilityAB} capturing the 1:1 periodic motion in terms of the variables $\dot{Z}_{k-1}$, $q$ and $\phi$. Note that in \cite{Serdukova2019StabilityAB} $\Delta t_{k-1} = T_u$ and $\Delta t_{k} = T_d$.

For a fixed set of parameters, using \eqref{eq:DifOfF1betweenFirstTwoImpacts} we get an equation of $\dot{Z}_{k-1}$ in terms of $q$ (or $\Delta t_{k-1}$) and $\phi$.
\begin{align}\label{eq:RSE1to1_Zdotkm1}
      \dot{Z}_{k-1}& =   \frac{1}{r-1} \bigg[\frac{1}{\pi}\sin(2\pi nq+\pi t_{k-1}+\phi)- \frac{1}{\pi}\sin(\pi t_{k-1}+\phi) - \frac{-L^{+}T+(rL^{-}+L^{+})\Delta t_{k-1}}{1+r}\bigg]
\end{align}
By \eqref{eq:InitialZdotInTermsOft0} we obtain an equation of $\varphi=\mod{(\phi,2\pi)}$ in terms of $\dot{Z}_{k-1}$ and $q$ appearing in the definition of $\Delta t_{k-1}$ and $\Delta t_{k}$:
\begin{align}
\begin{split}
   \phi &= \arcsin{\bigg(\frac{\pi}{T}\bigg[-(r\Delta t_{k-1}-\Delta t_{k})\dot{Z}_{k-1}-L^{-}\frac{(\Delta t_{k-1})^2}{2}+L^{+}\frac{(\Delta t_{k})^2}{2}\bigg]\bigg)}-\pi t_{k-1} + 2N\pi\\ &\textrm{(for some }N\in\mathbb{N})\\
   \end{split}
   \end{align}
\begin{align}
   \varphi  &= \arcsin{\bigg(\frac{\pi}{T}\bigg[-(r\Delta t_{k-1}-\Delta t_{k})\dot{Z}_{k-1}-L^{-}\frac{(\Delta t_{k-1})^2}{2}+L^{+}\frac{(\Delta t_{k})^2}{2}\bigg]\bigg)}-\pi t_{k-1}
\end{align}
\end{enumerate}
By \eqref{eq:squaredSumOfF1F2Difs}, we obtain a third equation in terms of $\dot{Z}_{k-1}, \phi, q$, which we can also view as a quadratic equation in terms of $\dot{Z}_{k-1}$:
\begin{align*}
      \frac{1}{\pi^2}[\sin(2nq\pi+\pi t_{k-1}+\phi)-\sin(\pi t_{k-1}+\phi)]^2 + \frac{1}{\pi^4} [\cos(2nq\pi+\pi t_{k-1}+\phi)-\cos(\pi t_{k-1}+\phi)]^2 =\\
       = \bigg[(r-1)\dot{Z}_{k-1}+\frac{-L^{+}\Delta t_{k}+rL^{-}\Delta t_{k-1}}{1+r}\bigg]^2+ \bigg[-d + \frac{\Delta t_{k-1}\Delta t_{k}}{2T}\bigg(L^{+}\Delta t_{k}+L^{-}\Delta t_{k-1}+2(1+r)\dot{Z}_{k-1}\bigg)\bigg]^2
\end{align*}
Finally, using \eqref{eq:DifOfF1betweenFirstTwoImpacts} we can obtain expressions for the following partial derivatives of the associated Jacobian matrices as defined in \eqref{eq:JacobiansGeneral}, where $P_{\Gamma^{+}\dots\Gamma^{-}}=P_{\Gamma^{+}\Gamma^{-}}$:
\begin{align}\label{eq:dtkdZdotkm1_SM}
    \frac{\partial t_{k}}{\partial \dot{Z}_{k-1}} &= \frac{-r\Delta t_{k-1}}{\dot{Z}_{k-1}+\frac{L^{+}\Delta t_{k}+L^{-}\Delta t_{k-1}}{1+r}}\\\label{eq:dtkdtkm1_SM}
    \frac{\partial t_{k}}{\partial t_{k-1}}
    &= \frac{r\dot{Z}_{k-1}-f(t_{k-1})\Delta t_{k-1}+L^{-}\Delta t_{k-1}}{\dot{Z}_{k-1}+\frac{L^{+}\Delta t_{k}+L^{-}\Delta t_{k-1}}{1+r}}\\\label{eq:dtkp1dtk_SM} 
     \frac{\partial t_{k+1}}{\partial t_{k}} &= 
    \frac{r\dot{Z}_{k}-f(t_{k})\Delta t_{k-1}+L^{+}\Delta t_{k}}{\dot{Z}_{k}+\frac{L^{-}\Delta t_{k-1}+L^{+}\Delta t_{k}}{1+r}}\\\label{eq:dtkp1dZdotk_SM} 
    \frac{\partial t_{k+1}}{\partial \dot{Z}_{k}}
    &= \frac{-r\Delta t_{k}}{\dot{Z}_k+\frac{L^{+}\Delta t_{k}+L^{-}\Delta t_{k-1}}{1+r}} 
\end{align}
The expressions obtained above are equivalent to \eqref{eq:dtkdZdotkm1}, \eqref{eq:dtkdtkm1} in \cref{sec:StabAnalysis_SimpleOneOne} and \eqref{eq:dtkp1dtk},\eqref{eq:dtkp1dZdotk} in \cref{sec:JacobianTtoB}. The expressions for the remaining partial derivatives, namely $\frac{\partial \dot{Z}_{k+j}}{\partial t_{k+j-1}}$ and $\frac{\partial \dot{Z}_{k+j}}{\partial \dot{Z}_{k+j-1}}$, for $j=0,1$, are the same as \eqref{eq:dZdotkdtkm1},\eqref{eq:dZdotkdZdotkm1},\eqref{eq:dZdotkp1dtk} and \eqref{eq:dZdotkp1dZdotk}.

\begin{remark}
    Note that if $\mu_{k}=0$, and thus, $L^{+}=L^{-}=-\bar{g}$, \eqref{eq:SecondImpactZdotInTemsOfFirst}-\eqref{eq:dtkp1dZdotk_SM} correspond to the expressions derived in \cite{Serdukova2019StabilityAB}. 
\end{remark}
\section{Analysis of the 1:1$^{\Sigma}_{s}$ periodic solution}\label{sec:SlidingSol_MS}

\subsection{Reduced system of equations for 1:1$^{\Sigma}_{s}$ periodic solution}\label{sec:RSESliding_Calculations}

Below we use steps (1)-(3) in the framework described in \cref{sec:Framework} to obtain the reduced system of equations \eqref{eq:RSESlidingSol} for the 1:1$^{\Sigma}_{s}$ periodic solution:
\begin{enumerate}
    \item  First impact to sliding onset, $(t_{k-1},d/2,\dot{Z}_{k-1})\to (t^{1}_{k-1},Z^{1}_{k-1},0)$, from \eqref{eq:SlidingOnsetEquationZdotFromOmega2} we have:
    \begin{align*}
           0 &= -r\dot{Z}_{k-1}+F_1(t^{1}_{k-1})-F_1(t^{0}_{k-1})-L^{-}\Delta t^{0}_{k-1}
    \end{align*}
    This equation corresponds to \eqref{eq:RSESlidingSol_1} and involves $\Delta t^{0}_{k-1}$, $\phi$, $\dot{Z}_{k-1}$.
  \item Sliding onset to sliding exit, $(t^{1}_{k-1},Z^{1}_{k-1},0)\to (t^{2}_{k-1},Z^{2}_{k-1},0)$, from \eqref{eq:SlidingConditionInMS_2} we have: 
  \begin{align*}
      \mod(\pi t^{2}_{k-1}+\phi,2\pi) = \arccos{(L^{-})}
  \end{align*}
  This equation corresponds to \eqref{eq:RSESlidingSol_2} and involves $\Delta t^{j-1}_{k-1}$, for $j=1,2$, and $\phi$.
     \item Sliding exit to second impact, $(t^{2}_{k-1},Z^{2}_{k-1},0)\to (t^{3}_{k-1},-d/2,\dot{Z}^{3}_{k-1})$, we substitute \eqref{eq:SlidingOnsetEquationZFromOmega2} and \eqref{eq:ZEqSlidingExit_SlidingSol} into \eqref{eq:SlidingtoSecondImpactZdotFromOmega2}, to obtain \eqref{eq:RSESlidingSol_3} (repeated below): 
    \begin{align*}
         -d &=-r\dot{Z}_{k-1}\Delta t^{0}_{k-1}+F_2(t^{1}_{k-1})-F_2(t^{0}_{k-1})-F_1(t^{0}_{k-1})\Delta t^{0}_{k-1}-L^{-}\frac{(\Delta t^{0}_{k-1})^2}{2}\\
       &+F_2(t^{3}_{k-1})-F_2(t^{2}_{k-1})-F_1(t^{2}_{k-1})\Delta t^{2}_{k-1}-L^{-}\frac{(\Delta t^{2}_{k-1})^2}{2}
       \end{align*}
    \eqref{eq:RSESlidingSol_3} involves $\phi$, $\dot{Z}_{k-1}$, and $\Delta t^{j-1}_{k-1}$, for $j=1,2,3$.
    \item  Second impact to third impact, $(t_{k},-d/2,\dot{Z}_{k})\to (t_{k+1},d/2,\dot{Z}_{k+1})$, we substitute \eqref{eq:SlidingtoSecondImpactZdotFromOmega2} and periodicity conditions into \eqref{eq:SecondImpactZinSlidingSolution} to obtain \eqref{eq:RSESlidingSol_4} (repeated below), where $t_k=t^{3}_{k-1}$:
    \begin{align*}
        d&=-(r+1)\Delta t_{k}F_1(t^{3}_{k-1})+r\Delta t_{k}F_1(t^{2}_{k-1})+r\Delta t_{k}L^{-}\Delta t^{2}_{k-1}+F_2(t^{0}_{k-1})-F_2(t^{3}_{k-1})-L^{+}\frac{(\Delta t_{k})^2}{2}
    \end{align*}
    \eqref{eq:RSESlidingSol_4} involves $\phi$, $\dot{Z}_{k-1}$ and $\Delta t^{j-1}_{k-1}$, for $j=1,2,3$ (explicitly and implicitly through $\Delta t_{k} = T - \sum_{j=1}^{3} \Delta t^{j-1}_{k-1}$).
    \item Second impact to third impact, $(t_{k},-d/2,\dot{Z}_{k})\to (t_{k+1},d/2,\dot{Z}_{k+1})$, we substitute \eqref{eq:SecondImpactZinSlidingSolution}, \eqref{eq:PeriodicityCond_SlidingSol} and other periodicity conditions into \eqref{eq:SecondImpactZdotinSlidingSolution} to obtain \eqref{eq:RSESlidingSol_5} (repeated below), where $t_{k}=t^{3}_{k-1}$:
    \begin{align*}
      \dot{Z}_{k-1}& = -(r+1)F_1(t^{3}_{k-1})+rF_1(t^{2}_{k-1})+rL^{-}\Delta t^{2}_{k-1}+F_1(t^{0}_{k-1})-L^{+}\Delta t_{k}    
    \end{align*}
      \eqref{eq:RSESlidingSol_5} involves $\phi$, $\dot{Z}_{k-1}$, $\Delta t^{j-1}_{k-1}$, for $j=1,2,3$.
\end{enumerate}
\subsection{Stability analysis for the 1:1$^{\Sigma}_{s}$ periodic solution}
Here, we show some details regarding the calculations of the the partial derivatives involved in the Jacobians of the 1:1$^{\Sigma}_{s}$ periodic solution. In particular, we show the details for finding the partial derivatives $\frac{\partial Z^{1}_{k-1}}{\partial \dot{Z}_{k-1}}$ and $\frac{\partial Z^{1}_{k-1}}{\partial t^{0}_{k-1}}$ associated with the motion from $\Gamma^{+}$ to $\Sigma$ from $(t^{0}_{k-1},d/2,\dot{Z}_{k-1})$ to $(t^{1}_{k-1},Z^{1}_{k-1},0)$: 
 \begin{align*}
    \frac{\partial Z^{1}_{k-1}}{\partial \dot{Z}_{k-1}} & = -r\dot{Z}_{k-1}\frac{\partial t^{1}_{k-1}}{\partial \dot{Z}_{k-1}}-r\Delta t^{0}_{k-1}+F_1(t^{1}_{k-1})\frac{\partial t^{1}_{k-1}}{\partial \dot{Z}_{k-1}}-F_1(t^{0}_{k-1})\frac{\partial t^{1}_{k-1}}{\partial \dot{Z}_{k-1}}-L^{-}\Delta t^{0}_{k-1}\frac{\partial t^{1}_{k-1}}{\partial \dot{Z}_{k-1}}\\
    & = \bigg[-r\dot{Z}_{k-1}+F_1(t^{1}_{k-1})-F_1(t^{0}_{k-1})-L^{-}\Delta t^{0}_{k-1}\bigg]\frac{\partial t^{1}_{k-1}}{\partial \dot{Z}_{k-1}} - r\Delta t^{0}_{k-1}= - r\Delta t^{0}_{k-1}\,,
\end{align*}
since $-r\dot{Z}_{k-1}+F_1(t^{1}_{k-1})-F_1(t^{0}_{k-1})-L^{-}\Delta t^{0}_{k-1}=0$, when evaluated at $(\dot{Z}^{0}_{k-1},\phi, \Delta t^{0}_{k-1},\Delta t^{1}_{k-1}, \Delta t^{2}_{k-1})$ capturing the 1:1$^{\Sigma}_{s}$ periodic solution.
 \begin{align*}
    \frac{\partial Z^{1}_{k-1}}{\partial t^{0}_{k-1}} & = -r\dot{Z}_{k-1}\bigg(\frac{\partial t^{1}_{k-1}}{\partial t^{0}_{k-1}}-1\bigg)+F_1(t^{1}_{k-1})\frac{\partial t^{1}_{k-1}}{\partial t^{0}_{k-1}}-F_1(t^{0}_{k-1})-F_1(t^{0}_{k-1})\bigg(\frac{\partial t^{1}_{k-1}}{\partial t^{0}_{k-1}}-1\bigg)\\
    & -f(t^{0}_{k-1})\Delta t^{0}_{k-1}-L^{-}\Delta t^{0}_{k-1}\bigg(\frac{\partial t^{1}_{k-1}}{\partial t^{0}_{k-1}}-1\bigg)\\
    &=  \bigg[-r\dot{Z}_{k-1}+F_1(t^{1}_{k-1})-F_1(t^{0}_{k-1})-L^{-}\Delta t^{0}_{k-1}\bigg]\bigg(\frac{\partial t^{1}_{k-1}}{\partial \dot{Z}_{k-1}}-1\bigg) + F_1(t^{1}_{k-1})-F_1(t^{0}_{k-1}) - f(t^{0}_{k-1})\Delta t^{0}_{k-1}\\
    &=  0 + F_1(t^{1}_{k-1})-F_1(t^{0}_{k-1}) -  f(t^{0}_{k-1})\Delta t^{0}_{k-1}
\end{align*}
\section{1:1$^{\Sigma}_{c}$ periodic solution}
\label{sec:LoopingSolAnalysis}
This solution is described by the following composition of maps: $P_{\Gamma^{-}\Gamma^{+}}\circ P_{\Sigma\Gamma^{-}}\circ P^{+}_{\Sigma\Sigma}\circ P_{\Gamma^{+}\Sigma}$. 
\subsection{Analytical expressions for the 1:1$^{\Sigma}_{c}$ solution}
We outline the key time points that describe this solution using the notation introduced in \eqref{eq:EqMotionVel}, \eqref{eq:EqMotionPosition}:
\begin{itemize}
    \item First impact on $\Gamma^{+}$ $\left(Z=\frac{d}{2}\right)$ occurs at $t=t^{0}_{k-1}$, where $t^{0}_{k-1}=t_{k-1}$.
    \item First crossing at $\Sigma$ $(\dot{Z}=0)$ at $t= t^{1}_{k-1}$.
    \item The solution evolves based on the vector field in $\Sigma^{+}$.
    \item Second crossing at $\Sigma$ $(\dot{Z}=0)$ at $t=t^{2}_{k-1}$.
     \item Second impact on $\Gamma^{-}$ $\left(Z=-\frac{d}{2}\right)$ occurs at $t=t^{3}_{k-1}=t_{k}$.
     \item Third impact on $\Gamma^{+}$ $\left(Z=\frac{d}{2}\right)$ occurs at $t=t^{1}_{k}=t_{k+1} = t_{k-1}+T$.
\end{itemize}
Since no sliding occurs while crossing the boundary $\Sigma$, the following conditions should hold at the two crossing times, $t^{1}_{k-1}$ and $t^{2}_{k-1}$.
\begin{align}\label{eq:LoopingSol_FirstCrossingCond}
    \cos(\pi t^{1}_{k-1} +\phi) > L^{+} \,, \cos(\pi t^{1}_{k-1} +\phi) > L^{-}\\\label{eq:LoopingSol_SecondCrossingCond}
     \cos(\pi t^{2}_{k-1} +\phi) < L^{+}  \,, \cos(\pi t^{2}_{k-1} +\phi) < L^{-}
\end{align}
\noindent based on conditions \eqref{eq:CrossingConditionsMinusToPlus_1}, \eqref{eq:CrossingConditionsMinusToPlus_2} and \eqref{eq:CrossingConditionsPlusToMinus_1}, \eqref{eq:CrossingConditionsPlusToMinus_2} found in  \cref{sec:SlidingStick}.

\subsubsection*{From first impact to crossing, $P_{\Gamma^{+}\Sigma}: (t_{k-1},d/2,\dot{Z}_{k-1})\to(t^{1}_{k-1},Z^{1}_{k-1},0)$}
\begin{align}
\label{eq:CrossingEquationZdotFromOmega2A}
     0&=\dot{Z}^{1}_{k-1} = -r\dot{Z}_{k-1}+F_1(t^{1}_{k-1})-F_1(t^{0}_{k-1})-L^{-}\Delta t^{0}_{k-1}\\\label{eq:CrossingEquationZFromOmega2A}
     \begin{split}
           Z^{1}_{k-1} &= \frac{d}{2}-r\dot{Z}_{k-1}\Delta t^{0}_{k-1}+F_2(t^{1}_{k-1})-F_2(t^{0}_{k-1})-F_1(t^{0}_{k-1})\Delta t^{0}_{k-1}-L^{-}\frac{(\Delta t^{0}_{k-1})^2}{2}
     \end{split}
\end{align}
\subsubsection*{From first crossing to second crossing, $P^{+}_{\Sigma\Sigma}: (t^{1}_{k-1},Z^{1}_{k-1},0)\to(t^{2}_{k-1},Z^{2}_{k-1},0)$}
\begin{align}
\label{eq:CrossingSolutionZdotFromOmega1}
    0=\dot{Z}^{2}_{k-1} & = F_1(t^{2}_{k-1})-F_1(t^{1}_{k-1})-L^{+}\Delta t^{1}_{k-1}\\\label{eq:CrossingSolutionZFromOmega1}
    Z^{2}_{k-1} & = Z^{1}_{k-1} + F_2(t^{2}_{k-1}) -F_2(t^{1}_{k-1}) - F_1(t^{1}_{k-1})\Delta t^{1}_{k-1}-\frac{L^{+}}{2}(\Delta t^{1}_{k-1})^2
\end{align}
\subsubsection*{From second crossing to second impact, $P_{\Sigma\Gamma^{-}}: (t^{2}_{k-1},Z^{2}_{k-1},0)\to(t_k,-d/2,\dot{Z}_k)$}
\begin{align}
\label{eq:SecondCrossingtoSecondImpactZdotFromOmega2}
      \dot{Z}^{3}_{k-1} &= F_1(t^{3}_{k-1})-F_1(t^{2}_{k-1})-L^{-}\Delta t^{2}_{k-1}\\\label{eq:SecondCrossingtoSecondImpactZFromOmega2}
    -\frac{d}{2} &= Z^{2}_{k-1}+F_2(t^{3}_{k-1})-F_2(t^{2}_{k-1})-F_1(t^{2}_{k-1})\Delta t^{2}_{k-1}-L^{-}\frac{(\Delta t^{2}_{k-1})^2}{2}
\end{align}
\subsubsection*{From second to third impact,  $P_{\Gamma^{+}\Gamma^{-}}: (t_k,-d/2,\dot{Z}_k)\to(t_{k+1},d/2,\dot{Z}_{k+1})$}

We use equations \eqref{eq:SecondImpactZdot} and \eqref{eq:SecondImpactZ} for $t_{k} = t^{3}_{k-1}$, $\dot{Z}_{k} = \dot{Z}^{3}_{k-1}$ and $t_{k+1} = t^{1}_{k}$.
\begin{align}
\label{eq:CrossingSolutionSecondImpactZdot}
     \dot{Z}_{k+1} &= -r\dot{Z}^{3}_{k-1}+F_1(t^{1}_{k})-F_1(t^{3}_{k-1})-L^{+}\Delta t_{k}\\\label{eq:CrossingSolutionSecondImpactZ}
    d &= -r\dot{Z}^{3}_{k-1}\Delta t_{k}+F_2(t^{1}_{k})-F_2(t^{3}_{k-1})-F_1(t^{3}_{k-1})\Delta t_{k}-L^{+}\frac{(\Delta t_{k})^2}{2}
\end{align}
\subsubsection*{Periodicity condition}
\begin{align}\label{eq:PeriodicityCond_CrossingSol}
    \dot{Z}_{k+1}&=\dot{Z}_{k-1}
\end{align}
\subsubsection{Reduced system of equations for the 1:1$^{\Sigma}_{c}$ periodic solution}
Using steps (1)-(3) in the framework described in \cref{sec:Framework}, we can reduce system \eqref{eq:CrossingEquationZdotFromOmega2A}-\eqref{eq:PeriodicityCond_CrossingSol} to the following subsystem \eqref{eq:RSELoopingSolution} and obtain a solution 
$(\dot{Z}_{k-1},\phi, \Delta t^{0}_{k-1},\Delta t^{1}_{k-1},\Delta t^{2}_{k-1})$, where $\Delta t^{j-1}_{k-1} = t^{j}_{k-1}-t^{j-1}_{k-1}$, $j=1,2,3$, that determines a 1:1$^{\Sigma}_{c}$ periodic solution as described in this section. We demonstrate these steps below:
\begin{subequations}\label{eq:RSELoopingSolution}
\begin{align}\label{eq:RSELoopingSolution_a}
      0 &= -r\dot{Z}_{k-1}+F_1(t^{1}_{k-1})-F_1(t^{0}_{k-1})-L^{-}\Delta t^{0}_{k-1}\\\label{eq:RSELoopingSolution_b}
    0 & = F_1(t^{2}_{k-1})-F_1(t^{1}_{k-1})-L^{+}\Delta t^{1}_{k-1}\\\label{eq:RSELoopingSolution_c}
    -d &=-r\dot{Z}_{k-1}\Delta t^{0}_{k-1}-F_2(t^{0}_{k-1})-F_1(t^{0}_{k-1})\Delta t^{0}_{k-1}-L^{-}\frac{(\Delta t^{0}_{k-1})^2}{2}\\ \nonumber
            & - F_1(t^{1}_{k-1})\Delta t^{1}_{k-1}-\frac{L^{+}}{2}(\Delta t^{1}_{k-1})^2
       +F_2(t^{3}_{k-1})-F_1(t^{2}_{k-1})\Delta t^{2}_{k-1}-L^{-}\frac{(\Delta t^{2}_{k-1})^2}{2}\\\label{eq:RSELoopingSolution_d}
       d&=-(r+1)\Delta t_{k}F_1(t^{3}_{k-1})+r\Delta t_{k}F_1(t^{2}_{k-1})+r\Delta t_{k}L^{-}\Delta t^{2}_{k-1}+F_2(t^{0}_{k-1})-F_2(t^{3}_{k-1})-L^{+}\frac{(\Delta t_{k})^2}{2}\\\label{eq:RSELoopingSolution_e}
        \dot{Z}_{k-1}& = -(r+1)F_1(t^{3}_{k-1})+rF_1(t^{2}_{k-1})+rL^{-}\Delta t^{2}_{k-1}+F_1(t^{0}_{k-1})-L^{+}\Delta t_{k}   
\end{align}
\end{subequations}
Note that once again the necessary variables are $(\dot{Z}_{k-1},\phi,\Delta t^{0}_{k-1})$, since $\Delta t^{1}_{k-1}$ and $\Delta t^{2}_{k-1}$ are implicit functions of $\dot{Z}_{k-1},\phi$ and $\Delta t^{0}_{k-1}$. Using the framework described in \cref{sec:Framework} we obtain the reduced system of equations \eqref{eq:RSELoopingSolution} for the 1:1$^{\Sigma}_{c}$ solution as shown below:
\begin{enumerate}
    \item From first impact to crossing, $(t_{k-1},d/2,\dot{Z}_{k-1})\to(t^{1}_{k-1},Z^{1}_{k-1},0)$, from \eqref{eq:CrossingEquationZdotFromOmega2A} we obtain \eqref{eq:RSELoopingSolution_a}, which involves $\Delta t^{0}_{k-1}$, $\phi$ (through $F_1$), $\dot{Z}_{k-1}$.
    \item From first crossing to second crossing, $(t^{1}_{k-1},Z^{1}_{k-1},0)\to(t^{2}_{k-1},Z^{2}_{k-1},0)$, from \eqref{eq:CrossingSolutionZdotFromOmega1} we obtain \eqref{eq:RSELoopingSolution_b}, which involves $\Delta t^{0}_{k-1}$, $\Delta t^{1}_{k-1}$, $\phi$.
     \item From second crossing to second impact, $(t^{2}_{k-1},Z^{2}_{k-1},0)\to(t_k,-d/2,\dot{Z}_k)$, we substitute  \eqref{eq:CrossingEquationZFromOmega2A} and \eqref{eq:CrossingSolutionZFromOmega1} into \eqref{eq:SecondCrossingtoSecondImpactZdotFromOmega2} to obtain \eqref{eq:RSELoopingSolution_c}.      \eqref{eq:RSELoopingSolution_c} involves $\Delta t^{j-1}_{k-1}$, for $j=1,\dots,3$, $\phi$, and $\dot{Z}^{0}_{k-1}$.
    \item From second to third impact, $(t_k,-d/2,\dot{Z}_k)\to(t_{k+1},d/2,\dot{Z}_{k+1})$, we substitute \eqref{eq:SecondCrossingtoSecondImpactZdotFromOmega2}, and periodicity conditions into \eqref{eq:CrossingSolutionSecondImpactZ} to obtain \eqref{eq:RSELoopingSolution_d}.     \eqref{eq:RSELoopingSolution_d} involves $\phi$, $\Delta t^{j-1}_{k-1}$, for $j=1,\dots,3$ (explicitly and implicitly through $\Delta t_{k}=T-\sum_{j=1}^{3}\Delta t^{j-1}_{k-1}$).
    \item From second to third impact, $(t_k,-d/2,\dot{Z}_k)\to(t_{k+1},d/2,\dot{Z}_{k+1})$, we substitute \eqref{eq:SecondCrossingtoSecondImpactZdotFromOmega2},  \eqref{eq:PeriodicityCond_CrossingSol} and other periodicity conditions into \eqref{eq:CrossingSolutionSecondImpactZdot} to obtain \eqref{eq:RSELoopingSolution_e}. \eqref{eq:RSELoopingSolution_e} involves $\phi$, $\dot{Z}_{k-1}$  and $\Delta t^{j-1}_{k-1}$, for $j=1,\dots,3$.
\end{enumerate}
\begin{remark}
 We use an adaptation for period-doubling of \eqref{eq:RSELoopingSolution} to obtain a 1:1$^{\Sigma}_{c}/2T$ solution in \Cref{fig:UnphysicalSolutions_Looping}. This system consists of nine equations.
\end{remark}
\subsection{Stability analysis of the 1:1$^{\Sigma}_{c}$ solution}
The 1:1$^{\Sigma}_{c}$ periodic solution is captured by the following composition of maps: $P_{\Gamma^{-}\Gamma^{+}}\circ P_{\Sigma \Gamma^{-}}\circ P^{+}_{\Sigma\Sigma}\circ P_{\Gamma^{+}\Sigma}$. The composition $P_{\Sigma \Gamma^{-}}\circ P^{+}_{\Sigma\Sigma}\circ P_{\Gamma^{+}\Sigma}$ maps the first impact at $\Gamma^{+}$ to the second impact at $\Gamma^{-}$.

We need the following quantities for the Jacobian of the map that takes in the $k-1^{\textrm{st}}$ impact and yields the $k^{\textrm{th}}$ impact: $\frac{\partial t_{k}}{\partial t_{k-1}}=\frac{\partial t_{k}}{\partial t^{0}_{k-1}}$, $\frac{\partial t_{k}}{\partial \dot{Z}_{k-1}}=\frac{\partial t_{k}}{\partial \dot{Z}_{k-1}}$, $\frac{\partial \dot{Z}_{k}}{\partial t_{k-1}}=\frac{\partial \dot{Z}_{k}}{\partial t^{0}_{k-1}}$, $\frac{\partial \dot{Z}_{k}}{\partial \dot{Z}_{k-1}}=\frac{\partial \dot{Z}_{k}}{\partial \dot{Z}_{k-1}}$.

The impact times and velocities are dependent not only on the previous impact times and velocities, but also on times and relative displacements between impacts. By the Chain Rule, we get:
\begin{align}\label{eq:dtkdtkm1_LoopingSol}
\begin{split}
    \frac{\partial t_{k}}{\partial t_{k-1}}= \frac{\partial t_{k}}{\partial t^{0}_{k-1}}&= \frac{\partial t_{k}}{\partial t^{2}_{k-1}} \frac{\partial t^{2}_{k-1}}{\partial t^{1}_{k-1}} \frac{\partial t^{1}_{k-1}}{\partial t^{0}_{k-1}}\\
    &+ \frac{\partial t_{k}}{\partial Z^{2}_{k-1}}\bigg(\frac{\partial Z^{2}_{k-1}}{\partial t^{2}_{k-1}}\frac{\partial t^{2}_{k-1}}{\partial t^{1}_{k-1}} \frac{\partial t^{1}_{k-1}}{\partial t^{0}_{k-1}}+  \frac{\partial Z^{2}_{k-1}}{\partial t^{1}_{k-1}} \frac{\partial t^{1}_{k-1}}{\partial t^{0}_{k-1}}+ \frac{\partial Z^{2}_{k-1}}{\partial Z^{1}_{k-1}} \frac{\partial Z^{1}_{k-1}}{\partial t^{0}_{k-1}} +\frac{\partial Z^{2}_{k-1}}{\partial Z^{1}_{k-1}} \frac{\partial Z^{1}_{k-1}}{\partial t^{1}_{k-1}}\frac{\partial t^{1}_{k-1}}{\partial t^{0}_{k-1}}\bigg)
    \end{split}
\end{align}
    \begin{align}\label{eq:dtkdZdotkm1_LoopingSol}
    \begin{split}
     \frac{\partial t_{k}}{\partial \dot{Z}_{k-1}}&=  \frac{\partial t_{k}}{\partial t^{2}_{k-1}} \frac{\partial t^{2}_{k-1}}{\partial t^{1}_{k-1}} \frac{\partial t^{1}_{k-1}}{\partial \dot{Z}_{k-1}} \\
     &+ \frac{\partial t_{k}}{\partial Z^{2}_{k-1}}\bigg(\frac{\partial Z^{2}_{k-1}}{\partial t^{2}_{k-1}}\frac{\partial t^{2}_{k-1}}{\partial t^{1}_{k-1}} \frac{\partial t^{1}_{k-1}}{\partial \dot{Z}_{k-1}}+  \frac{\partial Z^{2}_{k-1}}{\partial t^{1}_{k-1}} \frac{\partial t^{1}_{k-1}}{\partial \dot{Z}_{k-1}}+ \frac{\partial Z^{2}_{k-1}}{\partial Z^{1}_{k-1}} \frac{\partial Z^{1}_{k-1}}{\partial \dot{Z}_{k-1}} +\frac{\partial Z^{2}_{k-1}}{\partial Z^{1}_{k-1}} \frac{\partial Z^{1}_{k-1}}{\partial t^{1}_{k-1}}\frac{\partial t^{1}_{k-1}}{\partial \dot{Z}_{k-1}}\bigg)
     \end{split}
     \end{align}
\begin{equation}\label{eq:dZdotkdtkm1_LoopingSol}
 \frac{\partial \dot{Z}_{k}}{\partial t_{k-1}}= \frac{\partial \dot{Z}_{k}}{\partial t^{0}_{k-1}}= \frac{\partial\dot{Z}_{k}}{\partial t^{2}_{k-1}}\frac{\partial t^{2}_{k-1}}{\partial t^{1}_{k-1}} \frac{\partial t^{1}_{k-1}}{\partial t^{0}_{k-1}} +\dfrac{\partial  \dot{Z}_{k}}{\partial t_{k}}\dfrac{\partial  t_{k}}{\partial t^{0}_{k-1}}  
 \end{equation}
 \begin{equation}\label{eq:dZdotkdZdotkm1_LoopingSol}
          \frac{\partial \dot{Z}_{k}}{\partial \dot{Z}_{k-1}}= \frac{\partial \dot{Z}_{k}}{\partial t^{2}_{k-1}} \frac{\partial t^{2}_{k-1}}{\partial t^{1}_{k-1}} \frac{\partial t^{1}_{k-1}}{\partial \dot{Z}_{k-1}} + \dfrac{\partial  \dot{Z}_{k}}{\partial t_{k}}\dfrac{\partial  t_{k}}{\partial \dot{Z}_{k-1}}
\end{equation}

\subsubsection{Detailed expressions for stability analysis of the 1:1$^{\Sigma}_{c}$ periodic solution}\label{sec:SM_LoopingSol_LSA_Details}

In this section, we show detailed calculations of the partial derivatives in the Jacobian matrix $DP_{\Gamma^{+}\dots\Gamma^{-}}$ for the 1:1$^{\Sigma}_{c}$ periodic solution, where $P_{\Gamma^{+}\dots\Gamma^{-}}=P_{\Sigma\Gamma^{+}}\circ P^{+}_{\Sigma\Sigma}P_{\Gamma^{+}\Sigma}$. To obtain the necessary expressions \eqref{eq:dtkdtkm1_LoopingSol}-\eqref{eq:dZdotkdZdotkm1_LoopingSol} in \cref{sec:LoopingSolAnalysis} we implicitly or explicitly differentiate \eqref{eq:CrossingEquationZdotFromOmega2A}-\eqref{eq:CrossingSolutionSecondImpactZ}.
\begin{enumerate}
    \item From $(t^{0}_{k-1},d/2,\dot{Z}_{k-1})$ to $(t^{1}_{k-1},Z^{1}_{k-1},0)$:\\
         \begin{align}
    \frac{\partial t^{1}_{k-1}}{\partial t^{0}_{k-1}} & =\dfrac{f(t^{0}_{k-1})-L^{-}}{f(t^{1}_{k-1})-L^{-}}\\
    \frac{\partial t^{1}_{k-1}}{\partial \dot{Z}_{k-1}} & =\dfrac{r}{f(t^{1}_{k-1})-L^{-}}\\
    \begin{split}
   \frac{\partial Z^{1}_{k-1}}{\partial \dot{Z}_{k-1}} & = -r\dot{Z}_{k-1}\frac{\partial t^{1}_{k-1}}{\partial \dot{Z}_{k-1}}-r\Delta t^{0}_{k-1}+F_1(t^{1}_{k-1})\frac{\partial t^{1}_{k-1}}{\partial \dot{Z}_{k-1}}-F_1(t^{0}_{k-1})\frac{\partial t^{1}_{k-1}}{\partial \dot{Z}_{k-1}}\\
    &-L^{-}\Delta t^{0}_{k-1}\frac{\partial t^{1}_{k-1}}{\partial \dot{Z}_{k-1}}\\
    & = \bigg[-r\dot{Z}_{k-1}+F_1(t^{1}_{k-1})-F_1(t^{0}_{k-1})-L^{-}\Delta t^{0}_{k-1}\bigg]\frac{\partial t^{1}_{k-1}}{\partial \dot{Z}_{k-1}} - r\Delta t^{0}_{k-1}\\
    & =  - r\Delta t^{0}_{k-1} \quad \text{(using  \eqref{eq:CrossingEquationZdotFromOmega2A})}
    \end{split}
    \end{align}
     \begin{align}
    \frac{\partial Z^{1}_{k-1}}{\partial t^{0}_{k-1}} & = -r\dot{Z}_{k-1}\bigg(\frac{\partial t^{1}_{k-1}}{\partial t^{0}_{k-1}}-1\bigg)+F_1(t^{1}_{k-1})\frac{\partial t^{1}_{k-1}}{\partial t^{0}_{k-1}}-F_1(t^{0}_{k-1})-F_1(t^{0}_{k-1})\bigg(\frac{\partial t^{1}_{k-1}}{\partial t^{0}_{k-1}}-1\bigg)\\\nonumber
    & -f(t^{0}_{k-1})\Delta t^{0}_{k-1}-L^{-}\Delta t^{0}_{k-1}\bigg(\frac{\partial t^{1}_{k-1}}{\partial t^{0}_{k-1}}-1\bigg)\\\nonumber
    &=  \bigg[-r\dot{Z}_{k-1}+F_1(t^{1}_{k-1})-F_1(t^{0}_{k-1})-L^{-}\Delta t^{0}_{k-1}\bigg]\bigg(\frac{\partial t^{1}_{k-1}}{\partial \dot{Z}_{k-1}}-1\bigg) \\\nonumber
    &+ F_1(t^{1}_{k-1})-F_1(t^{0}_{k-1}) - f(t^{0}_{k-1})\Delta t^{0}_{k-1}\\\nonumber
    & = F_1(t^{1}_{k-1})-F_1(t^{0}_{k-1}) -  f(t^{0}_{k-1})\Delta t^{0}_{k-1} \quad \text{(using \eqref{eq:CrossingEquationZdotFromOmega2A})}
    \end{align}
\item From $(t^{1}_{k-1},Z^{1}_{k-1},0)$ to $(t^{2}_{k-1},Z^{2}_{k-1},0)$:
{\begin{align}
        \frac{\partial t^{2}_{k-1}}{\partial t^{1}_{k-1}} & = \frac{f(t^{1}_{k-1})-L^{+}}{f(t^{2}_{k-1})-L^{+}} \\
        \frac{\partial t^{2}_{k-1}}{\partial Z^{1}_{k-1}} & = 0 \\ \frac{\partial Z^{2}_{k-1}}{\partial Z^{1}_{k-1}} & = 1
\end{align}}  
\begin{align}
        \frac{\partial Z^{2}_{k-1}}{\partial t^{1}_{k-1}} & =  \frac{\partial Z^{1}_{k-1}}{\partial t^{1}_{k-1}}+F_1(t^{2}_{k-1})\frac{\partial t^{2}_{k-1}}{\partial t^{1}_{k-1}}-F_1(t^{1}_{k-1})-F_1(t^{1}_{k-1})\bigg(\frac{\partial t^{2}_{k-1}}{\partial t^{1}_{k-1}}-1\bigg)-f(t^{1}_{k-1})\Delta t^{1}_{k-1}\\\nonumber
        &-L^{+}\Delta t^{1}_{k-1}\bigg(\frac{\partial t^{2}_{k-1}}{\partial t^{1}_{k-1}}-1\bigg)\\\nonumber
        & =  \frac{\partial Z^{1}_{k-1}}{\partial t^{1}_{k-1}} +  \bigg[F_1(t^{2}_{k-1})-F_1(t^{1}_{k-1})-L^{+}\Delta t^{1}_{k-1}\bigg]\bigg(\frac{\partial t^{2}_{k-1}}{\partial t^{1}_{k-1}}-1\bigg) \\\nonumber
        &+F_1(t^{2}_{k-1})-F_1(t^{1}_{k-1})-f(t^{1}_{k-1})\Delta t^{1}_{k-1}\\ \nonumber
        &=F_1(t^{2}_{k-1})-F_1(t^{1}_{k-1})-f(t^{1}_{k-1})\Delta t^{1}_{k-1}\,,
\end{align}
using \eqref{eq:CrossingSolutionZdotFromOmega1} and $\dfrac{\partial Z^{1}_{k-1}}{\partial t^{1}_{k-1}} = -r\dot{Z}_{k-1} + F_1(t^{1}_{k-1})-F_1(t^{0}_{k-1})-L^{-}\Delta t^{0}_{k-1} = 0$.
\item From $(t^{2}_{k-1},Z^{2}_{k-1},0)$ to $(t_{k},-d/2,\dot{Z}_{k})$:
 \begin{align}
    \frac{\partial t_{k}}{\partial t^{2}_{k-1}} &  = \frac{-\frac{\partial Z^{2}_{k-1}}{\partial t^{2}_{k-1}}+f(t^{2}_{k-1})(t_{k}-t^{2}_{k-1})-L^{-}(t_{k}-t^{2}_{k-1})}{F_1(t_{k})-F_1(t^{2}_{k-1})-L^{-}(t_{k}-t^{2}_{k-1})}\,,
\end{align}
where $\dfrac{\partial Z^{2}_{k-1}}{\partial t^{2}_{k-1}} = F_1(t^{2}_{k-1})-F_1(t^{1}_{k-1})-L^{+}\Delta t^{1}_{k-1} = 0$.
 \begin{align}
    \frac{\partial t_{k}}{\partial Z^{2}_{k-1}} & = -\frac{1}{F_1(t_{k})-F_1(t^{2}_{k-1})-L^{-}(t_{k}-t^{2}_{k-1})} 
\end{align}
 \begin{align}
    \frac{\partial \dot{Z}_{k}}{\partial t^{2}_{k-1}} & = f(t_{k})\frac{\partial t_{k}}{\partial t^{2}_{k-1}}-f(t^{2}_{k-1})-L^{-}\bigg(\frac{\partial t_{k}}{\partial t^{2}_{k-1}}-1\bigg) \\ \frac{\partial \dot{Z}_{k}}{\partial Z^{2}_{k-1}} & = 0
\end{align}  
\end{enumerate}
\section{Analysis of the 1:1$^{\Sigma}_{cs}$ periodic solution}\label{sec:LoopingSlidingSolAnalysis}

Below we describe the 1:1$^{\Sigma}_{cs}$ periodic motion represented by
the composition $P_{\Gamma^{-}\Gamma^{+}}\circ P_{\Sigma\Gamma^{-}} \circ P^{s}_{\Sigma\Sigma}\circ P^{+}_{\Sigma\Sigma}\circ P_{\Gamma^{+}\Sigma}$. 

\subsection{Analytical expressions for the 1:1$^{\Sigma}_{cs}$ periodic solution}\label{sec:FSE_LoopingSlidingSol}
First, we outline the key time points for a 1:1$^{\Sigma}_{cs}$ periodic motion using the notation in \eqref{eq:EqMotionVel},\eqref{eq:EqMotionPosition}:
\begin{itemize}
    \item First impact on $\Gamma^{+}$ $\left(Z=\frac{d}{2}\right)$ occurs at $t=t^{0}_{k-1}$, where $t^{0}_{k-1}=t_{k-1}$.
    \item Intersecting $\Sigma$ $(\dot{Z}=0)$ at $t= t^{1}_{k-1}$ (crossing).
    \item The solution evolves based on the vector field in $\Sigma^{+}$.
    \item Intersecting $\Sigma$ $(\dot{Z}=0)$ at $t=t^{2}_{k-1}$ (while $\dot{Z}$ decreases in $\Sigma^{+}$), this is also the sliding onset.
    \item Exiting $\Sigma$ $(\dot{Z}=0)$ at $t=t^{3}_{k-1}$ (while $\dot{Z}$ decreases in $\Sigma^{-}$), this is also the sliding exit time.
     \item Second impact on $\Gamma^{-}$ $\left(Z=-\frac{d}{2}\right)$ occurs at $t=t^{4}_{k-1}=t_{k}$.
     \item Third impact on $\Gamma^{+}$ $\left(Z=\frac{d}{2}\right)$ occurs at $t=t^{1}_{k}=t_{k+1}=t_{k}+\Delta t_{k}$.
 \end{itemize}

Note that the crossing time, $t^{1}_{k-}1$, the onset time, $t^{2}_{k-1}$, and exit time, $t^{3}_{k-1}$, of sliding satisfy the following:
\begin{align}\label{eq:CrossingConditionInSM_1}
    \mod(\pi t^{1}_{k-1}+\phi, 2\pi)&\in[0,\arccos(L^{+}))\cup (2\pi-\arccos(L^{-}),2\pi]\\
\label{eq:SlidingConditionInSM_1}
    \mod(\pi t^{2}_{k-1}+\phi, 2\pi)&\in[\arccos(L^{+}),\arccos(L^{-}))\\\label{eq:SlidingConditionInSM_2}
     \mod(\pi t^{3}_{k-1}+\phi, 2\pi)&=\arccos(L^{-})\,,
\end{align}
as determined by conditions \eqref{eq:CrossingConditionsMinusToPlus_1}-\eqref{eq:CrossingConditionsMinusToPlus_2}, \eqref{eq:SlidingOnset2a}-\eqref{eq:SlidingOnset2c} and \eqref{eq:NonPassToSemiPassCondToOmega2_1}-\eqref{eq:NonPassToSemiPassCondToOmega2_3}.
These key time points then appear in the full system of equations for the 1:1$^{\Sigma}_{cs}$ periodic orbit, following \eqref{eq:EqMotionVel}, \eqref{eq:EqMotionPosition}.

\subsubsection*{From first impact to crossing,  $P_{\Gamma^{+}\Sigma}: (t_{k-1},d/2,\dot{Z}_{k-1})\to(t^{1}_{k-1},Z^{1}_{k-1},0)$}

\begin{align}
\label{eq:CrossingEquationZdotFromOmega2CSS}
     0&=\dot{Z}^{1}_{k-1}= -r\dot{Z}^{0}_{k-1}+F_1(t^{1}_{k-1})-F_1(t^{0}_{k-1})-L^{-}\Delta t^{0}_{k-1}\\\label{eq:CrossingEquationZFromOmega2CSS}
    Z^{1}_{k-1} &= \frac{d}{2}-r\dot{Z}^{0}_{k-1}\Delta t^{0}_{k-1}+F_2(t^{1}_{k-1})-F_2(t^{0}_{k-1})-F_1(t^{0}_{k-1})\Delta t^{0}_{k-1}-L^{-}\frac{(\Delta t^{0}_{k-1})^2}{2}
\end{align}

\subsubsection*{From first crossing to second crossing,  $P^{+}_{\Sigma\Sigma}: (t^{1}_{k-1},Z^{1}_{k-1},0)\to(t^{2}_{k-1},Z^{2}_{k-1},0)$}

\begin{align}
\label{eq:CrossingSolutionZdotFromOmega1CSS}
    0=\dot{Z}^{2}_{k-1} & = F_1(t^{2}_{k-1})-F_1(t^{1}_{k-1})-L^{+}\Delta t^{1}_{k-1}\\\label{eq:CrossingSolutionZFromOmega1CSS}
    Z^{2}_{k-1} & = Z^{1}_{k-1} + F_2(t^{2}_{k-1}) -F_2(t^{1}_{k-1}) - F_1(t^{1}_{k-1})\Delta t^{1}_{k-1}-\frac{L^{+}}{2}(\Delta t^{1}_{k-1})^2
\end{align}
\subsubsection*{From second crossing/sliding onset to end of sliding motion, $P^{s}_{\Sigma\Sigma}: (t^{2}_{k-1},Z^{2}_{k-1},0)\to(t^{3}_{k-1},Z^{3}_{k-1},0)$}

\begin{align}\label{eq:CSS_ZdotSlidingExit}
    \dot{Z}^{3}_{k-1} &=  \dot{Z}^{2}_{k-1}\\\label{eq:CSS_ZSlidingExit}
   Z^{3}_{k-1} &= Z^{2}_{k-1}\\\label{eq:CSS_slidingexiteq}
   \cos(\pi t^{3}_{k-1}+\phi)-L^{-}&=0
\end{align}

\subsubsection*{From sliding exit to second impact, $P_{\Sigma\Gamma^{-}}: (t^{3}_{k-1},Z^{3}_{k-1},0)\to(t^{4}_{k-1},-d/2,\dot{Z}^{4}_{k-1})$}

\begin{align}
\label{eq:SecondCrossingtoSecondImpactZdotFromOmega2CSS}
     \dot{Z}^{4}_{k-1} &= F_1(t^{4}_{k-1})-F_1(t^{3}_{k-1})-L^{-}\Delta t^{3}_{k-1}\\\label{eq:SecondCrossingtoSecondImpactZdotFromOmega2CSS}
    -\frac{d}{2} &= Z^{3}_{k-1}+F_2(t^{4}_{k-1})-F_2(t^{3}_{k-1})-F_1(t^{3}_{k-1})\Delta t^{3}_{k-1}-L^{-}\frac{(\Delta t^{3}_{k-1})^2}{2}
\end{align}

\subsubsection*{From second to third impact, $P_{\Gamma^{-}\Gamma^{+}}: (t_{k},-d/2,\dot{Z}_{k})\to(t_{k+1},d/2,\dot{Z}_{k+1})$}

We use \eqref{eq:SecondImpactZ} and \eqref{eq:SecondImpactZdot} for $t_{k} = t^{4}_{k-1}$, $\dot{Z}_{k} = \dot{Z}^{4}_{k-1}$ and $t_{k+1} = t^{1}_{k}$.
\begin{align}
\label{eq:CrossingSolutionSecondImpactZdotCSS}
     \dot{Z}_{k+1} &= -r\dot{Z}_k+F_1(t_{k+1})-F_1(t_{k})-L^{+}\Delta t_{k}\\\label{eq:CrossingSolutionSecondImpactZCSS}
    d &= -r\dot{Z}_{k}\Delta t_{k}+F_2(t_{k+1})-F_2(t_{k})-F_1(t_{k})\Delta t_{k}-L^{+}\frac{(\Delta t_{k})^2}{2}
\end{align}

\subsubsection*{Periodicity conditions}

\begin{align}\label{eq:PeriodicityCond_LoopingSlidingSol}
    \dot{Z}_{k+1}&=\dot{Z}^{0}_{k-1}
\end{align}
\subsection{Reduced system of equations for the 1:1$^{\Sigma}_{cs}$ periodic solution}\label{sec:RSE_LoopingSlidingSol}
 Following steps (1)-(3) described in \cref{sec:Framework}, we can reduce system \eqref{eq:CrossingEquationZdotFromOmega2CSS}-\eqref{eq:PeriodicityCond_LoopingSlidingSol} to the following subsystem \eqref{eq:RSE_OneOneLoopSlidingSol} and obtain a solution 
$(\dot{Z}_{k-1},\phi, \Delta t^{0}_{k-1},\Delta t^{1}_{k-1},\Delta t^{2}_{k-1},\Delta t^{3}_{k-1})$, corresponding to a 
1:1$^{\Sigma}_{cs}$ periodic solution. Recall that $\Delta t^{j-1}_{k-1} = t^{j}_{k-1}-t^{j-1}_{k-1}$, $j=1,2,3,4$.

\begin{subequations}\label{eq:RSE_OneOneLoopSlidingSol}
    \begin{align}\label{eq:RSE_OneOneLoopSlidingSol_A}
        0 &= -r\dot{Z}_{k-1}+F_1(t^{1}_{k-1})-F_1(t^{0}_{k-1})-L^{-}\Delta t^{0}_{k-1}\\\label{eq:RSE_OneOneLoopSlidingSol_B}
        0 & = F_1(t^{2}_{k-1})-F_1(t^{1}_{k-1})-L^{+}\Delta t^{1}_{k-1}\\\label{eq:RSE_OneOneLoopSlidingSol_C}
        0& =\cos(\pi t^{3}_{k-1}+\phi)-L^{-}\\ \label{eq:RSE_OneOneLoopSlidingSol_D}
        -d &= -r\dot{Z}_{k-1}\Delta t^{0}_{k-1}-F_2(t^{0}_{k-1})-F_1(t^{0}_{k-1})\Delta t^{0}_{k-1}-L^{-}\frac{(\Delta t^{0}_{k-1})^2}{2} \\\nonumber
 &+ F_2(t^{2}_{k-1}) - F_1(t^{1}_{k-1})\Delta t^{1}_{k-1}-\frac{L^{+}}{2}(\Delta t^{1}_{k-1})^2\\\nonumber
 &+F_2(t^{4}_{k-1})-F_2(t^{3}_{k-1})-F_1(t^{3}_{k-1})\Delta t^{3}_{k-1}-L^{-}\frac{(\Delta t^{3}_{k-1})^2}{2}\\\label{eq:RSE_OneOneLoopSlidingSol_E}
  d &= -(r+1)\Delta t_{k}F_1(t^{4}_{k-1})+rF_1(t^{3}_{k-1})\Delta t_{k}+F_2(t^{0}_{k-1})-F_2(t^{4}_{k-1})\\\nonumber
 & +r\Delta t_{k}L^{-}\Delta t^{3}_{k-1}-L^{+}\frac{(\Delta t_{k})^2}{2}\\\label{eq:RSE_OneOneLoopSlidingSol_F}
  \dot{Z}_{k-1}&= - (r+1)F_1(t^{4}_{k-1})+rF_1(t^{3}_{k-1})+rL^{-}\Delta t^{3}_{k-1}+F_1(t^{0}_{k-1})-L^{+}\Delta t_{k}
    \end{align}
\end{subequations}
\begin{enumerate}
    \item From first impact to first crossing, $(t_{k-1},d/2,\dot{Z}_{k-1})\to(t^{1}_{k-1},Z^{1}_{k-1},0)$, from \eqref{eq:CrossingEquationZFromOmega2CSS} we obtain \eqref{eq:RSE_OneOneLoopSlidingSol_A}, which involves $\Delta t^{0}_{k-1}$, $\phi$, and $\dot{Z}_{k-1}$.
    \item  From first crossing to sliding onset, $(t^{1}_{k-1},Z^{1}_{k-1},0)\to (t^{2}_{k-1},Z^{2}_{k-1},0)$, from \eqref{eq:CrossingSolutionZdotFromOmega1CSS} we obtain \eqref{eq:RSE_OneOneLoopSlidingSol_B} which involves $\phi$, and $\Delta t^{j-1}_{k-1}$, for $j=1,2$.
\item From sliding onset to sliding exit, $(t^{2}_{k-1},Z^{2}_{k-1},0)\to (t^{3}_{k-1},Z^{3}_{k-1},0)$, from \eqref{eq:CSS_slidingexiteq} we obtain \eqref{eq:RSE_OneOneLoopSlidingSol_C} which involves $\phi$, and $\Delta t^{j-1}_{k-1}$, for $j=1,2,3$.
\item From sliding exit to second impact, $(t^{3}_{k-1},Z^{3}_{k-1},0)\to (t_{k},-d/2,\dot{Z}_{k})$, we substitute \eqref{eq:CrossingEquationZFromOmega2CSS}, \eqref{eq:CrossingSolutionZFromOmega1CSS}, and \eqref{eq:CSS_ZSlidingExit} into \eqref{eq:SecondCrossingtoSecondImpactZdotFromOmega2CSS} to obtain \eqref{eq:RSE_OneOneLoopSlidingSol_D}. \eqref{eq:RSE_OneOneLoopSlidingSol_D}
 involves $\phi$, $\dot{Z}_{k-1}$, and $\Delta t^{j-1}_{k-1}$, for $j=1,2,3,4$.
\item From second impact to third impact, $(t_{k},-d/2,\dot{Z}_{k})\to(t_{k+1},d/2,\dot{Z}_{k+1})$, we substitute \eqref{eq:SecondCrossingtoSecondImpactZdotFromOmega2CSS} into \eqref{eq:CrossingSolutionSecondImpactZCSS} to obtain \eqref{eq:RSE_OneOneLoopSlidingSol_E}. \eqref{eq:RSE_OneOneLoopSlidingSol_E}
involves $\phi$, $\Delta t^{j-1}_{k-1}$, for $j=1,\dots,4$.
     \item From second impact to third impact, $(t_{k},-d/2,\dot{Z}_{k})\to(t_{k+1},d/2,\dot{Z}_{k+1})$, we substitute \eqref{eq:SecondCrossingtoSecondImpactZdotFromOmega2CSS}, \eqref{eq:PeriodicityCond_LoopingSlidingSol}, and other periodicity conditions into \eqref{eq:CrossingSolutionSecondImpactZdotCSS} to obtain \eqref{eq:RSE_OneOneLoopSlidingSol_F}. \eqref{eq:RSE_OneOneLoopSlidingSol_F} involves $\phi$, $\dot{Z}_{k-1}$ and $\Delta t^{j-1}_{k-1}$, for $j=1,\dots,4$.
     \end{enumerate}
\subsection{Stability analysis for the 1:1$^{\Sigma}_{cs}$ periodic solution}\label{sec:LoopingSlidingSol_StabilityAnalysis}

We calculate the partial derivatives in the Jacobian of the matrix  $P_{\Gamma^{+}\Gamma^{-}}=P_{\Gamma^{-}\Gamma^{+}}\circ P_{\Sigma \Gamma^{-}}\circ P^{s}_{\Sigma\Sigma}\circ P^{+}_{\Sigma\Sigma}\circ P_{\Gamma^{+}\Sigma}$, which maps the first impact at $\Gamma^{+}$ to the second impact at $\Gamma^{-}$. The expressions for the Jacobian of $P_{\Gamma^{-}\Gamma^{+}}$ are given in \cref{sec:JacobianTtoB}.
     \begin{align}
     \begin{split}
     \frac{\partial t_{k}}{\partial t_{k-1}}= \frac{\partial t^{4}_{k-1}}{\partial t^{0}_{k-1}}&=  \frac{\partial t^{4}_{k-1}}{\partial Z^{3}_{k-1}} \frac{\partial Z^{3}_{k-1}}{\partial Z^{2}_{k-1}} \bigg(\frac{\partial Z^{2}_{k-1}}{\partial t^{2}_{k-1}}\dfrac{\partial t^{2}_{k-1}}{\partial t^{1}_{k-1}}\dfrac{\partial t^{1}_{k-1}}{\partial t^{0}_{k-1}}+\frac{\partial Z^{2}_{k-1}}{\partial t^{1}_{k-1}}\frac{\partial t^{1}_{k-1}}{\partial t^{0}_{k-1}} \\
     &+ \frac{\partial Z^{2}_{k-1}}{\partial Z^{1}_{k-1}}\frac{\partial Z^{1}_{k-1}}{\partial t^{0}_{k-1}} + \frac{\partial Z^{2}_{k-1}}{\partial Z^{1}_{k-1}}\frac{\partial Z^{1}_{k-1}}{\partial t^{1}_{k-1}}\frac{\partial t^{1}_{k-1}}{\partial t^{0}_{k-1}}\bigg)
     \end{split}\\
             \begin{split}\frac{\partial t^{4}_{k-1}}{\partial \dot{Z}_{k-1}}&=  \frac{\partial t^{4}_{k-1}}{\partial Z^{3}_{k-1}}\frac{\partial Z^{3}_{k-1}}{\partial Z^{2}_{k-1}}\bigg(\frac{\partial Z^{2}_{k-1}}{\partial t^{2}_{k-1}}\dfrac{\partial t^{2}_{k-1}}{\partial t^{1}_{k-1}}\dfrac{\partial t^{1}_{k-1}}{\partial \dot{Z}_{k-1}}\\
     &+\frac{\partial Z^{2}_{k-1}}{\partial t^{1}_{k-1}}\frac{\partial t^{1}_{k-1}}{\partial \dot{Z}_{k-1}} + \frac{\partial Z^{2}_{k-1}}{\partial Z^{1}_{k-1}}\frac{\partial Z^{1}_{k-1}}{\partial \dot{Z}_{k-1}} + \frac{\partial Z^{2}_{k-1}}{\partial Z^{1}_{k-1}}\frac{\partial Z^{1}_{k-1}}{\partial t^{1}_{k-1}}\frac{\partial t^{1}_{k-1}}{\partial \dot{Z}_{k-1}}\bigg)
     \end{split}\\
       \frac{\partial \dot{Z}^{4}_{k-1}}{\partial t^{0}_{k-1}}&= (f(t^{4}_{k-1})-L^{-})\frac{\partial t^{4}_{k-1}}{\partial t^{0}_{k-1}}\\
        \frac{\partial \dot{Z}^{4}_{k-1}}{\partial \dot{Z}_{k-1}}&= (f(t^{4}_{k-1})-L^{-})\frac{\partial t^{4}_{k-1}}{\partial \dot{Z}_{k-1}}
\end{align}
\begin{enumerate}
    \item From $(t^{0}_{k-1},d/2,\dot{Z}_{k-1})$ to $(t^{1}_{k-1},Z^{1}_{k-1},0)$:
    \begin{align}
    \frac{\partial t^{1}_{k-1}}{\partial t^{0}_{k-1}} &  = \dfrac{f(t^{0}_{k-1})-L^{-}}{f(t^{1}_{k-1})-L^{-}}\\
    \frac{\partial t^{1}_{k-1}}{\partial \dot{Z}_{k-1}} &  = \dfrac{r}{f(t^{1}_{k-1})-L^{-}}
\end{align}
 \begin{align}
    \frac{\partial Z^{1}_{k-1}}{\partial \dot{Z}_{k-1}} & = -r\dot{Z}_{k-1}\frac{\partial t^{1}_{k-1}}{\partial \dot{Z}_{k-1}}-r\Delta t^{0}_{k-1}+F_1(t^{1}_{k-1})\frac{\partial t^{1}_{k-1}}{\partial \dot{Z}_{k-1}}-F_1(t^{0}_{k-1})\frac{\partial t^{1}_{k-1}}{\partial \dot{Z}_{k-1}}-L^{-}\Delta t^{0}_{k-1}\frac{\partial t^{1}_{k-1}}{\partial \dot{Z}_{k-1}}\\\nonumber
    & = \bigg[-r\dot{Z}_{k-1}+F_1(t^{1}_{k-1})-F_1(t^{0}_{k-1})-L^{-}\Delta t^{0}_{k-1}\bigg]\frac{\partial t^{1}_{k-1}}{\partial \dot{Z}_{k-1}} - r\Delta t^{0}_{k-1}\\\nonumber
    & = - r\Delta t^{0}_{k-1} \quad  \text{(from \eqref{eq:CrossingEquationZdotFromOmega2CSS})}
\end{align}
\begin{align}
    \frac{\partial Z^{1}_{k-1}}{\partial t^{0}_{k-1}} & = -r\dot{Z}_{k-1}\bigg(\frac{\partial t^{1}_{k-1}}{\partial t^{0}_{k-1}}-1\bigg)+F_1(t^{1}_{k-1})\frac{\partial t^{1}_{k-1}}{\partial t^{0}_{k-1}}-F_1(t^{0}_{k-1})-F_1(t^{0}_{k-1})\bigg(\frac{\partial t^{1}_{k-1}}{\partial t^{0}_{k-1}}-1\bigg)\\\nonumber
    & -f(t^{0}_{k-1})\Delta t^{0}_{k-1}-L^{-}\Delta t^{0}_{k-1}\bigg(\frac{\partial t^{1}_{k-1}}{\partial t^{0}_{k-1}}-1\bigg)\\\nonumber
    &=  \bigg[-r\dot{Z}_{k-1}+F_1(t^{1}_{k-1})-F_1(t^{0}_{k-1})-L^{-}\Delta t^{0}_{k-1}\bigg]\bigg(\frac{\partial t^{1}_{k-1}}{\partial \dot{Z}_{k-1}}-1\bigg)  \\\nonumber
    &+ F_1(t^{1}_{k-1})-F_1(t^{0}_{k-1}) - f(t^{0}_{k-1})\Delta t^{0}_{k-1}\\\nonumber
    &=  F_1(t^{1}_{k-1})-F_1(t^{0}_{k-1}) -  f(t^{0}_{k-1})\Delta t^{0}_{k-1} \quad  \text{(from \eqref{eq:CrossingEquationZdotFromOmega2CSS})}
\end{align}
\item From $(t^{1}_{k-1},Z^{1}_{k-1},0)$ to $(t^{2}_{k-1},Z^{2}_{k-1},0)$:
    \begin{align}
        \frac{\partial t^{2}_{k-1}}{\partial t^{1}_{k-1}} & = \frac{f(t^{1}_{k-1})-L^{+}}{f(t^{2}_{k-1})-L^{+}}\\
        \frac{\partial t^{2}_{k-1}}{\partial Z^{1}_{k-1}} & = 0\\
         \frac{\partial Z^{2}_{k-1}}{\partial Z^{1}_{k-1}} & = 1
\end{align}
 \begin{align}
        \frac{\partial Z^{2}_{k-1}}{\partial t^{1}_{k-1}} & =  \frac{\partial Z^{1}_{k-1}}{\partial t^{1}_{k-1}}+F_1(t^{2}_{k-1})\frac{\partial t^{2}_{k-1}}{\partial t^{1}_{k-1}}-F_1(t^{1}_{k-1})-F_1(t^{1}_{k-1})\bigg(\frac{\partial t^{2}_{k-1}}{\partial t^{1}_{k-1}}-1\bigg)-f(t^{1}_{k-1})\Delta t^{1}_{k-1}\\\nonumber
        &-L^{+}\Delta t^{1}_{k-1}\bigg(\frac{\partial t^{2}_{k-1}}{\partial t^{1}_{k-1}}-1\bigg)\\\nonumber
        & =  \frac{\partial Z^{1}_{k-1}}{\partial t^{1}_{k-1}} + \bigg[F_1(t^{2}_{k-1})-F_1(t^{1}_{k-1})-L^{+}\Delta t^{1}_{k-1}\bigg]\bigg(\frac{\partial t^{2}_{k-1}}{\partial t^{1}_{k-1}}-1\bigg) \\\nonumber
        &+F_1(t^{2}_{k-1})-F_1(t^{1}_{k-1})-f(t^{1}_{k-1})\Delta t^{1}_{k-1}\\\nonumber
        & =  \frac{\partial Z^{1}_{k-1}}{\partial t^{1}_{k-1}} +   F_1(t^{2}_{k-1})-F_1(t^{1}_{k-1})-f(t^{1}_{k-1})\Delta t^{1}_{k-1}
\end{align}
where $\frac{\partial Z^{1}_{k-1}}{\partial t^{1}_{k-1}} = -r\dot{Z}_{k-1} + F_1(t^{1}_{k-1})-F_1(t^{0}_{k-1})-L^{-}\Delta t^{0}_{k-1} = 0$ and using \eqref{eq:CrossingSolutionZdotFromOmega1CSS}.
\item From $(t^{2}_{k-1},Z^{2}_{k-1},0)$ to $(t^{3}_{k-1},Z^{3}_{k-1},0)$:
 \begin{align}
        \frac{\partial t^{3}_{k-1}}{\partial t^{2}_{k-1}} = \frac{\partial t^{3}_{k-1}}{\partial Z^{2}_{k-1}}=  \frac{\partial Z^{3}_{k-1}}{\partial t^{2}_{k-1}}  & =0\\
         \frac{\partial Z^{3}_{k-1}}{\partial Z^{2}_{k-1}} & = 1
\end{align}
\item From $(t^{3}_{k-1},Z^{3}_{k-1},0)$ to $(t^{4}_{k-1},-d/2,\dot{Z}^{4}_{k-1})$:
     \begin{align}
    \frac{\partial t^{4}_{k-1}}{\partial t^{3}_{k-1}} =\frac{\partial \dot{Z}^{4}_{k-1}}{\partial t^{3}_{k-1}} = \frac{\partial \dot{Z}^{4}_{k-1}}{\partial Z^{3}_{k-1}} & = 0\\
    \frac{\partial t^{4}_{k-1}}{\partial Z^{3}_{k-1}} & = -\frac{1}{F_1(t^{4}_{k-1})-F_1(t^{3}_{k-1})-L^{-}\Delta t^{3}_{k-1}}
\end{align}
\end{enumerate}
\section{Eigenvalues corresponding to periodic solutions}
\begin{figure}[hbtp!]
    \centering
    \begin{subfigure}[b]{0.45\textwidth}
    \centering
         \includegraphics[scale=0.4]{Figures/GrazingPaper_Figure3b_ParSet_muk=0p5_SimpleBranchFullDescription2.pdf}
    \caption{}
    \label{fig:SimpleSolutionBranch_MS}
    \end{subfigure}
    \hfill
   \begin{subfigure}[b]{0.45\textwidth}
   \centering
         \includegraphics[scale=0.4]{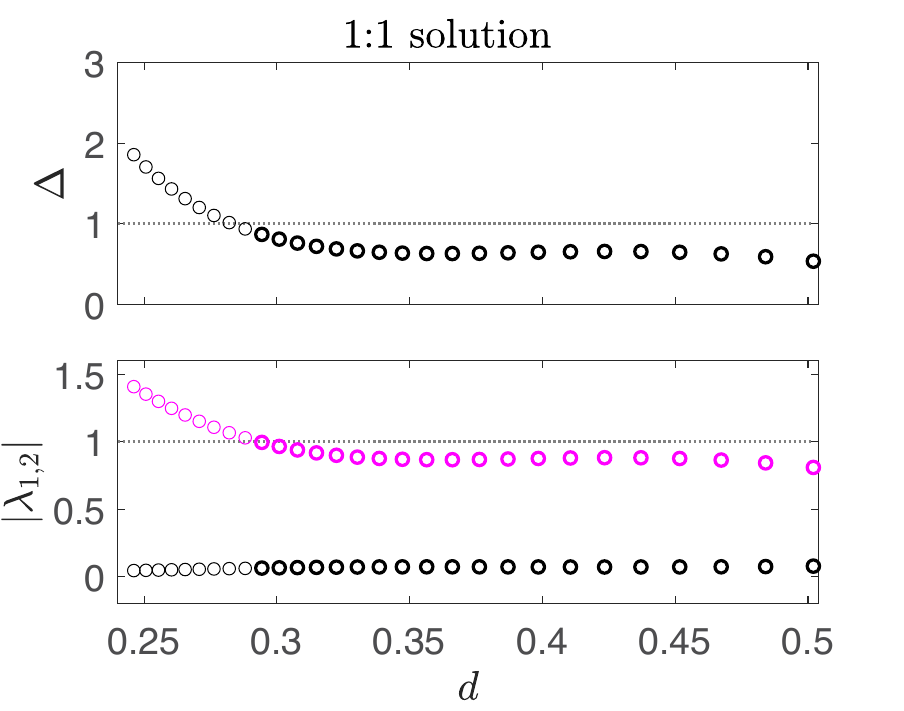}
    \caption{}
   \label{fig:SimpleSolutionBranch_Eigs}
    \end{subfigure}
    \hfill
    \begin{subfigure}[b]{0.45\textwidth}
    \centering
         \includegraphics[scale=0.4]{Figures/GrazingPaper_Figure3b_ParSet_muk=0p5_PDSimpleBranch_FullDescription2.pdf}
    \caption{}
    \label{fig:PDSimpleSolutionBranch_MS}
    \end{subfigure}
    \hfill
   \begin{subfigure}[b]{0.45\textwidth}
   \centering
         \includegraphics[scale=0.4]{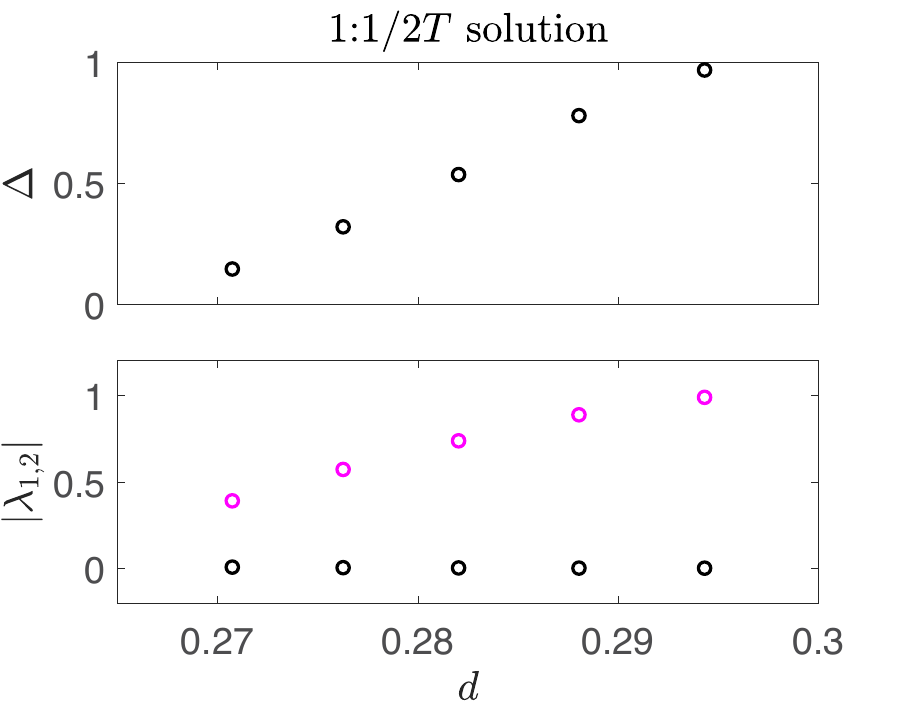}
    \caption{}
   \label{fig:PDSimpleSolutionBranch_Eigs}
    \end{subfigure}
      \hfill
    \begin{subfigure}[b]{0.45\textwidth}
    \centering
         \includegraphics[scale=0.4]{Figures/GrazingPaper_Figure3b_ParSet_muk=0p5_FullDescriptionOfLoopingBranch2.pdf}
    \caption{}
    \label{fig:LoopingSolutionBranch_MS}
    \end{subfigure}
    \hfill
    \begin{subfigure}[b]{0.45\textwidth}
    \centering
         \includegraphics[scale=0.4]{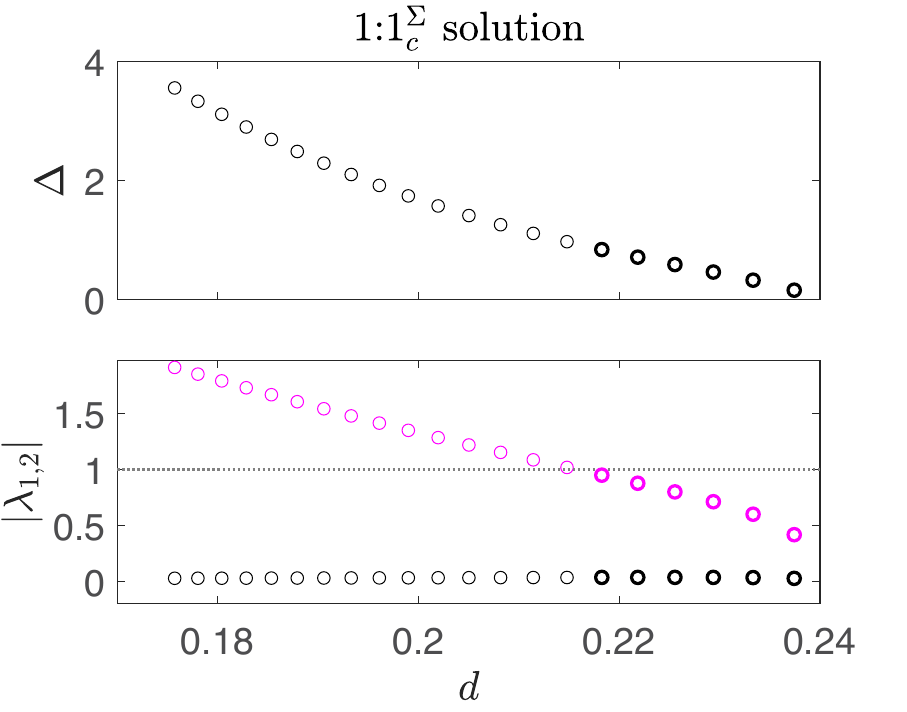}
    \caption{}
    \label{fig:LoopingSolutionBranch_Eigs}
    \end{subfigure}
    \caption{Analytical results for stable, unstable, and unphysical (a) 1:1 solutions by solving system \eqref{eq:RSE_SimpleSol}, (c) 1:1/$2T$ solutions as discussed in \cref{rem:AnalyticsFor2TSolutions}, and (e) 1:1$^{\Sigma}_{c}$ solutions by solving system \eqref{eq:RSELoopingSolution} plotted over the numerical results (blue and green dots) and using the parameters in \Cref{fig:ZdotBifDiagram_muk=0p5_AnalyticalSols_A}, i.e. $\mu_k=0.5$, $\beta = \pi/4$, $r=0.5$, $\omega =5\pi$, $s=0.5$, $A\in[3.1,14.5]$. Throughout the thick/thin circles correspond to a stable/unstable physical periodic solution (see corresponding eigenvalues in panels (b), (d), (f)). The various dashed lines, crosses (black in (a) and orange in (c)) and dashed-dotted lines correspond to unphysical solutions via a grazing bifurcation on $\Sigma$ (grazing-sliding) or $\Gamma^{+}$, since the switching or impacting dynamics of \eqref{eq:PWSsystem}-\eqref{eq:ImpactCondsRelFrame} are violated. (b),(d),(f): Discriminant (top panels) and eigenvalues (bottom panels) of the maps corresponding to the physical (b) 1:1, (d) 1:1$/2T$, and (f) 1:1$^{\Sigma}_{c}$ motions, respectively.}
\label{fig:LSA_OneOneBranch_MS}
\end{figure}

\Cref{fig:LSA_OneOneBranch_MS} shows the bifurcation diagrams of $\dot{Z}_k$ vs $d$ shown in \Cref{fig:LSA_OneOneBranch} (panels (a), (c), (e), repeated here from \Cref{fig:LSA_OneOneBranch} for the reader's convenience), and discriminants and eigenvalues (panels (b), (d), (f)) obtained from the stability analysis of the maps corresponding to physical 1:1 (panels (a)-(b)), 1:1$/2T$ (panels (c)-(d)), and 1:1$^{\Sigma}_{c}$ (panels (e)-(f)).
\newpage
\section*{Acknowledgments}
We would like to thank NSF, USA (Grant 
 CMMI/DMS 2009270) and EPSRC, UK (Grant 
 EP/V034391/1) for providing the financial support for this work.

\bibliographystyle{plain}
\bibliography{MS_References}

\begin{thebibliography}{10}

\bibitem{PrivComm_HSDY}
\url{https://sites.google.com/view/vienergyharvest}.

\bibitem{Awr}
Jan Awrejcewicz and Paweł Olejnik.
\newblock Analysis of dynamic systems with various friction laws.
\newblock {\em Applied Mechanics Reviews}, 58:389--411, 2005.

\bibitem{bapat1988impact}
CN~Bapat and C~Bapat.
\newblock Impact-pair under periodic excitation.
\newblock {\em Journal of Sound and Vibration}, 120(1):53--61, 1988.

\bibitem{CONE1995659}
K.M. Cone and R.I. Zadoks.
\newblock A numerical study of an impact oscillator with the addition of dry friction.
\newblock {\em Journal of Sound and Vibration}, 188(5):659--683, 1995.

\bibitem{DiBernardo_Review2008}
Mario di~Bernardo, Chris~J. Budd, Alan~R. Champneys, Piotr Kowalczyk, Arne~B. Nordmark, Gerard~Olivar Tost, and Petri~T. Piiroinen.
\newblock Bifurcations in nonsmooth dynamical systems.
\newblock {\em SIAM Review}, 50(4):629--701, 2008.

\bibitem{DULIN2022106983}
Sam Dulin, Kailee Lin, Larissa Serdukova, Rachel Kuske, and Daniil Yurchenko.
\newblock Improving the performance of a two-sided vibro-impact energy harvester with asymmetric restitution coefficients.
\newblock {\em International Journal of Mechanical Sciences}, 217:106983, 2022.

\bibitem{Jeffrey_Hogan2011}
M.~R. Jeffrey and S.~J. Hogan.
\newblock The geometry of generic sliding bifurcations.
\newblock {\em SIAM Review}, 53(3):505--525, 2011.

\bibitem{KUMAR2024118131}
Rahul Kumar, Rachel Kuske, and Daniil Yurchenko.
\newblock Exploring effective tet through a vibro-impact nonlinear energy sink over broad parameter regimes.
\newblock {\em Journal of Sound and Vibration}, 570:118131, 2024.

\bibitem{LI201981}
Chunliang Li, Jinjun Fan, Zhaoxia Yang, and Shan Xue.
\newblock On discontinuous dynamical behaviors of a 2-dof impact oscillator with friction and a periodically forced excitation.
\newblock {\em Mechanism and Machine Theory}, 135:81--108, 2019.

\bibitem{LI2021104001}
Haiqin Li, Ang Li, and Yunfa Zhang.
\newblock Importance of gravity and friction on the targeted energy transfer of vibro-impact nonlinear energy sink.
\newblock {\em International Journal of Impact Engineering}, 157:104001, 2021.

\bibitem{LIU201530}
Yang Liu, Ekaterina Pavlovskaia, Marian Wiercigroch, and Zhike Peng.
\newblock Forward and backward motion control of a vibro-impact capsule system.
\newblock {\em International Journal of Non-Linear Mechanics}, 70:30--46, 2015.
\newblock Nonlinear Dynamics in Engineering: Modelling, Analysis and Applications.

\bibitem{Luo2013VibroImpactDL}
Albert C.~J. Luo and Yu~Qi Guo.
\newblock Vibro-impact dynamics: Luo/vibro.
\newblock 2013.

\bibitem{LUO20081030}
Albert~C.J. Luo.
\newblock A theory for flow switchability in discontinuous dynamical systems.
\newblock {\em Nonlinear Analysis: Hybrid Systems}, 2(4):1030--1061, 2008.

\bibitem{Or12}
Yizhara Or and Elon Rimon.
\newblock Investigation of painlevé’s paradox and dynamic jamming during mechanism sliding motion.
\newblock {\em Nonlinear Dynamics}, 67:1647--1668, 2012.

\bibitem{Painleve}
Paul Painlevé.
\newblock Sur le lois frottement de glissemment.
\newblock {\em C. R. Acad. Sci.}, 121:112--115, 1895.

\bibitem{SERDUKOVA2021115811}
Larissa Serdukova, Rachel Kuske, and Daniil Yurchenko.
\newblock Post-grazing dynamics of a vibro-impacting energy generator.
\newblock {\em Journal of Sound and Vibration}, 492:115811, 2021.

\bibitem{Serdukova2019StabilityAB}
Larissa Serdukova, Rachel Kuske, and Daniil~V. Yurchenko.
\newblock Stability and bifurcation analysis of the period-t motion of a vibroimpact energy harvester.
\newblock {\em Nonlinear Dynamics}, 98:1807 -- 1819, 2019.

\bibitem{Serdukova2022FundamentalCO}
Larissa Serdukova, Rachel Kuske, and Daniil~V. Yurchenko.
\newblock Fundamental competition of smooth and non-smooth bifurcations and their ghosts in vibro-impact pairs.
\newblock {\em Nonlinear Dynamics}, 111:6129--6155, 2022.

\bibitem{SHAW1983129}
S.W. Shaw and P.J. Holmes.
\newblock A periodically forced piecewise linear oscillator.
\newblock {\em Journal of Sound and Vibration}, 90(1):129--155, 1983.

\bibitem{Hogan}
Hogan S.J. and Kristiansen K.U.
\newblock On the regularization of impact without collision: the painlevé paradox and compliance.
\newblock {\em Proc Math Phys Eng Sci.}, 473, 2017.

\bibitem{Stelter1991}
Peter Stelter and Walter Sextro.
\newblock {\em Bifurcations in Dynamic Systems with Dry Friction}, pages 343--347.
\newblock Birkh{\"a}user Basel, Basel, 1991.

\bibitem{YURCHENKO2017456}
D.~Yurchenko, Z.H. Lai, G.~Thomson, D.V. Val, and R.V. Bobryk.
\newblock Parametric study of a novel vibro-impact energy harvesting system with dielectric elastomer.
\newblock {\em Applied Energy}, 208:456--470, 2017.

\bibitem{Yurchenko_2017}
D~Yurchenko, D~V Val, Z~H Lai, G~Gu, and G~Thomson.
\newblock Energy harvesting from a de-based dynamic vibro-impact system.
\newblock {\em Smart Materials and Structures}, 26(10):105001, aug 2017.

\bibitem{ZhangFu2017}
Yanyan Zhang and Xilin Fu.
\newblock Flow switchability of motions in a horizontal impact pair with dry friction.
\newblock {\em Communications in Nonlinear Science and Numerical Simulation}, 44:89--107, 2017.

\end{thebibliography}
\end{document}